\documentclass[11pt]{article}%
\usepackage{amsfonts}
\usepackage{amsmath}
\usepackage{amssymb}
\usepackage{amsxtra}
\usepackage{graphicx}
\usepackage{geometry}
\usepackage{color}
\usepackage{colortbl}
\usepackage{caption}%
\providecommand{\U}[1]{\protect\rule{.1in}{.1in}}
\newtheorem{theorem}{Theorem}

\newtheorem{proposition}[theorem]{Proposition}
\newtheorem{remark}[theorem]{Remark}

\numberwithin{equation}{section}
\ifx\pdfoutput\relax\let\pdfoutput=\undefined\fi
\newcount\msipdfoutput
\ifx\pdfoutput\undefined\else
\ifcase\pdfoutput\else
\ifx\paperwidth\undefined\else
\ifdim\paperheight=0pt\relax\else\pdfpageheight\paperheight\fi
\ifdim\paperwidth=0pt\relax\else\pdfpagewidth\paperwidth\fi
\fi\fi\fi
\begin{document}

\title{Embedded eigenvalues of the Neumann problem in a strip\\with a box-shaped perturbation}
\author{G.Cardone\thanks{Universit\`{a} del Sannio, Department of Engineering, Corso
Garibaldi, 107, 82100 Benevento, Italy; email: giuseppe.cardone@unisannio.it},
T.Durante\thanks{University of Salerno, Department of Information Engineering,
Electrical Engineering and Applied Mathematics, Via Giovanni Paolo II, 132,
84084 Fisciano (SA), Italy; email: tdurante@unisa.it},
S.A.Nazarov\thanks{Mathematics and Mechanics Faculty, St. Petersburg State
University, 198504, Universitetsky pr., 28, Stary Peterhof, Russia;
Saint-Petersburg State Polytechnical University, Polytechnicheskaya ul., 29,
St. Petersburg, 195251, Russia; Institute of Problems of Mechanical
Engineering RAS, V.O., Bolshoj pr., 61, St. Petersburg, 199178, Russia; email:
srgnazarov@yahoo.co.uk.}}
\maketitle

\begin{abstract}
\medskip We consider the spectral Neumann problem for the Laplace operator in
an acoustic waveguide $\Pi_{l}^{\varepsilon}$ obtained from a straight unit
strip by a low box-shaped perturbation of size $2l\times\varepsilon,$ where
$\varepsilon>0$ is a small parameter. We prove the existence of the length
parameter $l_{k}^{\varepsilon}=\pi k+O(\varepsilon)$ with any $k=1,2,3,...$
such that the waveguide $\Pi_{l_{k}^{\varepsilon}}^{\varepsilon}$ supports a
trapped mode with an eigenvalue $\lambda_{k}^{\varepsilon}=\pi^{2}-4\pi
^{4}l^{2}\varepsilon^{2}+O(\varepsilon^{3})$ embedded into the continuous
spectrum. This eigenvalue is unique in the segment $[0,\pi^{2}]$ and is absent
in the case $l\neq l_{k}^{\varepsilon}.$ The detection of this embedded
eigenvalue is based on a criterion for trapped modes involving an artificial
object, the augmented scattering matrix. The main technical difficulty is
caused by corner points of the perturbed wall $\partial\Pi_{l}^{\varepsilon}$
and we discuss available generalizations for other piecewise smooth boundaries.

\medskip

Keywords: acoustic waveguide, Neumann problem, embedded eigenvalues,
continuous spectrum, box-shaped perturbation, asymptotics

\medskip

MSC: 35P05, 47A75, 49R50, 78A50.

\end{abstract}

\section{Introduction\label{sect1}}

\subsection{Formulation of problems\label{sect1.1}}

\begin{figure}
[ptb]
\begin{center}
\ifcase\msipdfoutput
\includegraphics[
height=0.9634in,
width=3.992in
]%
{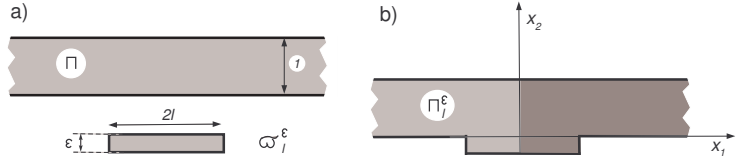}%
\else
\includegraphics[
height=0.9634in,
width=3.992in
]%
{C:/Users/Desktop/Documents/Google Drive/Articoli/Articoli-in-Preparazione/CardoneDuranteNazarov1/graphics/NewCaDuNa__1.pdf}%
\fi
\caption{The waveguide with a box-shaped perturbation (a) and its fragments
(b)}%
\label{f1}%
\end{center}
\end{figure}

In the union $\Pi_{l}^{\varepsilon}$, fig. \ref{f1}, b and a, of the unit
straight strip%
\begin{equation}
\Pi=\left\{  x=\left(  x_{1},x_{2}\right)  \in\mathbb{R}^{2},\text{ }x_{1}%
\in\mathbb{R},\text{ }x_{2}\in\left(  0,1\right)  \right\}  \label{00}%
\end{equation}
and a rectangle of length $2l>0$ and a small width $\varepsilon>0,$%
\begin{equation}
\varpi_{l}^{\varepsilon}=\left\{  x:\left\vert x_{1}\right\vert <l,\text{
\ }x_{2}\in\left(  -\varepsilon,0\right]  \right\}  , \label{0}%
\end{equation}
we consider the spectral Neumann problem%
\begin{align}
-\Delta u^{\varepsilon}\left(  x\right)   &  =\lambda^{\varepsilon
}u^{\varepsilon}\left(  x\right)  ,\text{ \ }x\in\Pi_{l}^{\varepsilon}=\Pi
\cup\varpi_{l}^{\varepsilon},\label{1}\\
\partial_{\nu}u^{\varepsilon}\left(  x\right)   &  =0,\ x\in\partial\Pi
_{l}^{\varepsilon}, \label{2}%
\end{align}
%

where $\Delta=\nabla\cdot\nabla$ is the Laplace operator, $\nabla
=\operatorname{grad},$ $\lambda^{\varepsilon}$ is the spectral parameter and
$\partial_{\nu}=\nu\cdot\nabla$ is the directional derivative, $\nu$ stands
for the unit outward normal defined everywhere at the boundary $\partial
\Pi_{l}^{\varepsilon},$ except for corner points, i.e. vertices of the
rectangle (\ref{0}). Since a solution of the problem (\ref{1}), (\ref{2}) may
get singularities at these points, the problem ought to be reformulated as the
integral identity \cite{Lad}%
\begin{equation}
\left(  \nabla u^{\varepsilon},\nabla v^{\varepsilon}\right)  _{\Pi
_{l}^{\varepsilon}}=\lambda^{\varepsilon}\left(  u^{\varepsilon}%
,v^{\varepsilon}\right)  _{\Pi_{l}^{\varepsilon}},\ \forall v^{\varepsilon}\in
H^{1}\left(  \Pi_{l}^{\varepsilon}\right)  , \label{3}%
\end{equation}
where $\left(  \ ,\ \right)  _{\Pi_{l}^{\varepsilon}}$ is the natural scalar
product in the Lebesgue space $L^{2}\left(  \Pi_{l}^{\varepsilon}\right)  $
and $H^{1}\left(  \Pi_{l}^{\varepsilon}\right)  $ stands for the Sobolev
space. The symmetric bilinear form on the left-hand side of (\ref{3}) is
closed and positive in $H^{1}\left(  \Pi_{l}^{\varepsilon}\right)  $ so that
problem (\ref{1}), (\ref{2}) is associated \cite[Ch 10]{BiSo} with a positive
self-adjoint operator $A_{l}^{\varepsilon}$ in $L^{2}\left(  \Pi
_{l}^{\varepsilon}\right)  $ whose spectrum $\wp=\wp_{co}$ is continuous and
covers the closed positive semi-axis $\overline{\mathbb{R}}_{+}=\left[
0,+\infty\right)  .$ The domain $\mathcal{D}\left(  A_{l}^{\varepsilon
}\right)  $ of $A_{l}^{\varepsilon},$ of course, belongs to $H^{1}\left(
\Pi_{l}^{\varepsilon}\right)  $ but is bigger than $H^{2}\left(  \Pi
_{l}^{\varepsilon}\right)  $ due to singularities of solutions at the corner
points, see, e.g., \cite[Ch.2]{NaPl}. The point spectrum $\wp_{po}$ of
$A_{l}^{\varepsilon}$ can be non-empty and the main goal of our paper is to
single out a particular value of the length parameter $l$ such that the
operator $A_{l}^{\varepsilon}$ wins an eigenvalue $\lambda_{l}^{\varepsilon
}\in\wp_{po}$ embedded into the continuous spectrum. The corresponding
eigenfunction $u_{l}^{\varepsilon}\in H^{1}\left(  \Pi_{l}^{\varepsilon
}\right)  $ decays exponentially at infinity and is called a trapped mode, cf.
\cite{LM} and \cite{Urcell}.

Our central result formulated below in Theorem \ref{TheoremEX}, roughly
speaking, demonstrates that an eigenvalue $\lambda_{l}^{\varepsilon}$ exists
in the interval $\left(  0,\pi^{2}\right)  \subset\wp_{co}$ for
$l^{\varepsilon}\approx\pi k$ with $k\in\mathbb{N}=\left\{  1,2,3,...\right\}
$ only. Asymptotics of $\lambda^{\varepsilon}$ and $l^{\varepsilon}$ are
constructed, too.

Problem (\ref{1}), (\ref{2}) is a model of an acoustic waveguide with hard
walls, cf. \cite{Acous}, but is also related in a natural way to the linear
theory of surface water-waves, cf. \cite{KuMaVa}. Indeed, the velocity
potential $\Phi^{\varepsilon}\left(  x,z\right)  $ satisfies the Laplace
equation in the channel $\Xi_{l,d}^{\varepsilon}=\Pi_{l}^{\varepsilon}%
\times\left(  -d,0\right)  \subset\mathbb{R}^{3}\ni\left(  x,z\right)  $ of
depth $d>0$ with the Neumann condition (no normal flow) at its vertical walls
and horizontal bottom as well as the spectral Steklov condition (the kinetic
one) on the free horizontal surface%
\[
\partial_{z}\Phi^{\varepsilon}\left(  x,0\right)  =\Lambda^{\varepsilon}%
\Phi^{\varepsilon}\left(  x,0\right)  ,\text{ \ }x\in\Pi_{l}^{\varepsilon}.
\]
After factoring out the dependence on the vertical variable $z$,
\begin{equation}
\Phi^{\varepsilon}\left(  x,z\right)  =u^{\varepsilon}\left(  x\right)
\left(  e^{z\lambda^{\varepsilon}}+e^{-\left(  z+2d\right)  \lambda
^{\varepsilon}}\right)  , \label{4}%
\end{equation}
see, e.g., \cite{Vas1} and \cite{LM}, the water-wave problem reduces to the
two-dimensional Neumann problem (\ref{1}), (\ref{2}) for the function
$u^{\varepsilon}$ in (\ref{4}) and the parameter $\lambda^{\varepsilon}$
determined from the equation%
\[
\Lambda^{\varepsilon}=\lambda^{\varepsilon}\frac{1-e^{-2d\lambda^{\varepsilon
}}}{1+e^{-2d\lambda^{\varepsilon}}}=\lambda^{\varepsilon}\tanh\left(
d\lambda^{\varepsilon}\right)  .
\]
We will not discuss separately this interpretation of our problem but in the
next section present some asymptotic formulas for eigenvalues of the Laplace
operator with either Dirichlet, or mixed boundary conditions.

\subsection{Asymptotics of eigenvalues\label{sect1.2}}

\begin{figure}
[ptb]
\begin{center}
\ifcase\msipdfoutput
\includegraphics[
height=0.7014in,
width=4.3716in
]%
{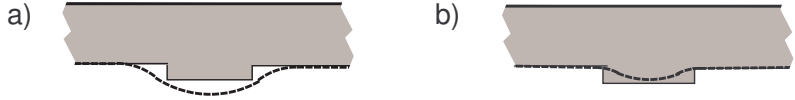}%
\else
\includegraphics[
height=0.7014in,
width=4.3716in
]%
{C:/Users/Desktop/Documents/Google Drive/Articoli/Articoli-in-Preparazione/CardoneDuranteNazarov1/graphics/CaDuNa2__2.pdf}%
\fi
\caption{The box-shaped perturbation enveloping (a) and entering (b) a regular
perturbation.}%
\label{f2}%
\end{center}
\end{figure}

\begin{figure}
[ptb]
\begin{center}
\ifcase\msipdfoutput
\includegraphics[
height=1.2496in,
width=3.5777in
]%
{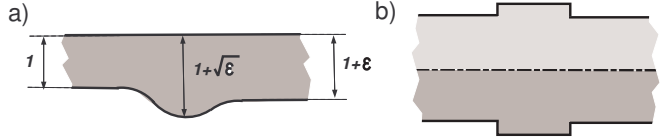}%
\else
\includegraphics[
height=1.2496in,
width=3.5777in
]%
{C:/Users/Desktop/Documents/Google Drive/Articoli/Articoli-in-Preparazione/CardoneDuranteNazarov1/graphics/CaDuNa3__3.pdf}%
\fi
\caption{A non-local perturbation (a) and a symmetric box-shaped perturbation
(b)}%
\label{f3}%
\end{center}
\end{figure}

Imposing the Dirichlet condition%
\begin{equation}
u^{\varepsilon}\left(  x\right)  =0,\text{ \ \ }x\in\partial\Pi_{l}%
^{\varepsilon}, \label{5}%
\end{equation}
instead of the Neumann condition (\ref{2}), creates the positive cut-off value
$\lambda_{\dagger}=\pi^{2}$ of the continuous spectrum $\wp^{D}=\left[
\pi^{2},+\infty\right)  $ of the Dirichlet problem (\ref{1}), (\ref{5}) which
provides an adequate model of a quantum waveguide, cf. \cite{Keijo}. The
interval $\left(  0,\pi^{2}\right)  $ stays now below the continuous spectrum
and therefore may contain eigenvalues composing the discrete spectrum
$\wp_{di}^{D}$ of the problem. As follows from a result in \cite{Sim}, the
multiplicity $\#\wp_{di}^{D}$ is equal to $1$ for a small $\varepsilon>0$.
Although the paper \cite{Sim} deals with a regular (smooth) perturbation of
the wall, it is possible to select two smooth shallow pockets boxes as in fig.
\ref{f2}, a and b, and to extend the existence and uniqueness result in
\cite{Sim} for the box-shaped perturbations by means of the max-min principle,
see, e.g.,\cite[Thm 10.2.2]{BiSo}. However, the attendant asymptotic formula
\begin{equation}
\lambda_{l}^{\varepsilon}=\pi^{2}-4\pi^{4}\varepsilon^{2}l^{2}+O\left(
\varepsilon^{3}\right)  ,\text{ \ }\varepsilon\rightarrow0, \label{6}%
\end{equation}
cannot be supported by these results because an application in \cite{Sim} of a
change of variables which transforms $\Pi_{l}^{\varepsilon}$ into $\Pi$
requires certain smoothness properties of the boundary $\partial\Pi
_{l}^{\varepsilon}$ which are naturally absent, fig. \ref{f1}, b. In Section
\ref{sect8.3} we will explain how our approach helps to justify formula
(\ref{6}). Local perturbations of quantum waveguides in $\mathbb{R}^{n}$,
$n\geq2$, are intensively investigated during last two decades and many
important results on the existence and asymptotic behavior of their discrete
spectrum have been published. We mention a few of them, namely \cite{ExDu, Ex,
Masl} for the slightly curved and twisted cylindrical waveguides,
\cite{sha1,na480,na561} for cranked waveguides and \cite{Gad, Gru} for the
Laplacian perturbed by a small second-order differential operator with
compactly supported coefficients. We also refer to \cite{na534} for non-local
perturbations, fig. \ref{f3}, a and to \cite{BoBuCa1, BoBuCa4, Sim} for
alternation of the Dirichlet and Neumann boundary conditions.%

In the literature one finds much less results on eigenvalues embedded into the
continuous spectrum, cf. the review papers \cite{BBD1, BBD2, LM}. First of all
we describe an elegant method \cite{Vas1} which is based on imposing an
artificial Dirichlet condition and had become of rather wide use in proving
the existence of embedded eigenvalues but only in symmetric waveguides.%

Let us consider an auxiliary mixed boundary value problem and supply the
Helmholtz equation (\ref{1}) with the Neumann condition on the lower lagged
wall and the Dirichlet condition on the upper straight wall, see fig.
\ref{f1}, b,%

\begin{equation}
u^{\varepsilon}\left(  x_{1},1\right)  =0,\text{ \ \ }x_{1}\in\mathbb{R}%
,\ \ \ \ \ \ \partial_{\nu}u^{\varepsilon}\left(  x\right)  =0,\text{
\ \ }x\in\partial\Pi_{l}^{\varepsilon},\text{ }x_{2}<1. \label{7}%
\end{equation}
Problem (\ref{1}), (\ref{7}), has the continuous spectrum $\wp_{co}%
^{M}=\left[  \pi^{2}/4,+\infty\right)  $ and in Section \ref{sect8.1} we will
show the existence of only one eigenvalue%
\begin{equation}
\lambda_{l}^{\varepsilon}=\frac{\pi^{2}}{4}(1-\pi^{2}l^{2}\varepsilon
^{2})+O(\varepsilon^{3}\left(  1+\left\vert \ln\varepsilon\right\vert \right)
^{2}),\text{ \ }\varepsilon\rightarrow+0, \label{8}%
\end{equation}
in the discrete spectrum $\wp_{di}^{M}\subset\left(  0,\pi^{2}/4\right)  $.
Following \cite{Vas1} we extend the corresponding eigenfunction $u_{l}%
^{\varepsilon}\left(  x_{1},x_{2}\right)  $ as an odd function in $x_{2}-1$
from $\Pi_{l}^{\varepsilon}$ onto the bigger waveguide $\widehat{\Pi}%
_{l}^{\varepsilon}$ drawn in fig. \ref{f3}, b and obtained as the union of the
strip $\mathbb{R}\times\left(  0,2\right)  $ and the box $\left(  -l,l\right)
\times\left(  -\varepsilon,2+\varepsilon\right)  $. Owing to the Dirichlet
condition in (\ref{7}) at the midline of $\widehat{\Pi}_{l}^{\varepsilon}$,
this extension $\widehat{u}_{l}^{\varepsilon}\left(  x_{1},x_{2}\right)  $ is
a smooth function everywhere in $\widehat{\Pi}_{l}^{\varepsilon},$ except at
corner points and inherits from $u_{l}^{\varepsilon}\left(  x_{1}%
,x_{2}\right)  $ an exponential decay as $x_{1}\rightarrow\pm\infty$. Clearly%
\begin{equation}
-\Delta\widehat{u}_{l}^{\varepsilon}\left(  x\right)  =\lambda_{l}%
^{\varepsilon}\widehat{u}_{l}^{\varepsilon}\left(  x\right)  ,\text{ \ }%
x\in\widehat{\Pi}_{l}^{\varepsilon},\text{ \ }\partial_{\nu}\widehat{u}%
_{l}^{\varepsilon}\left(  x\right)  =0,\text{ \ \ \ }x\in\partial\widehat{\Pi
}_{l}^{\varepsilon}, \label{888}%
\end{equation}
and, thus, $\widehat{u}_{l}^{\varepsilon}$ is an eigenfunction of the Neumann
problem (\ref{888}) while the corresponding eigenvalue (\ref{8}) belongs to
the continuous spectrum $\widehat{\wp}_{co}=\left[  0,+\infty\right)  $ of
this problem.

We emphasize that the method \cite{Vas1} requires the mirror symmetry of the
waveguide and cannot be applied to the asymmetric waveguide $\Pi
_{l}^{\varepsilon}$ in fig. \ref{f1}, b. The detected embedded eigenvalue
$\lambda_{l}^{\varepsilon}$ of the Neumann problem (\ref{888}) is stable with
respect to small symmetric perturbations of the waveguide walls but any
violation of the symmetry may lead it out from the spectrum and turn it into a
point of complex resonance, cf. \cite{Vas4} and, e.g., \cite{na546}.

The intrinsic instability of embedded eigenvalues requests for special
techniques to detect them as well as to construct their asymptotics. In the
present paper we use a criterion for the existence of trapped modes (see
\cite{na275} and Theorem \ref{TheoremA} below) and a concept of enforced
stability of eigenvalues in the continuous spectrum, cf. \cite{na489, na546}.

\subsection{Reduction of the problem\label{sect1.3}}

In view of the mirror symmetry about the $x_{2}$-axis, notice the difference
with the above mentioned assumption in \cite{Vas1}, we truncate the waveguide
$\Pi_{l}^{\varepsilon}$ and consider the Neumann problem%
\begin{align}
-\Delta u_{+}^{\varepsilon}\left(  x\right)   &  =\lambda_{+}^{\varepsilon
}u_{+}^{\varepsilon}\left(  x\right)  ,\text{ \ }x\in\Pi_{l+}^{\varepsilon
},\label{9}\\
\partial_{\nu}u_{+}^{\varepsilon}\left(  x\right)   &  =0,\text{ \ \ \ }%
x\in\partial\Pi_{l+}^{\varepsilon}, \label{10}%
\end{align}
is its right half (overshaded in fig. \ref{f1}, b)%
\begin{equation}
\Pi_{l+}^{\varepsilon}=\left\{  x\in\Pi_{l}^{\varepsilon}:x_{1}>0\right\}
=\left\{  x:x_{2}\in(-\varepsilon,0)\text{ for }x_{1}\in\left(  0,l\right)
,\text{ }x_{2}\in(0,1)\text{ for }x_{1}\geq l\right\}  . \label{11}%
\end{equation}
Clearly, the even in $x_{1}$ extension of an eigenfunction $u_{+}%
^{\varepsilon}$ of problem (\ref{9}), (\ref{10}) becomes an eigenfunction of
the original problem (\ref{1}), (\ref{2}). Searching for an eigenvalue%
\begin{equation}
\lambda^{\varepsilon}\in\left(  0,\pi^{2}\right)  , \label{12}%
\end{equation}
we will show in Section \ref{sect7} that, first, problem (\ref{9}), (\ref{10})
cannot get more than one eigenvalue in $\left(  0,\pi^{2}\right)  $ and,
second, the mixed boundary value problem in (\ref{11}) with the Dirichlet
condition at the truncation segment $\left\{  x:x_{1}=0,\text{ \ }x_{2}%
\in(-\varepsilon,1)\right\}  $ instead of the Neumann condition as in
(\ref{10}), does not have eigenvalues (\ref{12}). These mean that an
eigenfunction of problem (\ref{1}), (\ref{2}) associated with the eigenvalue
(\ref{12}) is always even in the variable $x_{1}$. In this way, we will be
able to describe the part $\wp_{po}\cap\left(  0,\pi^{2}\right)  $ of the
point spectrum in the entire waveguide $\Pi_{l}^{\varepsilon}$. In what
follows we skip the subscript $l$. Hence, we regard (\ref{3}) as an integral
identity serving for problem (\ref{9}), (\ref{10}) in $\Pi_{+}^{\varepsilon
}:=\Pi_{l+}^{\varepsilon}$ and denote $A_{+}^{\varepsilon}$ the corresponding
self-adjoint operator in $L^{2}\left(  \Pi_{+}^{\varepsilon}\right)  $, cf.
Section \ref{sect1.1}.

\subsection{Embedded eigenvalues\label{sect1.4}}

In the absence of the mirror symmetry about a midline of a waveguide the
modern literature gives much less results on the existence of embedded
eigenvalues. A distinguishing feature of an eigenvalue in the continuous
spectrum is its intrinsic instability with respect to a variation of the
waveguide shape while all eigenvalues in the continuous spectrum stay stable.
In this way to detect an eigenvalue in the Neumann waveguide $\Pi
_{l}^{\varepsilon}$, see (\ref{1}), (\ref{2}), a fine tuning of the
parameter-dependent shape is needed, namely the length $2l=2l(\varepsilon)$ of
the perturbation box (\ref{0}) must be chosen specifically in dependence of
its height $\varepsilon$.

\begin{figure}
[ptb]
\begin{center}
\ifcase\msipdfoutput
\includegraphics[
height=0.6512in,
width=4.235in
]%
{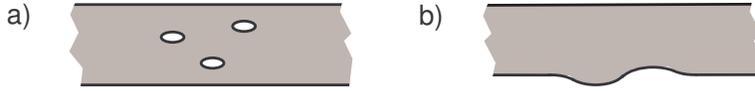}%
\else
\includegraphics[
height=0.6512in,
width=4.235in
]%
{C:/Users/Desktop/Documents/Google Drive/Articoli/Articoli-in-Preparazione/CardoneDuranteNazarov1/graphics/CaDuNa5__4.pdf}%
\fi
\caption{Singular (a) and regular (b) perturbations}%
\label{f5}%
\end{center}
\end{figure}

A method to construct particular waveguide shapes which support embedded
eigenvalues, was developed in \cite{na246, na252, na476, na489, na546} on the
basis of a sufficient condition \cite{na275} for the existence of
exponentially decaying solutions trapped modes in waveguides to elliptic
problems in domains with cylindrical outlets to infinity. As a result, several
examples of eigenvalues in the continuous spectrum were constructed without
requiring a geometrical symmetry of waveguides which are obtained from the
straight unit strip by either singular \cite{na246, na252, na533}, or regular
\cite{na476, na489, na546} perturbations of the boundary, see fig. \ref{f5}, a
and b, respectively; for a non-local smooth perturbation, see \cite{na521}. To
this end, a notion of the augmented scattering matrix \cite{na275} was used
together with certain traditional asymptotic procedures in domains with small
holes and cavities, cf. \cite[Ch.4, 5 and 2]{MaNaPl}, or in domains with
smoothly varied boundaries, cf. \cite[Ch. XII, \S 6.5]{Kato}.%

It is very important in the above-mentioned approach to detect embedded
eigenvalues that the asymptotic procedures in use are completed in such a way
they provide both, a formal derivation of asymptotic expansions and an
operator reformulation of the diffraction problem in the framework of the
perturbation theory of linear operators, see, e.g., \cite{HilPhi, Kato}, which
helps to conclude a smooth (actually analytic) dependence of scattering
matrices on geometrical parameters describing the waveguide shape. Moreover,
this permits to reformulate the sufficient condition \cite{na275} for the
existence of trapped modes as a nonlinear abstract equation and to fulfil the
condition by means of the contraction principle, cf. Sections \ref{sect2.3}
and \ref{sect4.1}, \ref{sect4.2} below.

The box-shaped perturbation (\ref{0}) of the strip (\ref{00}) can be regarded
as a combination of regular and singular perturbations, respectively outside
and inside neighborhoods of the corner points but unfortunately the authors do
not know a tool to reduce the problem (\ref{1}), (\ref{2}) (or (\ref{3}) in
the variational form) to an abstract equation in\ a fixed (independent of
$\varepsilon$) Banach space and to confirm necessary properties of the
scattering matrices in the waveguide $\Pi_{l}^{\varepsilon}=\Pi\cup\varpi
_{l}^{\varepsilon}$. Thus, in order to support each step of our procedure to
detect an embedded eigenvalue and to establish its uniqueness, we have to
obtain a certain new result for the problems (\ref{1}), (\ref{2}) and/or
(\ref{9}), (\ref{10}) which cannot be deduced from the general perturbation
theory. Although several approaches and tricks proposed in our paper work also
for other shapes like in fig. \ref{f6}, we focus the analysis on the
particular shape in fig. \ref{f1}, b.

\subsection{Architecture of the paper\label{sect.1.6}}

We proceed in Section \ref{sect2} with introducing different waves in $\Pi
_{+}^{\varepsilon}$: oscillatory and exponential for $\lambda^{\varepsilon}%
\in\left(  0,\pi^{2}\right)  $\ and linear in $x_{1}$ at the threshold
$\lambda^{\varepsilon}=\pi^{2}$. Then on the basis of the Mandelstam energy
principle, cf. \cite[\S \ 5.3]{NaPl} we perform the classification $\left\{
\text{incoming/outgoing}\right\}  $ for the introduced waves and impose two
types, physical and artificial, of radiation conditions at infinity. The
corresponding diffraction problems give rise to two scattering matrices
$s^{\varepsilon}$ and $S^{\varepsilon}$. Due to the restriction of the
boundary value problem (\ref{1}), (\ref{2}) onto the semi-infinite waveguide
(\ref{11}) the matrix $s^{\varepsilon}$ reduces to classical scalar reflexion
coefficient but the augmented scattering matrix $S^{\varepsilon}$ is of size
$2\times2$ because the artificial radiation conditions involve the exponential
waves in addition to the oscillatory waves. The above-mentioned criterion for
the existence of embedded eigenvalues is formulated in terms of the matrix
$S^{\varepsilon}$, see Theorem \ref{TheoremA} and note that its proof is
completed by Theorem \ref{TheoremMA} about solutions of the problem (\ref{9}),
(\ref{10}) with a fast exponential decay. In Section \ref{sect3} we construct
formal asymptotic expansions of the augmented scattering matrix which are
justified in\ Section \ref{sect6.4}. In order to detect an embedded eigenvalue
in Section \ref{sect4} we need the main asymptotic and first correction terms.
The two-fold nature of the box-shaped perturbation manifests itself in
different ans\"{a}tze for the diagonal entries $S_{11}^{\varepsilon}$ and
$S_{00}^{\varepsilon}$ of the matrix $S^{\varepsilon}$. In the first case the
asymptotic procedure looks like as for a regular perturbation of the boundary,
that is the boundary layer phenomenon does not influence the main asymptotic
term $S_{11}^{\varepsilon}$ in the expansion%
\begin{equation}
S_{11}^{\varepsilon}=S_{11}^{0}+\widehat{S}_{11}^{\varepsilon}. \label{S11}%
\end{equation}
In the second case the correction term $S_{\infty}^{^{\prime}}$ in the
expansion $S_{00}^{\varepsilon}=1+\varepsilon S_{00}^{^{\prime}}+\widehat
{S}_{00}^{\varepsilon}$ results from the boundary layer phenomenon while the
regular expansion affects higher-order terms only. It should be emphasized
that the augmented scattering matrix is unitary and symmetric and the main
asymptotic term in the expansion%
\begin{equation}
S_{01}^{\varepsilon}=S_{10}^{\varepsilon}=\varepsilon^{1/2}S_{10}^{^{\prime}%
}+\widehat{S}_{10}^{\varepsilon} \label{S01}%
\end{equation}
of the anti-diagonal entries can be computed by means of both the approaches.

In Section \ref{sect4} we first reduce the criterion $S_{11}^{\varepsilon}=-1$
from Theorem \ref{TheoremA} to an abstract equation and second solve it with
the help of the contraction principle. Finally we formulate Theorems
\ref{TheoremEX} and \ref{TheoremUN} on the existence and uniqueness of the
embedded eigenvalue. These assertions are proved in the next three sections.
In Section \ref{sect5} formulations of the problem (\ref{1}), (\ref{2}) in the
Kondratiev \ spaces (Therem \ref{TheoremKA}) and weighted spaces with detached
asymptotics (Theorem \ref{TheoremKMM}) are presented as well as the operator
formulation of the radiation condition at infinity. At the same time, key
results for the particular box-shaped perturbation (\ref{0}), are displayed in
Sections \ref{sect5.3} and \ref{sect5.5} where we verify the absence of
trapped modes with a fast decay and clarify the dependence of majorants in a
priori estimates for solutions on the small and spectral parameters
$\varepsilon\in\left(  0,\varepsilon_{0}\right)  $ and $\lambda\in\left(
0,\pi^{2}\right)  $.

In Section \ref{sect6} we evaluate remainders in the asymptotic formulas
(\ref{S11})-(\ref{S01}) for the augmented scattering matrix Theorem
\ref{TheoremASY} while the boundary layer phenomenon brings additional powers
of $\left\vert \ln\varepsilon\right\vert $ into bounds of estimates. Other
necessary properties of the matrix are clarified in Section \ref{sect7} where
the uniqueness of the embedded eigenvalue is verified, too. Again kinks of the
perturbation profile seriously complicate the analysis.

Conclusive remarks are collected in Section \ref{sect8} where, in particular,
we discuss essential simplifications of the analysis within the discrete
spectrum and a certain hardship for detecting eigenvalues near higher
thresholds of the continuous spectrum.

\section{Augmented scattering matrix and a criterion for trapped
modes\label{sect2}}

\subsection{Classification of waves.\label{sect2.1}}

For the spectral problem (\ref{12}), the limit ($\varepsilon=0$) problem
(\ref{1}), (\ref{2}) in the straight strip $\Pi=\mathbb{R}\times\left(
0,1\right)  $ has two oscillating waves%
\begin{equation}
w_{0}^{\varepsilon\pm}\left(  x\right)  =\left(  2k\right)  ^{-1/2}e^{\pm
ik^{\varepsilon}x_{1}},\text{ \ }k^{\varepsilon}=\sqrt{\lambda^{\varepsilon}}.
\label{13}%
\end{equation}
By the Sommerfeld principle, see, e.g., \cite{Keijo, Acous, Wilcox}, the waves
$w_{0}^{\varepsilon+}$ and $w_{0}^{\varepsilon-}$ are outgoing and incoming,
respective, in the one-sided waveguide (\ref{11}) according to their
wavenumbers $+k^{\varepsilon}$ and $-k^{\varepsilon}.$ In this way, the
inhomogeneous problem (\ref{9}), (\ref{10})%
\begin{equation}
-\Delta u^{\varepsilon}\left(  x\right)  -\lambda^{\varepsilon}u^{\varepsilon
}\left(  x\right)  =f^{\varepsilon}\left(  x\right)  ,\text{ \ }x\in\Pi
_{+}^{\varepsilon},\ \ \ \ \ \ \partial_{\nu}u^{\varepsilon}\left(  x\right)
=0,\text{\ \ }x\in\partial\Pi_{+}^{\varepsilon}, \label{14}%
\end{equation}
with, for example, compactly supported right-hand side $f^{\varepsilon}$ is
supplied with the radiation condition%
\begin{equation}
u^{\varepsilon}\left(  x\right)  =c_{0}^{\varepsilon}w_{0}^{\varepsilon
+}\left(  x\right)  +\widetilde{u}^{\varepsilon}\left(  x\right)  ,\text{
\ \ }\widetilde{u}^{\varepsilon}\left(  x\right)  =O(e^{-x_{1}\sqrt{\pi
^{2}-\lambda^{\varepsilon}}}). \label{15}%
\end{equation}
In Section \ref{sect5} we will give an operator formulation of problem
(\ref{14}), (\ref{15}).

At the threshold%
\begin{equation}
\lambda=\pi^{2} \label{16}%
\end{equation}
in addition to the oscillating waves%
\begin{equation}
w_{0}^{0\pm}\left(  x\right)  =\left(  2\pi\right)  ^{-1/2}e^{\pm i\pi x_{1}}
\label{17}%
\end{equation}
there appear the staying and growing waves
\[
w_{1}^{0}\left(  x\right)  =\cos\left(  \pi x_{2}\right)  ,\ w_{1}^{1}\left(
x\right)  =x_{1}\cos\left(  \pi x_{2}\right)
\]
which cannot be classified by the Sommerfield principle because of zero
wavenumber. However, as was observed in \cite[\S \ 5.3]{NaPl}, that, first,
waves (\ref{13}) verify the relations%
\begin{equation}
q_{R}\left(  w_{p}^{\varepsilon\pm},w_{q}^{\varepsilon\pm}\right)  =\pm
i\delta_{p,q},\text{ \ \ \ }q_{R}\left(  w_{p}^{\varepsilon\pm},w_{q}%
^{\varepsilon\mp}\right)  =0 \label{19}%
\end{equation}
and, second, the linear combinations%
\begin{equation}
w_{1}^{0\pm}\left(  x\right)  =\left(  x_{1}\mp i\right)  \cos\left(  \pi
x_{2}\right)  \label{20}%
\end{equation}
together with waves (\ref{17}) fulfil formulas (\ref{19}) at $\varepsilon=0$
as well. Here, $\delta_{p,q}$ is the Kronecker symbol, $p,q=0$ in the first
case and $p,q=0,1$ in the second case, and $q_{R}$ is a symplectic, that is
sesquilinear and anti-Hermitian form%
\[
q_{R}\left(  w,v\right)  =\int_{0}^{1}\left(  \overline{v\left(
R,x_{2}\right)  }\frac{\partial w}{\partial x_{1}}\left(  R,x_{2}\right)
-w\left(  R,x_{2}\right)  \overline{\frac{\partial v}{\partial x_{1}}\left(
R,x_{2}\right)  }\right)  dx_{2}.
\]
This form emerges from the Green formula in the truncated waveguide $\Pi
_{+}^{\varepsilon}\left(  R\right)  =\{x\in\Pi_{+}^{\varepsilon},$ $x_{1}%
\in\left(  0,R\right)  \}$ and therefore does not depend on the length
parameter $R>l$ for any of introduced waves and their linear combinations.
Hence, we skip the subscript $R$ in (\ref{19}) and (\ref{20}).

For waves (\ref{19}), the sign of $\operatorname{Im}q\left(  w_{0}%
^{\varepsilon\pm},w_{0}^{\varepsilon\pm}\right)  $ coincides with the sign of
the wavenumber and therefore indicates the propagation direction. Analogously,
we call the wave $w_{1}^{0+}$ outgoing and the wave $w_{1}^{0-}$ incoming in
the waveguide $\Pi_{+}^{\varepsilon}$ so that the problem (\ref{14}) with
$\lambda^{\varepsilon}=\pi^{2}$ ought to be supplied with the following
threshold radiation condition%
\begin{equation}
u^{\varepsilon}\left(  x\right)  =c_{0}^{\varepsilon}w_{0}^{0+}\left(
x\right)  +c_{1}^{\varepsilon}w_{1}^{0+}\left(  x\right)  +\widetilde
{u}^{\varepsilon}\left(  x\right)  ,\text{ \ \ }\widetilde{u}\left(  x\right)
=O(e^{_{-\sqrt{3}\pi x_{1}}}) \label{23}%
\end{equation}
In Section \ref{sect5} we will prove that this formulation of the problem at
the threshold (\ref{16}) provides an isomorphism in its operator setting.

As was demonstrated in \cite[\S \ 5.6]{NaPl}, the form $q$ is closely related
to the Umov-Poyting vector \cite{Poynt, Umov} so that both radiation
conditions (\ref{15}) and (\ref{23}) arise from the Mandelstam (energy)
principle, see \cite[\S 5.3]{NaPl}.

\subsection{Scattering matrices and exponential waves\label{sect2.2}}

In the case (\ref{12}) the incoming wave in (\ref{13}) generates the following
solution of problem (\ref{9}), (\ref{10}):%
\begin{equation}
\zeta_{0}^{\varepsilon}\left(  x\right)  =w_{0}^{\varepsilon-}\left(
x\right)  +s_{00}^{\varepsilon}w_{0}^{\varepsilon+}\left(  x\right)
+\widetilde{\zeta}_{0}^{\varepsilon}\left(  x\right)  \label{24}%
\end{equation}
Here, the remainder $\widetilde{\zeta}_{0}^{\varepsilon}$ decays as
$O(e^{-x_{1}\sqrt{\pi^{2}-\lambda^{\varepsilon}}})$ and $s_{00}^{\varepsilon}$
is the reflection coefficient which satisfies $\left\vert s_{00}^{\varepsilon
}\right\vert =1$ due to conservation of energy.

In the same way, in the case (\ref{16}) we can determine the solutions
\[
\zeta_{p}^{\varepsilon}\left(  x\right)  =w_{p}^{0-}\left(  x\right)
+s_{0p}^{0}w_{0}^{0+}\left(  x\right)  +s_{1p}^{0}w_{1}^{0+}\left(  x\right)
+\widetilde{\zeta}_{p}^{0}\left(  x\right)
\]
where $p=1$, $\widetilde{\zeta}_{p}^{0}\left(  x\right)  =O(e^{-x_{1}\sqrt
{3}\pi})$ and the coefficients $s_{qp}^{0}$ form the (threshold) scattering
matrix $s^{0}$ of size $2\times2.$ According to the normalization and
orthogonality conditions (\ref{19}) for waves (\ref{17}), (\ref{20}) and the
relation $w_{p}^{0+}\left(  x\right)  =\overline{w_{p}^{0-}\left(  x\right)
}$, the matrix $s^{0}$ is unitary and symmetric (cf. \cite[\S 2]{na489}) that
is
\begin{equation}
\left(  s^{0}\right)  ^{-1}=\left(  s^{0}\right)  ^{\ast},\text{ \ \ }%
s^{0}=\left(  s^{0}\right)  ^{\top} \label{26}%
\end{equation}
where $\top$ stands for transposition and $\left(  s^{0}\right)  ^{\ast
}=\left(  \overline{s^{0}}\right)  ^{\top}$ is the adjoint matrix.

The reflection coefficient $s_{00}^{\varepsilon}$ ought to be regarded as a
scattering matrix of size $1\times1$ in view of the only one couple of waves
(\ref{13}) which are able to drive energy along the waveguide (\ref{11}). For
example, dealing with the next couple of waves%
\begin{equation}
v_{1}^{\varepsilon\pm}\left(  x\right)  =\left(  k_{1}^{\varepsilon}\right)
^{-1/2}e^{\pm k_{1}^{\varepsilon}x_{1}}\cos\left(  \pi x_{2}\right)  ,\text{
\ \ \ }k_{1}^{\varepsilon}=\sqrt{\pi^{2}-\lambda^{\varepsilon}} \label{27}%
\end{equation}
which are decaying $(-)$ and growing $\left(  +\right)  ,$ one readily finds
that
\begin{equation}
q\left(  v_{1}^{\varepsilon\pm},v_{1}^{\varepsilon\pm}\right)  =0 \label{28}%
\end{equation}
but%
\begin{equation}
q\left(  v_{1}^{\varepsilon+},v_{1}^{\varepsilon-}\right)  =-q\left(
v_{1}^{\varepsilon-},v_{1}^{\varepsilon+}\right)  =1 \label{29}%
\end{equation}
As was observed in \cite[\S 5.6]{NaPl} and mentioned above, formula (\ref{28})
annihilates the projection on the $x_{3}$-axis of the Umov-Poyting vector and
therefore waves (\ref{27}) cannot drive energy. In the papers \cite{na246,
na252} (see \cite{na275} for general elliptic systems) the linear combinations
of exponential waves%
\begin{equation}
w_{1}^{\varepsilon\pm}\left(  x\right)  =2^{-1/2}\left(  v_{1}^{\varepsilon
+}\left(  x\right)  \mp iv_{1}^{\varepsilon-}\left(  x\right)  \right)
\label{30}%
\end{equation}
were introduced. It is remarkable that, thanks to (\ref{28}) and (\ref{29}),
waves (\ref{13}) and (\ref{30}) enjoy conditions (\ref{19}) with $p,q=0,1.$
The latter allows us to determine the solutions%
\begin{equation}
Z_{p}^{\varepsilon}\left(  x\right)  =w_{p}^{\varepsilon-}\left(  x\right)
+S_{0p}^{\varepsilon}w_{0}^{\varepsilon+}\left(  x\right)  +S_{1p}%
^{\varepsilon}w_{1}^{\varepsilon+}\left(  x\right)  +\widetilde{Z}%
_{p}^{\varepsilon}\left(  x\right)  ,\ \widetilde{Z}_{p}^{\varepsilon}\left(
x\right)  =O(e^{-x_{1}\sqrt{4\pi^{2}-\lambda^{\varepsilon}}}),\text{\ }p=0,1,
\label{31}%
\end{equation}
to compose the coefficient matrix $S^{\varepsilon}=\left(  S_{qp}%
^{\varepsilon}\right)  $ and to assure its unitary and symmetry property, cf.
(\ref{26}). Moreover, since $w_{p}^{\varepsilon+}\left(  x\right)
=\overline{w_{p}^{\varepsilon-}\left(  x\right)  }$, this matrix is symmetric,
see \cite[\S 2]{na489} again.

In Section \ref{sect5} we will give an operator formulation of the problem
(\ref{14}), at $\lambda^{\varepsilon}\in\left(  0,\pi^{2}\right)  $ with the
radiation condition%
\begin{equation}
U^{\varepsilon}\left(  x\right)  =c_{0}^{\varepsilon}w_{0}^{\varepsilon
+}\left(  x\right)  +c_{1}^{\varepsilon}w_{1}^{\varepsilon+}\left(  x\right)
+\widetilde{U}^{\varepsilon}\left(  x\right)  ,\ \ \ \widetilde{U}%
^{\varepsilon}\left(  x\right)  =O(e^{-x_{1}\sqrt{4\pi^{2}-\lambda
^{\varepsilon}}}). \label{32}%
\end{equation}
We recognize this condition as artificial because the right-hand side of
(\ref{32}) involves the exponentially growing wave $w_{1}^{\varepsilon
+}\left(  x\right)  $, see (\ref{30}) and (\ref{27}).

\subsection{A criterion for trapped modes.\label{sect2.3}}

A reason to consider problem (\ref{14}), (\ref{32}) and the augmented
scattering matrix $S^{\varepsilon}$ is explained by the following assertion.

\begin{theorem}
\label{TheoremA} Problem (\ref{9}),(\ref{10}) with the spectral parameter
(\ref{12}) has a trapped mode $u^{\varepsilon}\in H^{1}\left(  \Pi
_{+}^{\varepsilon}\right)  $ if and only if
\begin{equation}
S_{11}^{\varepsilon}=-1. \label{33}%
\end{equation}

\end{theorem}

In other words, equation (\ref{33}) provides a criterion for the existence of
a trapped mode in the one-sided waveguide (\ref{11}).

A proof of Theorem \ref{TheoremA} can be found, e.g., in \cite{na275} and
\cite[Thm 2]{na489} but, since the criterion (\ref{33}) plays the central rule
in our analysis, we here give the condensed proof.

The unitary property of $S^{\varepsilon}$ demonstrates that
\begin{equation}
S_{11}^{\varepsilon}=-1\qquad\Rightarrow\qquad S_{10}^{\varepsilon}%
=S_{01}^{\varepsilon}=0. \label{34}%
\end{equation}
Thus the solution (\ref{31}) with $p=1$ becomes a trapped mode because
formulas (\ref{30}) and (\ref{27}) assure that%
\[
Z_{1}^{\varepsilon}\left(  x\right)  =Z_{p}^{\varepsilon}\left(  x\right)
=w_{p}^{\varepsilon-}\left(  x\right)  -w_{1}^{\varepsilon+}\left(  x\right)
+\widetilde{Z}_{p}^{\varepsilon}\left(  x\right)  =-2^{1/2}iv_{1}%
^{\varepsilon-}\left(  x\right)  +\widetilde{Z}_{1}^{\varepsilon}\left(
x\right)  =O(e^{-x_{1}k_{1}^{\varepsilon}}).
\]
Hence, (\ref{33}) is a sufficient condition. To verify the necessity, we first
assume that the decomposition%
\begin{equation}
U^{\varepsilon}\left(  x\right)  =c^{\varepsilon}v^{\varepsilon-}\left(
x\right)  +\widetilde{U}^{\varepsilon}\left(  x\right)  \label{35}%
\end{equation}
of a trapped mode $U^{\varepsilon}\in H^{1}\left(  \Pi_{+}^{\varepsilon
}\right)  $ has a coefficient $c^{\varepsilon}\neq0.$

Then $U^{\varepsilon}$ becomes a linear combination of solutions (\ref{31}),
namely, according to (\ref{30}), we have%
\begin{align*}
U^{\varepsilon}  &  =C_{0}^{\varepsilon}Z_{0}^{\varepsilon}+C_{1}%
^{\varepsilon}Z_{1}^{\varepsilon}=C_{0}^{\varepsilon}w_{0}^{\varepsilon
-}+\left(  S_{00}^{\varepsilon}C_{0}^{\varepsilon}+S_{01}^{\varepsilon}%
C_{1}^{\varepsilon}\right)  w_{0}^{\varepsilon+}+\\
&  +2^{-1/2}\left(  v_{1}^{+}-iv_{1}^{-}\right)  C_{1}^{\varepsilon}%
+2^{-1/2}\left(  v_{1}^{+}+iv_{1}^{-}\right)  \left(  S_{10}^{\varepsilon
}C_{0}^{\varepsilon}+S_{11}^{\varepsilon}C_{1}^{\varepsilon}\right)
+\widetilde{U}^{\varepsilon}.
\end{align*}
Owing to the exponential decay of $U^{\varepsilon}$, coefficients of the
oscillating waves $w_{0}^{\varepsilon+}$ must vanish so that $C_{0}%
^{\varepsilon}=0,$\ 

$S_{01}^{\varepsilon}C_{1}^{\varepsilon}=0.$ Moreover, coefficients of the
exponential waves $v_{0}^{\varepsilon+}$ and $v_{0}^{\varepsilon-}$,
respectively, are $2^{-1/2}\left(  S_{11}^{\varepsilon}+1\right)
C_{1}^{\varepsilon}=0$\ and\ $2^{-1/2}\left(  S_{11}^{\varepsilon}-1\right)
C_{1}^{\varepsilon}=c^{\varepsilon}.$ Recalling our assumption $c^{\varepsilon
}\neq0,$ we see that $C_{1}^{\varepsilon}=-2^{-1/2}c^{\varepsilon}\neq0$ and,
therefore, (\ref{33}) holds true.

If $c^{\varepsilon}=0$ in (\ref{35}), the trapped mode $U^{\varepsilon}\left(
x\right)  $ gains very fast decay rate $O(e^{-x_{1}\sqrt{4\pi^{2}%
-\lambda^{\varepsilon}}}).$ In Section \ref{sect5.3} we will show with a new
argument that such trapped modes do not exist for a small $\varepsilon.$

\begin{remark}
\label{RemarkSO}. The relationship between the augmented scattering matrix and
the reflection coefficient in (\ref{24}) looks as follows:%
\begin{equation}
s_{00}^{\varepsilon}=S_{00}^{\varepsilon}-S_{01}^{\varepsilon}\left(
S_{11}^{\varepsilon}+1\right)  ^{-1}S_{10}^{\varepsilon}, \label{36}%
\end{equation}
see, e.g., \cite[Thm 3]{na489}. Note that, in view of (\ref{34}), the last
term in (\ref{36}) becomes null in the case $S_{11}^{\varepsilon}=-1$ when
$s_{00}^{\varepsilon}=S_{00}^{\varepsilon}.\ \ \ \ \ \boxtimes$
\end{remark}

\section{Formal asymptotics of the augmented scattering matrix\label{sect3}}

\subsection{Step-shaped perturbation of boundaries.\label{sect3.1}}

In this section we derive asymptotic expansions by means of a formal
asymptotic analysis and postpone their justification to Section \ref{sect6}.

Perturbation of the straight waveguide drawn in fig. \ref{f1},b and in fig.
\ref{f3}, b, ought to be regarded as a combination of regular and singular
perturbations, see, e.g., \cite[Ch.XII, \S \ 6.5]{Kato} and \cite[Ch. 2 and
4]{MaNaPl}, respectively. For a regular perturbation of the boundary, an
appropriate change of variables, which differs from the identity in magnitude
$O\left(  \varepsilon\right)  $ only, is usually applied in order to convert
the perturbed domain into the reference one. In this way, differential
operators in the problem gain small perturbations but asymptotics can be
constructed with the help of standard iterative procedures like decompositions
of a perturbed operator in the Neumann series.

\begin{figure}
[ptb]
\begin{center}
\ifcase\msipdfoutput
\includegraphics[
height=1.6475in,
width=3.7948in
]%
{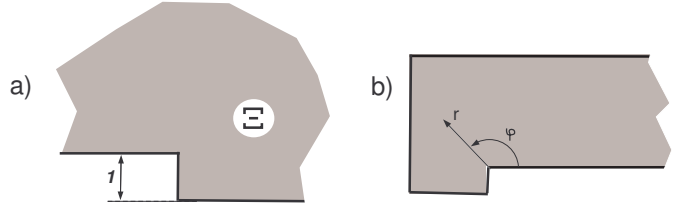}%
\else
\includegraphics[
height=1.6475in,
width=3.7948in
]%
{C:/Users/Desktop/Documents/Google Drive/Articoli/Articoli-in-Preparazione/CardoneDuranteNazarov1/graphics/CaDuNa4__5.pdf}%
\fi
\caption{A distorted half-plane (a) to describe the boundary layer near the
ledge (b)}%
\label{f4}%
\end{center}
\end{figure}

Singular perturbations of boundaries need much more delicate analysis because
they require for a description of asymptotics in the stretched coordinates
which, for the domain $\Pi_{+}^{\varepsilon}=\Pi_{l+}^{\varepsilon}$, see
(\ref{11}), take the form
\begin{equation}
\xi=\left(  \xi_{1},\xi_{2}\right)  =\varepsilon^{-1}\left(  x_{1}%
-l,x_{2}\right)  . \label{38}%
\end{equation}
Notice that the change $x\longrightarrow\xi$ and setting $\varepsilon=0$
transform $\Pi_{+}^{\varepsilon}$ into the upper half-plane $\mathbb{R}%
_{+}^{2}$ with a semi-infinite step, fig. \ref{f4}, a,
\begin{equation}
\Xi=\left\{  \xi\in\mathbb{R}^{2}:\xi_{2}>0\text{ \ for \ }\xi_{1}\leq0\text{
and }\xi_{2}>-1\text{ \ for \ }\xi_{1}>0\right\}  . \label{39}%
\end{equation}
%

As a result, the singular perturbation of the waveguide wall gives rise to the
boundary layer phenomenon described by solutions to the following problem:%
\begin{equation}
-\Delta_{\xi}v\left(  \xi\right)  =0,\text{ \ }\xi\in\Xi,\text{ \ \ }%
\partial_{\nu\left(  \xi\right)  }v\left(  \xi\right)  =g\left(  \xi\right)
,\text{ \ \ \ }\xi\in\partial\Xi. \label{40}%
\end{equation}
The Laplace equation is caused by the relation $\Delta_{x}+\lambda
=\varepsilon^{-2}\Delta_{\xi}+\lambda$ which singles out the Laplacian as the
main asymptotic part of the Helmholtz operator. The Neumann problem (\ref{40})
with a compactly supported datum $g$ admits a solution $v\left(  \xi\right)
=O(\left\vert \xi\right\vert ^{-1})$ as $\left\vert \xi\right\vert
\longrightarrow+\infty$ provided $\int_{\partial\Xi}g\left(  \xi\right)
ds_{\xi}=0.$ Otherwise, a solution grows at infinity like $c\ln\left\vert
\xi\right\vert $ and loses the intrinsic decay property of a boundary layer so
that traditional asymptotic procedures become much more sophisticated, see
\cite[Ch. 2, 4]{MaNaPl} and \cite{Ilin}. However, we will see that in our
particular problem the boundary perturbation does not affect the main
asymptotic term and the first correction term does not include a boundary layer.

\subsection{Asymptotic procedure\label{sect3.2}}

We search for an eigenvalue of the problem (\ref{9}), (\ref{10}) in the form%
\begin{equation}
\lambda^{\varepsilon}=\pi^{2}-\varepsilon^{2}\mu\label{42}%
\end{equation}
where the correction coefficient $\mu>0$ is to be found in Section
\ref{sect4}. Recalling the normalization factors in (\ref{27}) and (\ref{13})%
\begin{equation}
\left(  k_{1}^{\varepsilon}\right)  ^{-1/2}=\varepsilon^{-1/2}\mu
^{-1/4}+O(\varepsilon^{1/2}),\text{ \ \ \ }\left(  2k^{\varepsilon}\right)
^{-1/2}=\left(  2\pi\right)  ^{-1/2}+O(\varepsilon^{2}), \label{43}%
\end{equation}
we guess at the following asymptotic ans\"{a}tze for entries of the augmented
scattering matrix%
\begin{equation}
S_{11}^{\varepsilon}=S_{11}^{0}+\varepsilon S_{11}^{\prime}+\widetilde{S}%
_{11}^{\varepsilon},\text{ \ \ }S_{01}^{\varepsilon}=\varepsilon^{1/2}%
S_{01}^{0}+\varepsilon^{3/2}S_{01}^{\prime}+\widetilde{S}_{01}^{\varepsilon}
\label{44}%
\end{equation}
but aim to calculate the terms $S_{p1}^{0}$ and $S_{p1}^{\prime}$ only. We
further estimate the remainders in Section \ref{sect6}.

Using definitions of waves (\ref{13}) and (\ref{27}), (\ref{30}), we take
relations (\ref{43}) and (\ref{44}) into account and rewrite the decomposition
(\ref{31}) of the special solution $Z_{1}^{\varepsilon}$ as follows:%
\begin{align}
Z_{1}^{\varepsilon}\left(  x\right)   &  =\varepsilon^{-1/2}\left(
4\mu\right)  ^{-1/4}\cos\left(  \pi x_{2}\right)  (1+i+S_{11}^{0}\left(
1-i\right)  +\label{45}\\
&  +\varepsilon\left(  S_{11}^{\prime}\left(  1-i\right)  +x_{1}\sqrt{\mu
}\left(  1-i+S_{11}^{0}\left(  1+i\right)  \right)  +...\right)  +\nonumber\\
&  +\varepsilon^{1/2}\left(  2\pi\right)  ^{-1/2}\left(  S_{01}^{0}%
+\varepsilon S_{01}^{\prime}+...\right)  \left(  e^{i\pi x_{1}}+...\right)
.\nonumber
\end{align}
Here and everywhere in this section, ellipses stand for lower-order terms
inessential in our formal asymptotic analysis. In (\ref{45}), the Taylor
formula%
\begin{equation}
e^{k_{1}^{\varepsilon}x_{1}}\mp ie^{-k_{1}^{\varepsilon}x_{1}}=\left(  1\mp
i\right)  +\varepsilon x_{1}\sqrt{\mu}\left(  1\pm i\right)  +O\left(
\varepsilon^{2}x_{1}^{2}\right)  \label{46}%
\end{equation}
was applied so that expansion (\ref{46}) becomes meaningful under the
restriction $x_{1}<R$ with a fixed $R$, i.e. for $x\in\Pi^{\varepsilon}\left(
R\right)  .$

In view of the above observation we employ the method of matched asymptotic
expansions, cf. \cite{VanDyke,Ilin}, in the interpretation \cite{na489,
na457}. Namely, we regard (\ref{45}) as an outer expansion and introduce the
inner expansion%
\begin{equation}
Z_{1}^{\varepsilon}\left(  x\right)  =\varepsilon^{-1/2}Z_{1}^{0}\left(
x\right)  +\varepsilon^{1/2}Z_{1}^{\prime}\left(  x\right)  +... \label{47}%
\end{equation}
At the same time, coefficients of $\varepsilon^{-1/2}$ and $\varepsilon^{1/2}$
on the right-hand side of (\ref{45}) exhibit a behavior at infinity of the
terms $Z_{1}^{0}$ and $Z_{1}^{\prime}$ in (\ref{47}) because the upper bound
$R$ for $x_{1}$ can be chosen arbitrary large. Thus, they must satisfy%
\begin{align}
Z_{1}^{0}\left(  x\right)   &  =\left(  4\mu\right)  ^{-1/4}\cos\left(  \pi
x_{2}\right)  \left(  1+i+S_{11}^{0}\left(  1-i\right)  \right)
+...\label{Z0}\\
Z_{1}^{\prime}\left(  x\right)   &  =\left(  4\mu\right)  ^{-1/4}\cos\left(
\pi x_{2}\right)  S_{11}^{\prime}\left(  1-i\right)  +x_{1}\sqrt{\mu}\left(
1-i+S_{11}^{0}\left(  1+i\right)  \right)  +S_{01}^{0}\left(  2\pi\right)
^{-1/2}e^{i\pi x_{1}}+... \label{Z1}%
\end{align}

The formal passage to $\varepsilon=0$ transforms the waveguide (\ref{11}) into
the semi-infinite strip $\Pi_{+}^{0}=\mathbb{R\times}\left(  0,1\right)  $
while due to (\ref{42}) the Neumann problem (\ref{9}), (\ref{10}) converts
into%
\begin{align}
-\Delta u^{0}\left(  x\right)   &  =\pi^{2}u^{0}\left(  x\right)  ,\text{
\ }x\in\Pi_{+}^{0},\label{481}\\
\partial_{\nu}u^{0}\left(  x\right)   &  =0,\text{ \ \ }x\in\partial\Pi
_{+}^{0}. \label{482}%
\end{align}
This limit problem has two bounded solutions%
\begin{align}
u_{0}^{0}\left(  x\right)   &  =\frac{1}{2}\left(  e^{-i\pi x_{1}}+e^{i\pi
x_{1}}\right)  =\cos\left(  \pi x_{1}\right)  ,\label{491}\\
u_{1}^{0}\left(  x\right)   &  =\cos\left(  \pi x_{2}\right)  . \label{492}%
\end{align}
Comparing with (\ref{Z0}), we set%
\begin{equation}
Z_{1}^{0}\left(  x\right)  =\left(  4\mu\right)  ^{-1/4}\left(  1+i+S_{11}%
^{0}\left(  1-i\right)  \right)  \cos\left(  \pi x_{2}\right)  . \label{50}%
\end{equation}
Since $\lambda^{\varepsilon}=\pi^{2}+O\left(  \varepsilon^{2}\right)  ,$ the
function $Z_{1}^{\prime}$ also satisfies the homogeneous equation (\ref{481})
but the Neumann condition becomes inhomogeneous because of the boundary
perturbation. For $p=1$, we have
\begin{equation}
-\Delta Z_{p}^{\prime}\left(  x\right)  =\pi^{2}Z_{p}^{\prime}\left(
x\right)  ,\text{ \ \ }x\in\Pi_{+}^{0},\ \ \ \ \ \ \ \ \ \ \ \ \partial_{\nu
}Z_{p}^{\prime}\left(  x\right)  =g_{p}\left(  x\right)  ,\text{ \ \ }%
x\in\partial\Pi_{+}^{0}. \label{51}%
\end{equation}
To determine the datum $g_{1},$ we observe that the function (\ref{50})
satisfies the Neumann condition (\ref{10}) everywhere on $\partial
\Pi^{\varepsilon}$, except at the lower side $\Upsilon^{\varepsilon}=\left\{
x:x_{1}\in\left(  0,l\right)  ,\text{ \ \ }x_{2}=-\varepsilon\right\}  $ of
the box-shaped perturbation $\varpi_{+}^{\varepsilon}=\left(  0,l\right)
\times\left(  -\varepsilon,0\right]  $ in (\ref{11}). Furthermore, we obtain
\begin{align}
\partial_{\nu}Z_{1}^{p}\left(  x_{1},-\varepsilon\right)   &  =-\partial
_{2}Z_{1}^{0}\left(  x_{1},-\varepsilon\right)  =\left(  4\mu\right)
^{-1/4}\pi\sin\left(  -\pi\varepsilon\right)  \left(  1+i+S_{11}^{0}\left(
1-i\right)  \right)  =\label{5111}\\
&  =-\varepsilon\left(  4\mu\right)  ^{-1/4}\pi^{2}\left(  1+i+S_{11}%
^{0}\left(  1-i\right)  \right)  +O\left(  \varepsilon^{3}\right)
=:-\varepsilon G_{1}^{\prime}+O\left(  \varepsilon^{3}\right) \nonumber
\end{align}
and, hence,%
\begin{equation}
g_{p}\left(  x\right)  =\left\{
\begin{array}
[c]{l}%
0,\text{ \ \ }x\in\partial\Pi_{+}^{0}\setminus\overline{\Upsilon^{0}},\\
G_{p}^{\prime},\text{ \ \ }x\in\Upsilon^{0}.
\end{array}
\right.  \label{52}%
\end{equation}
Although the Neumann datum (\ref{52}) is not smooth and has a jump at the
point $x=\left(  l,0\right)  ,$ the problem (\ref{51}) with $p=1$ has a
solution in $H_{loc}^{1}(\overline{\Pi_{+}^{0}})$ such that
\begin{equation}
Z_{p}^{\prime}\left(  x\right)  =C_{p}e^{i\pi x_{1}}+\left(  C_{p}^{0}%
+x_{1}C_{p}^{1}\right)  \cos\left(  \pi x_{2}\right)  +\widetilde{Z}%
_{p}^{\prime}\left(  x\right)  ,\text{ \ }\widetilde{Z}_{p}^{\prime}\left(
x\right)  =O(e^{-x_{1}\sqrt{3}\pi}). \label{5222}%
\end{equation}
This fact is a direct consequence of the elliptic theory in domains with
cylindrical outlets to infinity (see the key works \cite{AgNi, Ko, MaPl1,
MaPl2} and, e.g., the monographs \cite{KoMaRo, NaPl}, but also may be derived
by the Fourier method after splitting $\Pi_{+}^{0}$ into the rectangle
$\left(  0,l\right)  \times\left(  0,1\right)  $ and the semi-strip $\left(
l,+\infty\right)  \times\left(  0,1\right)  .$ A simple explanation how to
apply the above-mentioned theory can be found in the introductory chapter 2 of
\cite{NaPl}, the review \cite{na417} and Section \ref{sect5} of this paper.
The solution (\ref{5222}) is defined up to the term $c\cos\left(  \pi
x_{2}\right)  $ and, therefore, the coefficient $C^{0}$ can be taken
arbitrary. Other coefficients in (\ref{5222}) are computed by application of
the Green formula in the long ($R$ is big) rectangle $\Pi^{0}\left(  R\right)
=\left(  0,R\right)  \times\left(  0,1\right)  .$

Indeed, we send $R$ to $+\infty$ and obtain%
\begin{align}
0  &  =\underset{R\rightarrow+\infty}{\lim}\int_{\Pi_{+}^{0}\left(  R\right)
}\left(  u_{1}^{0}\left(  x\right)  \left(  \Delta+\pi^{2}\right)
Z_{1}^{\prime}\left(  x\right)  -Z_{1}^{\prime}\left(  x\right)  \left(
\Delta+\pi^{2}\right)  u_{1}^{0}\left(  x\right)  \right)  dx=\label{54}\\
&  =\underset{R\rightarrow+\infty}{\lim}\int_{0}^{1}\cos\left(  \pi
x_{2}\right)  \partial_{1}Z_{1}^{\prime}\left(  R,x_{2}\right)  dx_{2}%
-\int_{0}^{l}\cos\left(  \pi0\right)  \partial_{2}Z_{1}^{\prime}\left(
x_{1},0\right)  dx_{1}=\nonumber\\
&  =\frac{1}{2}C_{1}^{1}+\left(  4\pi\right)  ^{-1/4}\pi^{2}l\left(
1+i+S_{11}^{0}\left(  1-i\right)  \right)  .\nonumber
\end{align}
In the same way we deal with the functions (\ref{491}) and (\ref{5222}) that
results in the equality%
\begin{gather}
0=\underset{R\rightarrow+\infty}{\lim}%
{\displaystyle\int_{0}^{1}}
\left.  \left(  \cos\left(  \pi x_{1}\right)  \partial_{1}Z_{1}^{\prime
}\left(  x\right)  dx_{2}-U_{1}^{\prime}\left(  x\right)  \partial_{1}%
\cos\left(  \pi x_{1}\right)  \right)  \right\vert _{x_{1}=R}dx_{2}%
-\label{55}\\
-%
{\displaystyle\int_{0}^{l}}
\cos\left(  \pi x_{1}\right)  \partial_{2}Z_{1}^{\prime}\left(  x_{1}%
,0\right)  dx_{1}=i\pi C_{1}+\left(  4\mu\right)  ^{-1/4}\pi^{2}\left(
1+i+S_{11}^{0}\left(  1-i\right)  \right)
{\displaystyle\int_{0}^{l}}
\cos\left(  \pi x_{1}\right)  dx_{1}.\nonumber
\end{gather}
Comparing (\ref{Z1}) and (\ref{5222}), we arrive at the relations
\begin{equation}
\left(  2\pi\right)  ^{-1/2}S_{01}^{0}=C_{1},\text{ \ \ }\left(  4\mu\right)
^{-1/4}\sqrt{\mu}\left(  1-i+S_{11}^{0}\left(  1+i\right)  \right)  =C_{1}%
^{1},\text{ \ \ }S_{11}^{\prime}\left(  1-i\right)  =C_{1}^{0}\text{\ }
\label{56}%
\end{equation}
which together with our calculations (\ref{55}) and (\ref{56}) give us the
following formulas:%
\begin{align}
S_{01}^{0}  &  =\left(  4\mu\right)  ^{-1/4}\left(  2\pi\right)  ^{1/2}\pi
i\left(  1+i+S_{11}^{0}\left(  1-i\right)  \right)  \int_{0}^{l}\cos\left(
\pi x_{1}\right)  dx_{1}=\label{57}\\
&  =-\left(  4\mu\right)  ^{-1/4}\left(  2\pi\right)  ^{1/2}\left(
1-i-S_{11}^{0}\left(  1+i\right)  \right)  \sin\left(  \pi l\right)
\nonumber\\
\sqrt{\mu}\left(  1-i+S_{11}^{0}\left(  1+i\right)  \right)   &  =-2\pi
^{2}l\left(  1+i+S_{11}^{0}\left(  1-i\right)  \right)  \Rightarrow
\label{58}\\
&  \Rightarrow S_{11}^{0}=-\frac{\sqrt{\mu}\left(  1-i\right)  2\pi
^{2}l\left(  1+i\right)  }{\sqrt{\mu}\left(  1+i\right)  2\pi^{2}l\left(
1-i\right)  }=-\frac{4\pi^{2}l\sqrt{\mu}+i\left(  4\pi^{4}l^{2}-\mu\right)
}{4\pi^{4}l^{2}+\mu}.\nonumber
\end{align}
We emphasize that $\mu=4\pi^{4}l^{2}\Rightarrow S_{11}^{0}=-1.$

The necessary computations are completed. It should be mentioned that, to
determine the correction terms $S_{11}^{\prime}$ and $S_{01}^{\prime}$ in the
ans\"{a}tze (\ref{44}), one has to make another step in our asymptotic
procedure, see the next section, but they are of no further use.

\subsection{The detailed asymptotic procedure\label{sect3.3}}

Let us construct asymptotics of the entries
\begin{equation}
S_{00}^{\varepsilon}=S_{00}^{0}+\varepsilon S_{00}^{\prime}+\widetilde{S}%
_{00}^{\varepsilon}\text{, \ \ \ }S_{10}^{\varepsilon}=\varepsilon^{1/2}%
S_{10}^{0}+\varepsilon^{3/2}S_{10}^{\prime}+\widetilde{S}_{10}^{\varepsilon}
\label{60}%
\end{equation}
in the augmented scattering matrix. These asymptotics are not of very
importance in our analysis of eigenvalues but are much more representational
than the analogous formulas (\ref{44}) because the ansatz for $Z_{0}%
^{\varepsilon}$ indeed involves the boundary layer concentrated near the ledge
of the box-shape perturbation in (\ref{11}).

Using (\ref{43}), (\ref{60}) and (\ref{46}), we rewrite the decomposition
(\ref{31}) of $Z_{0}^{\varepsilon}\left(  x\right)  $ as follows:%
\begin{align}
Z_{0}^{\varepsilon}\left(  x\right)   &  =\left(  2\pi\right)  ^{-1/2}\left(
e^{-i\pi x_{1}}+S_{00}^{0}e^{i\pi x_{1}}+\varepsilon S_{00}^{\prime}e^{i\pi
x_{1}}+...\right)  +\label{61}\\
&  +\left(  4\mu\right)  ^{-1/4}\left(  S_{10}^{0}+\varepsilon S_{10}^{\prime
}+...\right)  \cos\left(  \pi x_{2}\right)  \left(  1-i+\varepsilon x_{1}%
\sqrt{\mu}\left(  1+i\right)  +...\right)  .\nonumber
\end{align}
Then we accept in a finite part of $\Pi_{+}^{\varepsilon}$ the ansatz%
\begin{equation}
Z_{0}^{\varepsilon}\left(  x\right)  =Z_{0}^{0}\left(  x\right)  +\varepsilon
Z_{0}^{\prime}\left(  x\right)  +... \label{62}%
\end{equation}
Referring to (\ref{61}), we fix the behavior of its terms at infinity%
\begin{align}
Z_{0}^{0}\left(  x\right)   &  =\left(  2\pi\right)  ^{-1/2}\left(  e^{-i\pi
x_{1}}+S_{00}^{0}e^{i\pi x_{1}}\right)  +\left(  4\mu\right)  ^{-1/4}%
S_{10}^{0}\left(  1-i\right)  \cos\left(  \pi x_{2}\right)  +...,\label{Z00}\\
Z_{0}^{\prime}\left(  x\right)   &  =\left(  2\pi\right)  ^{-1/2}%
S_{00}^{\prime}e^{i\pi x_{1}}+\left(  4\mu\right)  ^{-1/4}\cos\left(  \pi
x_{2}\right)  (S_{10}^{\prime}\left(  1-i\right)  +x_{1}\sqrt{\mu}S_{10}%
^{0}\left(  1+i\right)  )+... \label{Z01}%
\end{align}
Comparing (\ref{Z00}) with (\ref{491}), (\ref{492}), we conclude that
\begin{equation}
Z_{0}^{0}\left(  x\right)  =\left(  2\pi\right)  ^{-1/2}\left(  e^{-i\pi
x_{1}}+e^{i\pi x_{2}}\right)  +\left(  4\mu\right)  ^{-1/4}S_{10}^{0}\left(
1-i\right)  \cos\left(  \pi x_{2}\right)  \label{63}%
\end{equation}
and, hence,%
\begin{equation}
S_{00}^{0}=1 \label{64}%
\end{equation}

Let us describe the correction terms in (\ref{60}). The function (\ref{63})
satisfies the equation (\ref{9}) with $\lambda^{\varepsilon}=\pi^{2}$ and
leaves the discrepancies
\begin{gather}
\partial_{\nu}Z_{0}^{0}\left(  x_{1},-\varepsilon\right)  =-\partial_{2}%
Z_{0}^{0}\left(  x_{1},-\varepsilon\right)  =-\varepsilon\left(  4\mu\right)
^{-1/4}\pi^{2}S_{10}^{0}\left(  1-i\right)  +O\left(  \varepsilon^{3}\right)
=\label{Z22}\\
=-\varepsilon G_{0}^{\prime}+O\left(  \varepsilon^{3}\right)  ,x_{1}\in\left(
0,l\right) \nonumber\\
\partial_{\nu}Z_{0}^{0}\left(  l,x_{2}\right)  =\partial_{1}Z_{0}^{0}\left(
l,x_{2}\right)  =-\left(  2\pi\right)  ^{1/2}\sin\left(  \pi l\right)  ,\text{
\ \ }x_{2}\in\left(  -\varepsilon,0\right)  \label{Z02}%
\end{gather}
in the boundary condition (\ref{10}) on the big $\Upsilon^{\varepsilon}$ and
small $\upsilon^{\varepsilon}=\left\{  x:x_{1}=l,\text{ }x_{2}\in\left(
-\varepsilon,0\right)  \right\}  $ sides of the rectangle $\varpi
_{+}^{\varepsilon}$, respectively. The discrepancy (\ref{Z22}) is similar to
(\ref{5111}) and appears as the datum (\ref{52}) in the problem (\ref{51})
with $p=0.$ To compensate for (\ref{Z02}), we need the boundary layer
\begin{equation}
V_{0}^{0}\left(  \xi\right)  =\left(  2\pi\right)  ^{1/2}\sin\left(  \pi
l\right)  v\left(  \xi\right)  \label{Z03}%
\end{equation}
where $\xi$ are stretched coordinates (\ref{38}) and $v$ is a solution of the
Neumann problem (\ref{40}) in the unbounded domain (\ref{39}) with the
right-hand side%
\[
g\left(  \xi\right)  =\left\{
\begin{array}
[c]{l}%
0,\text{ \ \ }\xi_{1}\neq0,\\
1,\text{ \ \ }\xi_{1}=0,\text{ \ \ }\xi_{2}\in\left(  -1,0\right)  ,
\end{array}
\right.  \ \text{for }\xi\in\partial\Xi.
\]
This solution, of course, can be constructed by an appropriate conformal
mapping but we only need its presentation at infinity%
\begin{equation}
v\left(  \xi\right)  =(B/\pi)\ln(1/\left\vert \xi\right\vert )+c+O\left(
1/\left\vert \xi\right\vert \right)  ,\text{ \ \ \ \ }\left\vert
\xi\right\vert \rightarrow\infty. \label{Z05}%
\end{equation}
The constant $c$ is arbitrary but the coefficient $B$ can be computed by the
Green formula in the truncated domain $\Xi\left(  R\right)  =\left\{  \xi
\in\Xi:\left\vert \xi\right\vert <R\right\}  $ with $R\rightarrow+\infty:$%
\begin{align}
0  &  =\underset{R\rightarrow+\infty}{\lim}R%
{\displaystyle\int_{0}^{\pi+\arcsin\left(  1/R\right)  }}
\dfrac{\partial v}{\partial\rho}\left(  \xi\right)  d\varphi+%
{\displaystyle\int_{0}^{1}}
\dfrac{\partial v}{\partial\xi_{1}}\left(  0,\xi_{2}\right)  d\xi
_{2}=\label{Z08}\\
&  =-\dfrac{B}{\pi}%
{\displaystyle\int_{0}^{\pi}}
d\varphi+\int_{0}^{1}d\xi_{2}=-B+1\Rightarrow B=1.\nonumber
\end{align}
Here, $\left(  \rho,\varphi\right)  $ is the polar coordinates system.

We fix $c=-\pi^{-1}\ln\varepsilon$ in (\ref{Z05}) and observe that
\begin{equation}
v\left(  \xi\right)  =(1/\pi)\left(  \ln\left(  1/\rho\right)  -\ln
\varepsilon\right)  +O\left(  1/\rho^{2}\right)  =(1/\pi)\ln\left(
1/r\right)  +O\left(  \varepsilon^{2}/r^{2}\right)  \label{Z06}%
\end{equation}
in the polar coordinate system $\left(  r,\varphi\right)  $ centered at the
point $x=\left(  l,0\right)  \in\partial\Pi^{\varepsilon},$ see fig. \ref{f4}, b.

Applying the method of matched asymptotic expansions in the traditional
manner, cf. \cite{Ilin, VanDyke}, \cite[Ch. 2]{MaNaPl}, as well as \cite[Ch.
5]{MaNaPl}, \cite{na178} for ledge-shaped perturbation of domains, we consider
(\ref{62}) as an outer expansion in a finite part of the waveguide while
$\varepsilon V_{0}^{0}\left(  \xi\right)  $ becomes the main term of the inner
expansion in the vicinity of the ledge of the box-shaped perturbation in
(\ref{11}). In view of (\ref{Z03}) and (\ref{Z06}), the standard matching
procedure proposes $Z_{0}^{\prime}$ as a singular solution of the homogeneous
problem (\ref{481}), (\ref{482}) with the following asymptotic condition at
the point $x=\left(  l,0\right)  \in\partial\Pi^{0}:$%
\begin{equation}
Z_{0}^{\prime}\left(  x\right)  =\left(  2/\pi\right)  ^{1/2}\sin\left(  \pi
l\right)  \ln(1/r)+c+O\left(  r\right)  \text{ \ \ \ \ for }r\rightarrow+0.
\label{Z07}%
\end{equation}

Arguing in the same way as for the function $Z_{1}^{\prime},$ we conclude that
the problem (\ref{51}), (\ref{Z07}) has a solution $Z_{0}^{\prime}$ which
admits the representation (\ref{5222}) with $p=0$ under the restriction
$x_{1}\geq R>l$ needed due to the logarithmic singularity in (\ref{Z07}). The
coefficient $C_{0}^{0}$ is arbitrary but $C_{0}$ and $C_{0}^{1}$ can be
computed again by means of the Green formula in the domain $\Pi_{+}^{0}\left(
R,\delta\right)  =\left\{  x\in\Pi_{+}^{0}:x_{1}<R,\text{\ }r>\delta\right\}
$ and the limit passage $R\rightarrow+\infty,$ $\delta\rightarrow+0,$ cf.
(\ref{54}), (\ref{55}) and (\ref{Z08}). Dealing with $u_{0}^{0}$ and
$Z_{0}^{\prime},$ we take into account the equality $\partial_{2}u_{0}%
^{0}\left(  x_{1},0\right)  =0$ and obtain that%
\begin{gather}
0=%
{\displaystyle\int_{\Pi\left(  R,\delta\right)  }}
\cos\left(  \pi x_{1}\right)  \left(  \Delta Z_{0}^{\prime}\left(  x\right)
+\pi^{2}Z_{0}^{\prime}\left(  x\right)  \right)
dx=\ \ \ \ \ \ \ \ \ \ \ \ \ \ \ \ \ \ \ \ \ \ \ \ \ \label{ZC}\\
=%
{\displaystyle\int_{0}^{1}}
\left.  \left(  \cos\left(  \pi x_{1}\right)  \dfrac{\partial Z_{0}^{\prime}%
}{\partial x_{1}}\left(  x\right)  +\pi\sin\left(  \pi x_{1}\right)
Z_{0}^{\prime}\left(  x\right)  \right)  \right\vert _{x_{1}=R}dx_{2}%
-\nonumber\\
-\delta%
{\displaystyle\int_{0}^{\pi}}
\left.  \left(  \cos\left(  \pi x_{1}\right)  \dfrac{\partial Z_{0}^{\prime}%
}{\partial r}\left(  x\right)  -Z_{0}^{\prime}\left(  x\right)  \dfrac
{\partial}{\partial r}\cos\left(  \pi x_{1}\right)  \right)  \right\vert
_{r=\delta}d\varphi=\nonumber\\
=i\pi C_{0}+\left(  \dfrac{\pi}{2}\right)  ^{1/2}\sin\left(  2\pi l\right)
+o\left(  1\right)  \Rightarrow C_{0}=i\left(  2\pi\right)  ^{-1/2}\sin\left(
2\pi l\right)  .\nonumber
\end{gather}
Inserting $u_{1}^{0}$ and $Z_{0}^{\prime}$ into the Green formula in $\Pi
_{+}^{0}\left(  R,\delta\right)  $, we take into account the inhomogeneous
boundary condition on $\Upsilon^{0}$ and derive that
\begin{equation}%
\begin{array}
[c]{c}%
0=%
{\displaystyle\int_{0}^{1}}
\left.  \cos\left(  \pi x_{2}\right)  \partial_{1}Z_{0}^{\prime}\left(
x\right)  \right\vert _{x_{1}=R}dx_{2}-%
{\displaystyle\int_{0}^{l-\delta}}
\left.  \cos\left(  \pi x_{2}\right)  \partial_{2}Z_{0}^{\prime}\left(
x\right)  \right\vert _{x_{2}=0}dx_{1}-\\
-\delta%
{\displaystyle\int_{0}^{\pi}}
\left.  \left(  \cos\left(  \pi x_{2}\right)  \partial_{r}Z_{0}^{\prime
}\left(  x\right)  -Z_{0}^{\prime}\left(  x\right)  \partial_{r}\cos\left(
\pi x_{2}\right)  \right)  \right\vert _{r=\delta}d\varphi=\\
=\dfrac{1}{2}C_{0}^{1}\left(  4\mu\right)  ^{-1/4}\pi^{2}lS_{10}^{0}\left(
1-i\right)  +\left(  2\pi\right)  ^{1/2}\sin\left(  \pi l\right)  +o\left(
1\right)  .
\end{array}
\label{ZC1}%
\end{equation}
We now compare coefficients in the expansions (\ref{5222}), $p=0$, and
(\ref{Z01}). According to the calculations (\ref{ZC}) and (\ref{ZC1}) we
derive the formulas%
\begin{align}
C_{0}  &  =\left(  2\pi\right)  ^{-1/2}S_{00}^{\prime}\Rightarrow
S_{00}^{\prime}=i\sin\left(  2\pi l\right)  ,\nonumber\\
C_{0}^{1}  &  =\left(  4\mu\right)  ^{-1/4}\sqrt{\mu}S_{10}^{0}\left(
1+i\right)  \Rightarrow\label{anti}\\
S_{10}^{0}  &  =-\dfrac{\left(  4\mu\right)  ^{1/4}2\left(  2\pi\right)
^{1/2}\sin\left(  \pi l\right)  }{\sqrt{\mu}\left(  1+i\right)  +2\pi
^{2}l\left(  1-i\right)  }=-\left(  4\mu\right)  ^{1/4}\left(  2\pi\right)
^{1/2}\frac{\sqrt{\mu}\left(  1-i\right)  +2\pi^{2}l\left(  1+i\right)  }%
{4\pi^{4}l^{2}+\mu}\sin\left(  \pi l\right)  .\nonumber
\end{align}
Calculation of coefficients in the ans\"{a}tze (\ref{44}) and (\ref{60}) is
completed. It is worth to underline that the expression (\ref{anti}) for the
main asymptotic term of $\varepsilon^{-1/2}S_{10}^{\varepsilon}=\varepsilon
^{-1/2}S_{01}^{\varepsilon}$ can be derived from the cumbersome relations
(\ref{57}) and (\ref{58}) as well.

\section{Detection of a trapped mode\label{sect4}}

\subsection{Reformulation of the criterion\label{sect4.1}}

We opt for the form%
\begin{equation}
\mu=4\pi^{4}l^{2}+\mathcal{\vartriangle}\mu,\text{ \ \ }l=\pi
k+\mathcal{\vartriangle}l \label{T1}%
\end{equation}
of the spectral and length parameters which support a trapped mode. Here,
$k\in\mathbb{N}$ is fixed but small $\mathcal{\vartriangle}\mu,$
$\mathcal{\vartriangle}l$ are to be determined. If $\mathcal{\vartriangle}%
\mu=0$ and $\mathcal{\vartriangle}l=0$, the equalities $S_{11}^{0}=-1$ and
$S_{01}^{0}=0$ hold due to (\ref{57}) and (\ref{anti}).

We purpose to choose the small increments $\mathcal{\vartriangle}\mu$ and
$\mathcal{\vartriangle}l$ in (\ref{T1}) such that the criterion (\ref{33}) for
the existence of a trapped mode is satisfied. Since $S_{11}^{\varepsilon}$ is
complex, the criterion furnishes two equations for two real parameters
$\mathcal{\vartriangle}\mu$ and $\mathcal{\vartriangle}l.$ It is convenient to
consider the other equations
\begin{equation}
\operatorname{Im}S_{11}^{\varepsilon}=0,\text{ \ \ \ }\operatorname{Re}%
S_{01}^{0}=0 \label{T2}%
\end{equation}
which, for small $\varepsilon$ and $\mathcal{\vartriangle}\mu,$
$\mathcal{\vartriangle}l,$ are equivalent to $S_{11}^{\varepsilon}=-1.$
Indeed, from formulas (\ref{58}), (\ref{anti}) and (\ref{64}) together with
estimates (\ref{est}) it follows that%
\begin{equation}
\left\vert S_{11}^{\varepsilon}+1\right\vert +\left\vert S_{00}^{\varepsilon
}-1\right\vert \leq c\left(  \varepsilon+\left\vert \mathcal{\vartriangle}%
\mu\right\vert +\left\vert \mathcal{\vartriangle}l\right\vert \right)
^{\delta},\text{ \ \ \ }\delta\in\left(  0,1\right)  . \label{T3}%
\end{equation}
Since $S^{\varepsilon}$ is unitary and symmetric, see Section \ref{sect2.2},
the second assumption in (\ref{T2}) means that $S_{01}^{\varepsilon}%
=S_{10}^{\varepsilon}=i\sigma$ with some $\sigma\in\mathbb{R}$ and,
furthermore,%
\begin{equation}
0=\overline{S_{00}^{\varepsilon}}S_{01}^{\varepsilon}+\overline{S_{10}%
^{\varepsilon}}S_{11}^{\varepsilon}=2i\sigma+O\left(  \left\vert
\sigma\right\vert \left(  \varepsilon+\left\vert \mathcal{\vartriangle}%
\mu\right\vert +\left\vert \mathcal{\vartriangle}l\right\vert \right)
^{\delta}\right)  . \label{T999}%
\end{equation}
Hence, $\sigma=0$ when $\mathcal{\vartriangle}\mu$, $\mathcal{\vartriangle}l$
and $\varepsilon>0$ are small so that $S_{01}^{\varepsilon}=0\Rightarrow
\left\vert S_{11}^{\varepsilon}\right\vert =-1$ due to (\ref{T3}),(\ref{T2}).

We have proved that (\ref{T2})$\Rightarrow$(\ref{33}) but the inverse
implication (\ref{33})$\Rightarrow$(\ref{T2}) is obvious.

\subsection{Solving the system of transcendental equations\label{sect4.2}}

By virtue of (\ref{58}), (\ref{S11}), and (\ref{T1}), the first equation
(\ref{T2}) turns into%
\begin{equation}
\mathcal{\vartriangle}\mu=-\varepsilon\left(  8\pi^{4}l^{2}%
+\mathcal{\vartriangle}\mu\right)  \operatorname{Im}\widehat{S}_{11}%
^{\varepsilon}. \label{T4}%
\end{equation}
The formulas (\ref{anti}), (\ref{S01}), and a simple algebraic calculation
convert the second equation (\ref{T2}) into%
\[
\sin l=\left(  4\mu\right)  ^{-1/4}\left(  2\pi\right)  ^{-1/2}\frac{4\pi
^{4}l^{2}+\mu}{4\pi^{2}l\sqrt{\mu}+2\mu}\varepsilon\operatorname{Re}%
\widehat{S}_{01}^{\varepsilon}%
\]
and thus%
\begin{equation}
\mathcal{\vartriangle}l=\arcsin\left(  \left(  -1\right)  ^{k}\left(
4\mu\right)  ^{-1/4}\left(  2\pi\right)  ^{-1/2}\frac{4\pi^{2}l^{2}+\mu}%
{2\pi^{2}l+\sqrt{\mu}}\varepsilon\operatorname{Re}\widehat{S}_{01}%
^{\varepsilon}\right)  . \label{T5}%
\end{equation}
We search for a solution $\left(  \mathcal{\vartriangle}\mu
,\mathcal{\vartriangle}l\right)  $ of the transcendental equations (\ref{T4}),
(\ref{T5}) in the closed disk%
\begin{equation}
\overline{\mathbb{B}_{\varrho}}=\left\{  \left(  \mathcal{\vartriangle}%
\mu,\mathcal{\vartriangle}l\right)  \in\mathbb{R}^{2}:\left\vert
\mathcal{\vartriangle}\mu\right\vert ^{2}+\left\vert \mathcal{\vartriangle
}l\right\vert ^{2}\leq\varrho^{2}\right\}  \label{T6}%
\end{equation}
and rewrite them in the condensed form%
\begin{equation}
\left(  \mathcal{\vartriangle}\mu,\mathcal{\vartriangle}l\right)
=T^{\varepsilon}\left(  \mathcal{\vartriangle}\mu,\mathcal{\vartriangle
}l\right)  \text{ \ \ in }\overline{\mathbb{B}_{\varrho}} \label{T7}%
\end{equation}
where $T^{\varepsilon}$ is a nonlinear operator involving asymptotic
remainders from formulas (\ref{S11}) and (\ref{S01}) for the augmented
scattering matrix $S^{\varepsilon}=S^{\varepsilon}\left(
\mathcal{\vartriangle}\mu,\mathcal{\vartriangle}l\right)  $. The estimates
(\ref{est}) and Proposition \ref{PropositionANAL} below demonstrate that the
operator is smooth in $\overline{\mathbb{B}_{\varrho}}$ with $\varrho
\leq\varrho_{0}$, $\varrho_{0}>0,$ and, furthermore,%
\[
\left\vert T^{\varepsilon}\left(  \mathcal{\vartriangle}\mu
,\mathcal{\vartriangle}l\right)  \right\vert =c_{\varrho}\varepsilon\left(
1+\left\vert \ln\varepsilon\right\vert \right)  ^{2}\text{ \ \ for }\left(
\mathcal{\vartriangle}\mu,\mathcal{\vartriangle}l\right)  \in\overline
{\mathbb{B}_{\varrho}}.
\]
Hence, for any $\varrho\leq\varrho_{0}$, there exists $\varepsilon\left(
\varrho\right)  >0$ such that $T^{\varepsilon}$ with $\varepsilon\in\left(
0,\varepsilon\left(  \varrho\right)  \right)  $ is a contraction operator in
the disk (\ref{T6}). By the Banach contraction principle, the abstract
equation (\ref{T7}) has a unique solution $\left(  \mathcal{\vartriangle}%
\mu,\mathcal{\vartriangle}l\right)  \in\overline{\mathbb{B}_{\varrho}}$ and
the estimate $\left\vert \mathcal{\vartriangle}\mu\right\vert +\left\vert
\mathcal{\vartriangle}l\right\vert \leq C\varepsilon\left(  1+\left\vert
\ln\varepsilon\right\vert \right)  ^{2}$ is valid. This solution depends on
$\varepsilon$ and determines the spectral and length parameters (\ref{T1})
supporting a trapped mode in the perturbed waveguide (\ref{11}) according to
the criterion (\ref{33}) from Theorem \ref{TheoremA} reformulated as (\ref{T2}).

\subsection{The main results\label{sect4.3}}

Based on the performed formal calculations, we will prove in the next three
sections the following existence and uniqueness theorems.

\begin{theorem}
\label{TheoremEX}Let $k\in\mathbb{N}$. There exist $\varepsilon_{k}$,
$c_{k}>0$ and $\mathcal{\vartriangle}\mu_{k}\left(  \varepsilon\right)  $,
$\mathcal{\vartriangle}l_{k}\left(  \varepsilon\right)  ,$ such that, for any
$\varepsilon\in\left(  0,\varepsilon_{k}\right)  $, the estimate%
\[
\left\vert \mathcal{\vartriangle}\mu_{k}\left(  \varepsilon\right)
\right\vert +\left\vert \mathcal{\vartriangle}l_{k}\left(  \varepsilon\right)
\right\vert \leq c_{k}\varepsilon\left(  1+\left\vert \ln\varepsilon
\right\vert \right)  ^{2}%
\]
is valid and the problem (\ref{1}), (\ref{2}) in the waveguide $\Pi_{l\left(
\varepsilon\right)  }^{\varepsilon}=\Pi\cup\varpi_{l_{k}\left(  \varepsilon
\right)  }^{\varepsilon}$ with the box-shaped perturbation (\ref{0}) of length
$2l_{k}\left(  \varepsilon\right)  =2\left(  \pi k+\mathcal{\vartriangle}%
_{k}l\left(  \varepsilon\right)  \right)  $ has an eigenvalue%
\begin{equation}
\lambda_{k}^{\varepsilon}=\pi^{2}-\varepsilon^{2}(4\pi^{4}\left(  \pi
k+\mathcal{\vartriangle}l_{k}\left(  \varepsilon\right)  \right)
^{2}+\mathcal{\vartriangle}\mu_{k}\left(  \varepsilon\right)  ). \label{T10}%
\end{equation}
The eigenvalue (\ref{T10}) is unique in the interval $\left(  0,\pi
^{2}\right)  $.
\end{theorem}

\begin{theorem}
\label{TheoremUN}Let $k\in\mathbb{N}$ and $\delta>0.$ There exist
$\varepsilon_{k}^{\delta}>0$ such that, for any $\varepsilon\in\left(
0,\varepsilon_{k}^{\delta}\right)  $, the waveguide $\Pi_{l}^{\varepsilon}$
with the length parameter
\begin{equation}
l\in\left[  \pi\left(  k-1\right)  +\delta,\pi\left(  k+1\right)
-\delta\right]  \label{T11}%
\end{equation}
does not support a trapped mode in the case $l\neq l_{k}\left(  \varepsilon
\right)  $ where $l_{k}\left(  \varepsilon\right)  $ is taken from Theorem
\ref{TheoremEX}.
\end{theorem}

\section{Weighted spaces with detached asymptotics\label{sect5}}

\subsection{Reformulation of the problem\label{sect5.1}}

Let $W_{\beta}^{1}\left(  \Pi_{+}^{\varepsilon}\right)  $ be the Kondratiev
(weighted Sobolev) space composed from functions $u^{\varepsilon}$ in
$H_{loc}^{1}\left(  \overline{\Pi_{+}^{\varepsilon}}\right)  $ with the finite
norm%
\begin{equation}
\left\Vert u^{\varepsilon};W_{\beta}^{1}\left(  \Pi_{+}^{\varepsilon}\right)
\right\Vert =\left\Vert e^{\beta x_{1}}u^{\varepsilon};H^{1}\left(  \Pi
_{+}^{\varepsilon}\right)  \right\Vert \label{K1}%
\end{equation}
where $\beta\in\mathbb{R}$ is the exponential weight index. If $\beta>0$,
functions in $W_{\beta}^{1}\left(  \Pi_{+}^{\varepsilon}\right)  $ decay at
infinity in the semi-infinite waveguide (\ref{11}) but in the case $\beta<0$ a
certain exponential growth is permitted while the rate of decay/growth is
ruled by $\beta.$ Clearly, $W_{0}^{1}\left(  \Pi_{+}^{\varepsilon}\right)
=H^{1}\left(  \Pi_{+}^{\varepsilon}\right)  $. The space $C_{c}^{\infty
}\left(  \overline{\Pi_{+}^{\varepsilon}}\right)  $ of smooth compactly
supported functions is dense in $W_{\beta}^{1}\left(  \Pi_{+}^{\varepsilon
}\right)  $ for any $\beta.$

By a solution of the problem (\ref{14}) in $W_{\sigma}^{1}\left(  \Pi
_{+}^{\varepsilon}\right)  $, $\sigma\in\mathbb{R}$, we understand a function
$u^{\varepsilon}\in W_{\sigma}^{1}\left(  \Pi_{+}^{\varepsilon}\right)  $
satisfying the integral identity%
\begin{equation}
\left(  \nabla u^{\varepsilon},\nabla v^{\varepsilon}\right)  _{\Pi
^{\varepsilon}}-\lambda^{\varepsilon}\left(  u^{\varepsilon},v^{\varepsilon
}\right)  _{\Pi^{\varepsilon}}=F^{\varepsilon}\left(  v^{\varepsilon}\right)
\text{ \ \ }\forall v^{\varepsilon}\in W_{-\sigma}^{1}\left(  \Pi
_{+}^{\varepsilon}\right)  \label{K2}%
\end{equation}
where $F^{\varepsilon}\in W_{-\sigma}^{1}\left(  \Pi_{+}^{\varepsilon}\right)
^{\ast}$ is an (anti)linear continuous functional on $W_{-\sigma}^{1}\left(
\Pi_{+}^{\varepsilon}\right)  $ and $\left(  \ ,\ \right)  _{\Pi^{\varepsilon
}}$ is an extension of the Lebesgue scalar product up to a duality between an
appropriate couple of weighted spaces. In view of (\ref{K1}) all terms in
(\ref{K2}) are defined correctly so that the problem (\ref{K2}) is associated
with the continuous mapping%
\[
W_{\sigma}^{1}\left(  \Pi_{+}^{\varepsilon}\right)  \ni u^{\varepsilon}%
\mapsto\mathcal{A}_{\sigma}^{\varepsilon}\left(  \lambda^{\varepsilon}\right)
u^{\varepsilon}=F^{\varepsilon}\in W_{-\sigma}^{1}\left(  \Pi_{+}%
^{\varepsilon}\right)  ^{\ast}.
\]
If $f^{\varepsilon}\in L_{\sigma}^{2}\left(  \Pi_{+}^{\varepsilon}\right)  $
that is $e^{\sigma z}f^{\varepsilon}\in L^{2}\left(  \Pi_{+}^{\varepsilon
}\right)  $, then the functional%
\[
v^{\varepsilon}\mapsto F^{\varepsilon}\left(  v^{\varepsilon}\right)  =\left(
f^{\varepsilon},v^{\varepsilon}\right)  _{\Pi_{+}^{\varepsilon}}%
\]
belongs to $W_{-\sigma}^{1}\left(  \Pi_{+}^{\varepsilon}\right)  ^{\ast}$.
Clearly, $\mathcal{A}_{-\sigma}^{\varepsilon}\left(  \lambda^{\varepsilon
}\right)  $ is the adjoint operator for $\mathcal{A}_{\sigma}^{\varepsilon
}\left(  \lambda^{\varepsilon}\right)  $.

\subsection{The Fredholm property, asymptotics and the index\label{sect5.2}}

Let us formulate some well-know results of the theory of elliptic problems in
domains with cylindrical outlets to infinity (see the key papers \cite{Ko,
MaPl1, MaPl2} and, e.g., the monographs \cite{NaPl, KoMaRo}). This theory
mainly deals with the classical (differential) formulation of boundary value
problems, however as was observed in \cite{na417}, passing to the weak
formulation involving integral identities of type (\ref{K2}) does not meet any
visible obstacle. The only disputable point, namely the dependence of
constants on the small parameter $\varepsilon$, we will be discussed in
Section \ref{sect5.5}.

\begin{theorem}
\label{TheoremKA}(see \cite{Ko}) Let $\lambda^{\varepsilon}\in\left(
0,\pi^{2}\right]  $.

1) The operator $\mathcal{A}_{\beta}^{\varepsilon}\left(  \lambda
^{\varepsilon}\right)  $ is Fredholm if and only if%
\begin{equation}
\beta\neq\beta_{0}:=0,\text{ \ \ }\beta\neq\beta_{\pm j}:=\pm\sqrt{\pi
^{2}j^{2}-\lambda^{\varepsilon}},\text{ \ \ }j\in\mathbb{N}. \label{K4}%
\end{equation}
In the case $\beta=\beta_{p}$ with $p\in\mathbb{Z}$ the range $\mathcal{A}%
_{\beta}^{\varepsilon}\left(  \lambda^{\varepsilon}\right)  W_{\beta}%
^{1}\left(  \Pi_{+}^{\varepsilon}\right)  $ is not closed subspace in
$W_{-\beta}^{1}\left(  \Pi_{+}^{\varepsilon}\right)  ^{\ast}.$

2) Let $\gamma\in\left(  \beta_{1},\beta_{2}\right)  $ and let $u^{\varepsilon
}\in W_{-\gamma}^{1}\left(  \Pi_{+}^{\varepsilon}\right)  $ be a solution of
problem (\ref{K2}) with the weight index $\sigma=-\gamma$ and the right-hand
side $F^{\varepsilon}\in W_{-\gamma}^{1}\left(  \Pi_{+}^{\varepsilon}\right)
^{\ast}\subset W_{\gamma}^{1}\left(  \Pi_{+}^{\varepsilon}\right)  $. Then the
asymptotic decomposition%
\begin{equation}
u^{\varepsilon}\left(  x\right)  =\widetilde{u}^{\varepsilon}\left(  x\right)
+\sum_{\pm}\left(  a_{\pm}^{\varepsilon}w_{0}^{\varepsilon\pm}\left(
x\right)  +b_{\pm}^{\varepsilon}w_{1}^{\varepsilon\pm}\left(  x\right)
\right)  \label{K5}%
\end{equation}
and the estimate%
\begin{equation}
\left(  \left\Vert \widetilde{u}^{\varepsilon};W_{\gamma}^{1}\left(
\Pi^{\varepsilon}\right)  \right\Vert ^{2}+\underset{\pm}{\sum}\left(
\left\vert a_{\pm}^{\varepsilon}\right\vert ^{2}+\left\vert b_{\pm
}^{\varepsilon}\right\vert ^{2}\right)  \right)  ^{1/2}\leq c_{\varepsilon
}\left(  \left\Vert F^{\varepsilon};W_{-\gamma}^{1}\left(  \Pi^{\varepsilon
}\right)  ^{\ast}\right\Vert +\left\Vert u^{\varepsilon};W_{-\gamma}%
^{1}\left(  \Pi^{\varepsilon}\right)  \right\Vert \right)  \label{K6}%
\end{equation}
are valid, where $\widetilde{u}^{\varepsilon}\in W_{\gamma}^{1}\left(
\Pi^{\varepsilon}\right)  $ is the asymptotic remainder, $a_{\pm}%
^{\varepsilon}$ and $b_{\pm}^{\varepsilon}$ are coefficients depending on
$F^{\varepsilon}$ and $u^{\varepsilon}$, the waves $w_{0}^{\varepsilon\pm}$
are given by (\ref{13}) and $w_{1}^{\varepsilon\pm}$ by (\ref{20}) for
$\lambda^{\varepsilon}=\pi^{2}$ but by (\ref{30}) for $\lambda^{\varepsilon
}\in\left(  0,\pi^{2}\right)  $. The factor $c_{\varepsilon}$ in (\ref{K6}) is
independent of $F^{\varepsilon}$ and $u^{\varepsilon}$ but may depend on
$\varepsilon\in\left[  0,\varepsilon_{0}\right]  .$
\end{theorem}

\bigskip

As was mentioned $\mathcal{A}_{\gamma}^{\varepsilon}\left(  \lambda
^{\varepsilon}\right)  ^{\ast}=\mathcal{A}_{-\gamma}^{\varepsilon}\left(
\lambda^{\varepsilon}\right)  $ and, hence, kernels and co-kernels of these
operators are in the relationship%
\begin{equation}
\ker\mathcal{A}_{\pm\gamma}^{\varepsilon}\left(  \lambda^{\varepsilon}\right)
=\text{\textrm{coker}}\mathcal{A}_{\mp\gamma}^{\varepsilon}\left(
\lambda^{\varepsilon}\right)  \label{K7}%
\end{equation}
In the next assertion we compare the indexes Ind$\mathcal{A}_{\pm\gamma
}^{\varepsilon}\left(  \lambda^{\varepsilon}\right)  =\dim\ker\mathcal{A}%
_{\pm\gamma}^{\varepsilon}\left(  \lambda^{\varepsilon}\right)  -\dim
$\textrm{coker}$\mathcal{A}_{\pm\gamma}^{\varepsilon}\left(  \lambda
^{\varepsilon}\right)  ;$ notice that Ind$\mathcal{A}_{\gamma}^{\varepsilon
}\left(  \lambda^{\varepsilon}\right)  =-$Ind$\mathcal{A}_{-\gamma
}^{\varepsilon}\left(  \lambda^{\varepsilon}\right)  $ according to (\ref{K7}).

\begin{theorem}
\label{TheoremKB} (see \cite[Thm. 3.3.3, 5.1.4 (4)]{NaPl}) If $\gamma
\in\left(  \beta_{1},\beta_{2}\right)  ,$ see (\ref{K4}), then%
\begin{equation}
\text{Ind}\mathcal{A}_{-\gamma}^{\varepsilon}\left(  \lambda^{\varepsilon
}\right)  =\text{Ind}\mathcal{A}_{\gamma}^{\varepsilon}\left(  \lambda
^{\varepsilon}\right)  +4 \label{K8}%
\end{equation}

\end{theorem}

We emphasize that the last $4$ is nothing but the number of waves detached in
(\ref{K5}). From (\ref{K7})-(\ref{K8}), it follows that
\begin{equation}
\text{Ind}\mathcal{A}_{-\gamma}^{\varepsilon}\left(  \lambda^{\varepsilon
}\right)  =-\text{Ind}\mathcal{A}_{\gamma}^{\varepsilon}\left(  \lambda
^{\varepsilon}\right)  =2 \label{K9}%
\end{equation}

\subsection{Absence of trapped modes with a fast decay rate.\label{sect5.3}}

In this section we prove that, for $\lambda\in(0,\pi^{2}]$ and $\gamma
\in\left(  \beta_{1},\beta_{2}\right)  $, there holds the formula%
\begin{equation}
\dim\ker\mathcal{A}_{\gamma}^{\varepsilon}\left(  \lambda\right)  =0
\label{M1}%
\end{equation}
which, in particular, completes the proof of Theorem \ref{TheoremA}, cf. our
assumption $c_{\varepsilon}\neq0$ for trapped mode (\ref{35}) while
$c_{\varepsilon}=0$ leads to $U^{\varepsilon}=\widetilde{U}^{\varepsilon}%
\in\ker\mathcal{A}_{\gamma}^{\varepsilon}\left(  \lambda\right)  $. Clearly,
$\ker\mathcal{A}_{\beta}^{0}\left(  \lambda\right)  =0$ for any $\beta>0$ that
is the limit problem (\ref{9}), (\ref{10}) in the straight semi-strip $\Pi
_{+}^{0}$ cannot have a trapped mode. However, as was mentioned in Section
\ref{sect1.4}, formula (\ref{M1}) does not follow by a standard perturbation
argument and, moreover, $\dim\ker\mathcal{A}_{\beta}^{\varepsilon}\left(
\lambda^{\varepsilon}\right)  >0$ for some $\beta\in\left(  0,\beta
_{1}\right)  $ and $\lambda^{\varepsilon}\in(0,\pi^{2})$.

\begin{theorem}
\label{TheoremMA}Let $\gamma\in\left(  \beta_{1},\beta_{2}\right)  $ be fixed.
There exists $\varepsilon_{0}>0$ such that, for $\varepsilon\in\left(
0,\varepsilon_{0}\right)  $ and $\lambda\in(0,\pi^{2}]$, formula (\ref{M1}) is valid.
\end{theorem}

\textbf{Proof.} Let us assume that, for some $\lambda\in\left(  0,\pi
^{2}\right]  $ and an infinitesimal positive sequence $\left\{  \varepsilon
_{k}\right\}  _{k\in\mathbb{N}}$, the homogeneous problem (\ref{9}),
(\ref{10}) has a solution $u^{\varepsilon_{k}}\in W_{\gamma}^{1}\left(
\Pi_{+}^{\varepsilon}\right)  $. We denote by $u_{0}^{\varepsilon_{k}}$ the
restriction of $u^{\varepsilon_{k}}$ onto the semi-strip $\Pi_{+}%
^{0}=\mathbb{R\times}\left(  0,1\right)  $. Under the normalization condition%
\begin{equation}
\left\Vert u^{\varepsilon_{k}};L^{2}\left(  \Pi_{+}^{\varepsilon}\left(
2l\right)  \right)  \right\Vert =1, \label{M0}%
\end{equation}
we are going to perform the limit passage $\varepsilon_{k}\rightarrow+0$ in
the integral identity%
\begin{equation}
\left(  \nabla u^{\varepsilon_{k}},\nabla v^{\varepsilon}\right)  _{\Pi
_{+}^{\varepsilon}}=\lambda^{\varepsilon}\left(  u^{\varepsilon}%
,v^{\varepsilon}\right)  _{\Pi_{+}^{\varepsilon}}, \label{M2}%
\end{equation}
where $v^{\varepsilon}$ is obtained from a test function $v\in C_{c}^{\infty
}\left(  \overline{\Pi_{+}^{0}}\right)  $ by the even extension over the
$x_{1}$-axis. If we prove that

$(i)$ $u_{0}^{\varepsilon_{k}}$ converges to $u_{0}^{0}\in W_{-\gamma}%
^{1}\left(  \Pi_{+}^{0}\right)  $ weakly in $W_{-\gamma}^{1}\left(  \Pi
_{+}^{0}\right)  $ and, therefore, strongly in $L^{2}\left(  \Pi_{+}%
^{0}\left(  2l\right)  \right)  ;$

$(ii)$ $\left\Vert u^{\varepsilon_{k}};L^{2}\left(  \Pi_{+}^{\varepsilon_{k}%
}\setminus\Pi_{+}^{0}\right)  \right\Vert \rightarrow0;$

$(iii)$ $\left(  \nabla u^{\varepsilon_{k}},\nabla v\right)  _{\Pi
^{\varepsilon}\setminus\Pi^{0}}\rightarrow0$ with any smooth function $v$ in
the rectangle $\left[  0,l\right]  \times\left[  -1,1\right]  ,$

\noindent then the limit passage in (\ref{M2}) and (\ref{M0}) gives
\begin{gather}
\left(  \nabla u_{0}^{0},\nabla v\right)  _{\Pi_{+}^{0}}=\lambda\left(
u_{0}^{0},v\right)  _{\Pi_{+}^{0}}\text{ \ \ }\forall v\in C_{c}^{\infty
}(\overline{\Pi_{+}^{0}}),\label{M3}\\
\left\Vert u_{0}^{0};L^{2}\left(  \Pi_{+}^{0}\left(  2l\right)  \right)
\right\Vert =1. \label{M4}%
\end{gather}
By a density argument, the integral identity (\ref{M3}) is valid with any
$v\in W_{-\gamma}^{1}\left(  \Pi_{+}^{0}\right)  $ and, therefore, $u_{0}%
^{0}=0$ because the limit problem in $\Pi_{+}^{0}$ cannot get a non-trivial
trapped mode.

Let us confirm facts $(i)-(iii)$. We write $\varepsilon$ instead of
$\varepsilon_{k}.$

First, we apply a local estimate, see, e.g., \cite{ADN1}, to the solution
$u^{\varepsilon}$ of the problem (\ref{9}), (\ref{10}) with $\lambda
u^{\varepsilon}$ as a given right-hand side:%
\begin{equation}
\left\Vert u^{\varepsilon};H^{2}\left(  \varpi^{\prime}\right)  \right\Vert
\leq c\lambda\left\Vert u^{\varepsilon};L^{2}\left(  \varpi^{\prime\prime
}\right)  \right\Vert . \label{M5}%
\end{equation}
Here, $\varpi^{\prime}=\left(  4l/3,5l/3\right)  \times\left(  0,1\right)  $
and $\varpi^{\prime\prime}=\left(  l,2l\right)  \times\left(  0,1\right)  $
are rectangles such that $\varpi^{\prime}\subset\varpi^{\prime\prime}%
\subset\Pi^{\varepsilon}\left(  2l\right)  $ and, therefore, the right-hand
side of (\ref{M5}) is less than $c\lambda$ according to (\ref{M0}).

Second, we split $u^{\varepsilon}$ as follows:%
\begin{equation}
u^{\varepsilon}=u_{l}^{\varepsilon}+u_{\infty}^{\varepsilon},\text{ \ \ }%
u_{l}^{\varepsilon}=\left(  1-\chi\right)  u^{\varepsilon},\text{
\ \ }u_{\infty}^{\varepsilon}=\chi u^{\varepsilon} \label{M55}%
\end{equation}
where $\chi\in C^{\infty}\left(  \mathbb{R}\right)  $ is a cut-off function,
$\chi\left(  x_{1}\right)  =1$ for $x_{1}\geq5l/3$ and $\chi\left(
x_{1}\right)  =0$ for $x_{1}\leq4l/3.$ The components in (\ref{M55}) satisfy
the integral identities%
\begin{gather}
\left(  \nabla u_{l}^{\varepsilon},\nabla v_{l}\right)  _{\Pi_{+}%
^{\varepsilon}\left(  2l\right)  }=\lambda\left(  \left(  1-\chi\right)
u^{\varepsilon},v_{l}\right)  _{\Pi_{+}^{\varepsilon}\left(  2l\right)
}+\left(  \nabla u^{\varepsilon},v_{l}\nabla\chi\right)  _{\varpi^{\prime}%
}-\label{M6}\\
-\left(  u^{\varepsilon}\nabla\chi,\nabla v_{l}\right)  _{\varpi^{\prime}%
}\text{ }\forall v_{l}\in H^{1}\left(  \Pi_{+}^{\varepsilon}\left(  2l\right)
\right)  ,\nonumber\\
\left(  \nabla u_{\infty}^{\varepsilon},\nabla v_{\infty}\right)
_{\Pi_{\infty}\left(  l\right)  }-\lambda\left(  u_{\infty}^{\varepsilon
},v_{\infty}\right)  _{\Pi_{\infty}\left(  l\right)  }=\label{M7}\\
=F_{\infty}^{\varepsilon}\left(  v_{\infty}\right)  :=\left(  u^{\varepsilon
}\nabla\chi,\nabla v_{\infty}\right)  _{\varpi^{\prime}}-\left(  \nabla
u^{\varepsilon},v_{\infty}\nabla\chi\right)  _{\varpi^{\prime}}\text{
\ }\forall v_{\infty}\in W_{-\gamma}^{1}\left(  \Pi_{\infty}\left(  l\right)
\right)  .\nonumber
\end{gather}
Third, inserting $v_{l}=u_{l}^{\varepsilon}$ into (\ref{M6}) and taking
(\ref{M0}), (\ref{M5}) into account yield%
\begin{equation}
\left\Vert \nabla u_{l}^{\varepsilon};L^{2}\left(  \Pi_{+}^{\varepsilon
}\left(  2l\right)  \right)  \right\Vert \leq c. \label{M8}%
\end{equation}
The problem (\ref{M7}) needs a bit more advanced argument. It is posed in the
semi-strip $\Pi_{\infty}\left(  l\right)  $ independent of $\varepsilon$ and,
thus, the following a priory estimate in the Kondratiev space, see \cite{Ko}
and, e.g.,\cite[Thm 5.1.4 (1)]{NaPl},%
\begin{align}
\left\Vert u_{\infty}^{\varepsilon};W_{-\gamma}^{1}\left(  \Pi_{\infty}\left(
l\right)  \right)  \right\Vert  &  \leq c_{1}\left(  \left\Vert F_{\infty
}^{\varepsilon};W_{\gamma}^{1}\left(  \Pi_{\infty}\left(  l\right)  \right)
^{\ast}\right\Vert +\left\Vert u_{\infty}^{\varepsilon};L^{2}\left(  \Pi
_{+}^{\varepsilon}\left(  2l\right)  \cap\Pi_{\infty}\left(  l\right)
\right)  \right\Vert \right) \label{M9}\\
&  \leq c_{2}\left(  \left\Vert u^{\varepsilon};L^{2}\left(  \varpi^{\prime
}\right)  \right\Vert +\left\Vert \nabla u^{\varepsilon};L^{2}\left(
\varpi^{\prime}\right)  \right\Vert +\left\Vert u^{\varepsilon};\Pi
_{+}^{\varepsilon}\left(  2l\right)  \right\Vert \right) \nonumber
\end{align}
involves some constants $c_{m}$ independent of $\varepsilon$. In this way,
formulas (\ref{M8}) and (\ref{M9}), (\ref{M5}), (\ref{M1}) assure that%
\begin{equation}
\left\Vert u^{\varepsilon};W_{-\gamma}^{1}\left(  \Pi_{+}^{\varepsilon
}\right)  \right\Vert \leq c. \label{M10}%
\end{equation}
Thus, the convergence in $\left(  i\right)  $ occurs along a subsequence which
is still denoted by $\left\{  \varepsilon_{k}\right\}  $.

The last step of our consideration uses integration in $t\in\left(
-\varepsilon,0\right)  $ and $x_{1}\in\left(  0,l\right)  $ of the
Newton-Leibnitz formula%
\[
\left\vert u^{\varepsilon}\left(  t,x_{2}\right)  \right\vert ^{2}=\int
_{t}^{t+1/2}\frac{\partial}{\partial x_{2}}\left(  \chi_{0}\left(
x_{2}\right)  \left\vert u^{\varepsilon}\left(  x_{1},x_{2}\right)
\right\vert ^{2}\right)  dx_{2}%
\]
where $\chi_{0}\in C^{\infty}\left(  \mathbb{R}\right)  $ is a cut-off
function, $\chi_{0}\left(  x_{2}\right)  =1$ for $x_{2}<1/6$ and $\chi
_{0}\left(  x_{2}\right)  =0$ for $x_{2}>1/3$. As a result, we obtain the
estimate%
\[
\int_{\Pi_{+}^{\varepsilon}\setminus\Pi_{+}^{0}}\left\vert u^{\varepsilon
}\left(  x\right)  \right\vert ^{2}dx\leq c\varepsilon\int_{\Pi_{+}%
^{\varepsilon}\left(  l\right)  }\left(  \left\vert \nabla u^{\varepsilon
}\left(  x\right)  \right\vert ^{2}+\left\vert u^{\varepsilon}\left(
x\right)  \right\vert ^{2}\right)  dx\leq C\varepsilon
\]
while referring to (\ref{M10}) again. This provides (ii) as well as (iii)
because
\begin{align*}
\left\vert \int_{0}^{l}\int_{-\varepsilon}^{0}\nabla u^{\varepsilon}\left(
x_{1},x_{2}\right)  \cdot\nabla v\left(  x_{1},x_{2}\right)  dx_{2}%
dx_{1}\right\vert  &  \leq\underset{x\in\overline{\Pi_{+}^{\varepsilon}%
}\setminus\Pi_{+}^{0}}{\max}\left\vert \nabla v\left(  x\right)  \right\vert
\left(  \text{meas}_{2}\left(  \Pi_{+}^{\varepsilon}\setminus\Pi^{0}\right)
\right)  ^{1/2}\left\Vert \nabla u^{\varepsilon};L^{2}\left(  \Pi
_{+}^{\varepsilon}\setminus\Pi_{+}^{0}\right)  \right\Vert \leq\\
&  \leq c_{v}\varepsilon^{1/2}l^{1/2}\left\Vert u^{\varepsilon};W_{\gamma}%
^{1}\left(  \Pi_{+}^{\varepsilon}\right)  \right\Vert \leq C_{v}%
\varepsilon^{1/2}%
\end{align*}
Theorem \ref{TheoremMA} is proved. $\boxtimes$

\subsection{Radiation conditions\label{sect5.4}}

Let $\lambda^{\varepsilon}\in\left(  0,\pi^{2}\right)  $ and $\gamma\in\left(
\beta_{1},\beta_{2}\right)  $, cf. Theorem \ref{TheoremKA}. The pre-image
$\mathfrak{W}_{\gamma}^{1}\left(  \Pi_{+}^{\varepsilon}\right)  $ of the
subspace $W_{-\gamma}^{1}\left(  \Pi_{+}^{\varepsilon}\right)  ^{\ast}$ in
$W_{\gamma}^{1}\left(  \Pi_{+}^{\varepsilon}\right)  ^{\ast}$ for the operator
$\mathcal{A}_{-\gamma}^{\varepsilon}\left(  \lambda^{\varepsilon}\right)  $
consists of functions in the form (\ref{K5}). Introducing the norm $\left\Vert
u^{\varepsilon};\mathfrak{W}_{\gamma}^{1}\left(  \Pi_{+}^{\varepsilon}\right)
\right\Vert $ as the left-hand side of (\ref{K6}) makes $\mathfrak{W}_{\gamma
}^{1}\left(  \Pi_{+}^{\varepsilon}\right)  $ a Hilbert space but this Hilbert
structure is of no use in our paper.

The restriction $\mathcal{B}_{\gamma}^{\varepsilon}\left(  \lambda
^{\varepsilon}\right)  $ of $\mathcal{A}_{-\gamma}^{\varepsilon}\left(
\lambda^{\varepsilon}\right)  $ onto $\mathcal{W}_{\gamma}^{1}\left(  \Pi
_{+}^{\varepsilon}\right)  \subset W_{-\gamma}^{1}\left(  \Pi_{+}%
^{\varepsilon}\right)  $ inherits all properties of $\mathcal{A}_{-\gamma
}^{\varepsilon}\left(  \lambda^{\varepsilon}\right)  ,$ in particular,
$\mathcal{B}_{\gamma}^{\varepsilon}\left(  \lambda^{\varepsilon}\right)  $ is
a Fredholm operator with Ind$\mathcal{B}_{\gamma}^{\varepsilon}\left(
\lambda^{\varepsilon}\right)  =$Ind$\mathcal{A}_{-\gamma}^{\varepsilon}\left(
\lambda^{\varepsilon}\right)  =2,$ see (\ref{K9}). Thus, the restriction
$\mathfrak{A}_{\gamma}^{\varepsilon}\left(  \lambda^{\varepsilon}\right)
_{out}$ of $\mathfrak{A}_{\gamma}^{\varepsilon}\left(  \lambda^{\varepsilon
}\right)  $ onto the subspace%
\begin{equation}
\mathcal{W}_{\gamma}^{1}\left(  \Pi_{+}^{\varepsilon}\right)  _{out}=\left\{
u^{\varepsilon}\in\mathcal{W}_{\gamma}^{1}\left(  \Pi_{+}^{\varepsilon
}\right)  :a_{-}^{\varepsilon}=b_{-}^{\varepsilon}=0\text{ in (\ref{K5}%
)}\right\}  \label{K11}%
\end{equation}
of codimension $2$ becomes of index zero.

\begin{theorem}
\label{TheoremKMM}Let $\lambda$, $\gamma$ and $\varepsilon$ be the same as in
Theorem \ref{TheoremMA}. Then the operator $\mathcal{B}_{\gamma}^{\varepsilon
}\left(  \lambda^{\varepsilon}\right)  $ actualizes the isomorphism%
\[
\mathcal{W}_{\gamma}^{1}\left(  \Pi_{+}^{\varepsilon}\right)  _{out}\approx
W_{-\gamma}^{1}\left(  \Pi_{+}^{\varepsilon}\right)  ^{\ast}.
\]

\end{theorem}

Since the decomposition (\ref{K5}) of a function $u^{\varepsilon}%
\in\mathfrak{W}_{\gamma}^{1}\left(  \Pi_{+}^{\varepsilon}\right)  _{out}$
loses the incoming waves $w_{0}^{\varepsilon-}$ and $w_{1}^{\varepsilon-}$ due
to the restriction in (\ref{K11}), $\mathcal{B}_{\gamma}^{\varepsilon}\left(
\lambda^{\varepsilon}\right)  $ has to be interpreted as an operator of the
problem (\ref{K11}) with the radiation condition (\ref{23}) at $\lambda
=\pi^{2}$ and (\ref{32}) at $\lambda\in(0,\pi^{2}).$ Theorem \ref{TheoremKMM}
says that such problem is uniquely solvable, while its solution in the form%
\begin{equation}
u^{\varepsilon}\left(  x\right)  =\widetilde{u}^{\varepsilon}\left(  x\right)
+a_{+}^{\varepsilon}w_{0}^{\varepsilon+}\left(  x\right)  +b_{+}^{\varepsilon
}w_{1}^{\varepsilon+}\left(  x\right)  \label{KK1}%
\end{equation}
obeys the estimate%
\begin{equation}
\left\Vert \widetilde{u}^{\varepsilon};W_{\gamma}^{1}\left(  \Pi
_{+}^{\varepsilon}\right)  \right\Vert +\left\vert a_{+}^{\varepsilon
}\right\vert +\left\vert b_{+}^{\varepsilon}\right\vert \leq C_{\varepsilon
}\left\Vert F^{\varepsilon};W_{-\gamma}^{1}\left(  \Pi_{+}^{\varepsilon
}\right)  \right\Vert . \label{KK2}%
\end{equation}

\subsection{Remark on the dependence of bounds on the small parameter
$\varepsilon$\label{sect5.5}}

If
\begin{equation}
\lambda^{\varepsilon}\in\left[  \delta,\pi^{2}-\delta\right]  \label{KK3}%
\end{equation}
with a fixed $\delta>0,$ the coefficient in the estimate (\ref{K6}) can be
chosen independent of $\varepsilon\in\left[  0,\varepsilon\left(
\delta\right)  \right]  $ with some $\varepsilon\left(  \delta\right)  >0$.
This fact originates in the smooth dependence of the waves (\ref{13}) and
(\ref{27}), (\ref{30}) on the parameter (\ref{KK3}) and the following
consideration. By multiplying $u^{\varepsilon}$ with the same cut-off function
$\chi$ as in (\ref{M55}), we reduce the problem (\ref{K2}) onto the semi-strip
$\Pi_{\infty}\left(  l\right)  $, namely, inserting $v^{\varepsilon}=\chi
v_{\infty}$ with any $v_{\infty}\in W_{\gamma}^{1}\left(  \Pi_{\infty}\left(
l\right)  \right)  $ as a test function, we obtain for $u_{\infty
}^{\varepsilon}=\chi u^{\varepsilon}$ the integral identity%
\begin{gather}
\left(  \nabla u_{\infty}^{\varepsilon},\nabla v_{\infty}\right)
_{\Pi_{\infty}\left(  l\right)  }-\lambda^{\varepsilon}\left(  u_{\infty
}^{\varepsilon},v_{\infty}\right)  _{\Pi_{\infty}\left(  l\right)  }%
=F_{\infty}^{\varepsilon}\left(  v_{\infty}\right)
:=\ \ \ \ \ \ \ \ \ \ \ \ \ \ \ \ \ \ \ \ \ \ \ \ \ \ \label{KK4}\\
:=F^{\varepsilon}\left(  \chi v_{\infty}\right)  -\left(  \nabla
u^{\varepsilon},v_{\infty}\nabla\chi\right)  _{\Pi_{\infty}\left(  l\right)
}+\left(  u^{\varepsilon}\nabla\chi,\nabla v_{\infty}\right)  _{\Pi_{\infty
}\left(  l\right)  }.\nonumber
\end{gather}
Moreover,%
\begin{align*}
\left\Vert F_{\infty}^{\varepsilon};W_{-\gamma}^{1}\left(  \Pi_{\infty}\left(
l\right)  \right)  ^{\ast}\right\Vert  &  \leq\left\Vert F^{\varepsilon
}\left(  \chi\cdot\right)  ;W_{-\gamma}^{1}\left(  \Pi_{\infty}\left(
l\right)  \right)  \right\Vert +c_{\chi}\left\Vert u^{\varepsilon}%
;H^{1}\left(  \Pi_{\infty}\left(  l\right)  \cap\Pi_{+}^{\varepsilon}\left(
2l\right)  \right)  \right\Vert \leq\\
&  \leq c\left(  \left\Vert F^{\varepsilon};W_{-\gamma}^{1}\left(  \Pi
_{+}^{\varepsilon}\right)  \right\Vert +\left\Vert u^{\varepsilon};W_{-\gamma
}^{1}\left(  \Pi_{+}^{\varepsilon}\right)  \right\Vert \right)  ,
\end{align*}
cf. the right-hand side of (\ref{K6}). By restriction (\ref{KK3}),
$\lambda^{\varepsilon}$ stays at a distance from the thresholds $\lambda
_{0}^{\dagger}=0$ and $\lambda_{1}^{\dagger}=\pi^{2}$ so that we may choose
the same weight index $\gamma$ for all legalized $\lambda^{\varepsilon}$.

Hence, a general result in \cite{Ko}, see also \cite[\S \ 3.2]{NaPl}, on the
basis of a perturbation argument provides a common factor $c^{\varepsilon
}=const$ in the estimate (\ref{K6}) for ingredients of the asymptotic
representation (\ref{K5}) of the solution $u_{\infty}^{\varepsilon}=\chi
u^{\varepsilon}$ to the problem (\ref{KK4}) in the $\varepsilon$-independent
domain $\Pi_{\infty}\left(  l\right)  .$ Since the weight $e^{\gamma x_{1}}$
is uniformly bounded in $\Pi_{+}^{\varepsilon}\left(  2l\right)  =\Pi
_{+}^{\varepsilon}\setminus\Pi_{\infty}\left(  2l\right)  $, the evident
relation%
\[
\left\Vert \widetilde{u}^{\varepsilon};W_{\gamma}^{1}\left(  \Pi
_{+}^{\varepsilon}\left(  2l\right)  \right)  \right\Vert \leq c\left\Vert
u^{\varepsilon};W_{-\gamma}^{1}\left(  \Pi_{+}^{\varepsilon}\left(  2l\right)
\right)  \right\Vert +\sum_{\pm}\left(  \left\vert a_{\pm}^{\varepsilon
}\right\vert +\left\vert b_{\pm}^{\varepsilon}\right\vert \right)
\]
allows us to extend the above mentioned estimate over the whole waveguide
$\Pi_{+}^{\varepsilon}$.

Similarly, in the case (\ref{KK3}) the factor $C^{\varepsilon}$ in (\ref{KK2})
can be fixed independent of $\varepsilon$, too.

The desired eigenvalue (\ref{42}) is located in the vicinity of the threshold
$\lambda_{1}^{\dagger}=\pi^{2}$ and the above consideration becomes
unacceptable. Moreover, the normalization factor $\left(  \pi^{2}%
-\lambda^{\varepsilon}\right)  ^{-1/4}$ in (\ref{27}) is big so that the
independence property of $c^{\varepsilon}$ and $C^{\varepsilon}$ is surely
lost. Thus, our immediate objective is to modify the estimates in order to
make them homotype for all small $\varepsilon>0$. We emphasize that a
modification of the normalization factor does not suffice because the waves
$e^{\pm k_{1}^{\varepsilon}x_{1}}\cos\left(  \pi x_{2}\right)  $ in (\ref{27})
become equal at $\varepsilon=0$.

We follow a scheme in \cite[\S 3]{na489} and define for $\lambda^{\varepsilon
}\in((0,\pi^{2})$ the linear combinations of the exponential waves (\ref{27})%
\begin{equation}
\mathbf{w}_{1}^{\pm}\left(  \lambda^{\varepsilon};x\right)  =(1/2)\cos\left(
\pi x_{2}\right)  ((1/k_{1}^{\varepsilon})(e^{k_{1}^{\varepsilon}x_{1}%
}-e^{-k_{1}^{\varepsilon}x_{1}})\mp i(e^{k_{1}^{\varepsilon}x_{1}}%
+e^{-k_{1}^{\varepsilon}x_{1}})), \label{W1}%
\end{equation}
cf. (\ref{30}). A direct calculation demonstrates that the new waves
(\ref{W1}) together with the old waves (\ref{13}),%
\begin{equation}
\mathbf{w}_{0}^{\pm}\left(  \lambda^{\varepsilon};x\right)  =w_{0}%
^{\varepsilon\pm}\left(  x\right)  =\left(  2k\right)  ^{-1/2}e^{\pm
ik^{\varepsilon}x_{1}}, \label{W2}%
\end{equation}
still satisfy the normalization and orthogonality conditions (\ref{19}) but
additionally are in the relationship%
\[
\mathbf{w}_{0}^{\pm}\left(  \lambda^{\varepsilon};x\right)  -w_{1}^{0\pm
}\left(  x\right)  =O((\pi^{2}-\lambda^{\varepsilon})x_{1}),\ \ \ \ \mathbf{w}%
_{1}^{\pm}\left(  \lambda^{\varepsilon};x\right)  -w_{0}^{0\pm}\left(
x\right)  =O((\pi^{2}-\lambda^{\varepsilon})^{1/2}x_{1}).
\]
In other words, the waves (\ref{W1}) and (\ref{W2}) smoothly become the waves
(\ref{20}) and (\ref{13}) introduced in Section \ref{sect2.1} at the threshold
$\lambda^{\varepsilon}=\pi^{2}$. The first property of $\mathbf{w}_{p}^{\pm
}\left(  \lambda^{\varepsilon};x\right)  $ allows us to repeat considerations
in Sections \ref{sect5.4}, \ref{sect2.2} and compose the space $\mathbf{W}%
_{\gamma}^{1}\left(  \Pi_{+}^{\varepsilon}\right)  _{out}$ of functions
satisfying the new, so-called artificial, radiation condition%

\begin{equation}
\mathbf{u}^{\varepsilon}\left(  x\right)  =\widetilde{\mathbf{u}}%
^{\varepsilon}\left(  x\right)  +\mathbf{a}_{+}^{\varepsilon}\mathbf{w}%
_{0}^{+}\left(  \lambda^{\varepsilon};x\right)  +\mathbf{b}_{+}^{\varepsilon
}\mathbf{w}_{1}^{+}\left(  \lambda^{\varepsilon};x\right)  ,\text{
\ \ }\widetilde{\mathbf{u}}^{\varepsilon}\in\mathbf{W}_{\gamma}^{1}\left(
\Pi_{+}^{\varepsilon}\right)  \label{U5}%
\end{equation}
cf. (\ref{KK1}), to determine the solutions $\mathbf{Z}_{p}^{\varepsilon
}\left(  \lambda^{\varepsilon};\cdot\right)  \in\mathbf{W}_{-\gamma}%
^{1}\left(  \Pi_{+}^{\varepsilon}\right)  $ of the homogeneous problem
(\ref{K2}), $\sigma=-\gamma,$%
\begin{gather}
\mathbf{Z}_{p}^{\varepsilon}\left(  \lambda^{\varepsilon};x\right)
=\widetilde{\mathbf{Z}}_{p}^{\varepsilon}\left(  \lambda^{\varepsilon
};x\right)  +\mathbf{w}_{p}^{-}\left(  \lambda^{\varepsilon};x\right)
+\mathbf{S}_{0p}^{\varepsilon}\left(  \lambda^{\varepsilon}\right)
\mathbf{w}_{0}^{-}\left(  \lambda^{\varepsilon};x\right)  +\mathbf{S}%
_{1p}^{\varepsilon}\left(  \lambda^{\varepsilon}\right)  \mathbf{w}_{1}%
^{+}\left(  \lambda^{\varepsilon};x\right)  ,\text{ }\label{U6}\\
\widetilde{\mathbf{Z}}_{p}^{\varepsilon}\left(  \lambda^{\varepsilon}%
;\cdot\right)  \in\mathbf{W}_{\gamma}^{1}\left(  \Pi_{+}^{\varepsilon}\right)
,\text{ }p=0,1,\nonumber
\end{gather}
cf. (\ref{31}) and to detect a unitary and symmetric artificial scattering
matrix $\mathbf{S}^{\varepsilon}\left(  \lambda^{\varepsilon}\right)  =\left(
\mathbf{S}_{qp}^{\varepsilon}\left(  \lambda^{\varepsilon}\right)  \right)
_{q,\text{ }p=0,1}$. At the same time, the second property of $\mathbf{w}%
_{p}^{\pm}\left(  \lambda^{\varepsilon};x\right)  $ assures that, for a fixed
$\varepsilon$, the operator%
\begin{equation}
\mathbf{B}_{\gamma}^{\varepsilon}\left(  \lambda^{\varepsilon}\right)
_{out}:\mathbf{W}_{\gamma}^{1}\left(  \Pi_{+}^{\varepsilon}\right)
_{out}\rightarrow W_{-\gamma}^{1}\left(  \Pi_{+}^{\varepsilon}\right)  ^{\ast}
\label{U7}%
\end{equation}
of the problem (\ref{K2}), $\sigma=-\gamma$, with the radiation condition
(\ref{U5}) depends continuously on the spectral parameter $\lambda
^{\varepsilon}\in\left(  \pi^{2}-\delta,\pi^{2}\right]  $, $\delta>0$, when
the domain of $\mathbf{B}_{\gamma}^{\varepsilon}\left(  \lambda^{\varepsilon
}\right)  _{out}$ is equipped with the norm
\begin{equation}
\left\Vert \mathbf{u}^{\varepsilon};\mathbf{W}_{\gamma}^{1}\left(  \Pi
_{+}^{\varepsilon}\right)  \right\Vert =\left\Vert \widetilde{\mathbf{u}%
}^{\varepsilon};W_{\gamma}^{1}\left(  \Pi_{+}^{\varepsilon}\right)
\right\Vert +\left\vert \mathbf{a}_{+}^{\varepsilon}\right\vert +\left\vert
\mathbf{b}_{+}^{\varepsilon}\right\vert \label{U10}%
\end{equation}
of a weighted space with detached asymptotics, cf. the left-hand side of
(\ref{KK3}).

Recalling our reasoning in Section \ref{sect5.3} and the beginning of this
section, we conclude that the operator (\ref{U7}) is an isomorphism while its
norm and the norm of the inverse are uniformly bounded in
\begin{equation}
\lambda\in\left[  \pi^{2}-\delta,\pi^{2}\right]  ,\text{ \ \ }\varepsilon
\in\left[  0,\varepsilon_{0}\right]  . \label{U8}%
\end{equation}
Furthermore, by the Fourier method, entries of the matrix $\mathbf{S}%
^{\varepsilon}\left(  \lambda\right)  $ can be expressed as weighted integrals
of solutions (\ref{U6}), this matrix is continuous in both arguments
(\ref{U8}) and the limit matrix%
\begin{equation}
\mathbf{S}^{0}\left(  \pi^{2}\right)  =\text{diag}\left\{  1,-1\right\}
\label{U9}%
\end{equation}
is nothing but augmented scattering matrix at the thresholds and its diagonal
form is due to the explicit solutions (\ref{491}) and (\ref{492}) in the
semi-strip $\Pi^{0}$,%
\[
Z_{0}^{0}\left(  x\right)  =(1/\sqrt{2\pi})\left(  e^{i\pi x_{1}}+e^{-i\pi
x_{1}}\right)  ,\ \ Z_{1}^{0}\left(  x_{1}\right)  =\cos\left(  \pi
x_{2}\right)  =(1/2i)\left(  \left(  x_{1}+i\right)  \cos\left(  \pi
x_{2}\right)  -\left(  x_{1}-i\right)  \cos\left(  \pi x_{2}\right)  \right)
.
\]

We resume that above consideration and find out a unique solution
$\mathbf{u}^{\varepsilon}\in\mathbf{W}_{\gamma}^{1}\left(  \Pi_{+}%
^{\varepsilon}\right)  _{out}\subset W_{-\gamma}^{1}\left(  \Pi_{+}%
^{\varepsilon}\right)  $ of the problem (\ref{K2}) with $\sigma=-\gamma$,
$F^{\varepsilon}\in W_{-\gamma}^{1}\left(  \Pi_{+}^{\varepsilon}\right)
^{\ast}$ and the artificial radiation condition (\ref{U5}). Moreover, the
estimate%
\begin{equation}
\left\Vert \mathbf{u}^{\varepsilon};\mathbf{W}_{\gamma}^{1}\left(  \Pi
_{+}^{\varepsilon}\right)  \right\Vert \leq c\left\Vert F^{\varepsilon
};W_{-\gamma}^{1}\left(  \Pi_{+}^{\varepsilon}\right)  ^{\ast}\right\Vert
\label{U12}%
\end{equation}
is valid, where $c$ is independent of both parameters (\ref{U8}).

We now search for a solution $u^{\varepsilon}\in W_{\gamma}^{1}\left(  \Pi
_{+}^{\varepsilon}\right)  _{out}$ of the same integral identity but the
radiation condition from Section \ref{sect5.4} in the form%
\begin{equation}
u^{\varepsilon}=\mathbf{u}^{\varepsilon}+\mathbf{c}_{0}^{\varepsilon
}\mathbf{Z}_{0}^{\varepsilon}+\mathbf{c}_{1}^{\varepsilon}\mathbf{Z}%
_{1}^{\varepsilon}. \label{U13}%
\end{equation}
The unknown coefficients $\mathbf{c}_{p}^{\varepsilon}$ should be fixed such
that the decomposition (\ref{KK1}) is satisfied. To this end, we insert into
the right-hand side of (\ref{U13}) formulas (\ref{U5}), (\ref{U6}) and
(\ref{W1}), (\ref{W2}), we compare the resultant coefficients of the waves
(\ref{13}), (\ref{27}) in (\ref{U13}) with those in (\ref{KK1}) and arrive at
the following systems of linear algebraic equations for the unknowns
$\mathbf{c}_{0}^{\varepsilon}$, $\mathbf{c}_{1}^{\varepsilon}$ and
$a_{+}^{\varepsilon}$, $b_{+}^{\varepsilon}:$%
\begin{align}
a_{+}^{\varepsilon}  &  =\mathbf{a}_{+}^{\varepsilon}+\mathbf{S}%
_{01}^{\varepsilon}\mathbf{c}_{1}^{\varepsilon},\text{ \ \ }0=\mathbf{c}%
_{0}^{\varepsilon},\label{U14}\\
\left(  2k_{1}^{\varepsilon}\right)  ^{1/2}b_{+}^{\varepsilon}  &  =\left(
1-ik_{1}^{\varepsilon}\right)  \mathbf{b}_{+}^{\varepsilon}+\left(  \left(
1+ik_{1}^{\varepsilon}\right)  +\left(  1-ik_{1}^{\varepsilon}\right)
\mathbf{S}_{11}^{\varepsilon}\right)  \mathbf{c}_{1}^{\varepsilon}%
,\label{U15}\\
\left(  2k_{1}^{\varepsilon}\right)  ^{1/2}b_{+}^{\varepsilon}  &  =\left(
1+ik_{1}^{\varepsilon}\right)  \mathbf{b}_{+}^{\varepsilon}+\left(  \left(
1-\varepsilon k_{1}^{\varepsilon}\right)  +\left(  1+ik_{1}^{\varepsilon
}\right)  \mathbf{S}_{11}^{\varepsilon}\right)  \mathbf{c}_{1}^{\varepsilon
}.\nonumber
\end{align}
Solving the system (\ref{U15}) with the help of the Cramer's rule, a simple
calculation gives the determinant%
\[
\left(  2k_{1}^{\varepsilon}\right)  ^{3/2}\left(  1-\mathbf{S}_{11}%
^{\varepsilon}\right)  =\left(  2\sqrt{\pi^{2}-\lambda^{\varepsilon}}\right)
^{3/2}i\left(  1-\mathbf{S}_{11}^{\varepsilon}\right)
\]
and the estimates%
\begin{equation}
\left\vert b_{+}^{\varepsilon}\right\vert \leq c\left(  \pi^{2}-\lambda
^{\varepsilon}\right)  ^{-1/2}\left\vert \mathbf{b}_{+}^{\varepsilon
}\right\vert ,\text{ \ \ \ }\left\vert \mathbf{c}_{1}^{\varepsilon}\right\vert
\leq c\left\vert \mathbf{b}_{+}^{\varepsilon}\right\vert \label{U16}%
\end{equation}
because $2\geq\left\vert 1-\mathbf{S}_{11}^{\varepsilon}\right\vert \geq1/2$
due to (\ref{U9}) and (\ref{U8}). In view of the first relation in (\ref{U14})
we obtain that%
\[
\left\vert a_{+}^{\varepsilon}\right\vert \leq c\left(  \left\vert
\mathbf{a}_{+}^{\varepsilon}\right\vert +\left\vert \mathbf{b}_{+}%
^{\varepsilon}\right\vert \right)  .
\]
Collecting formulas (\ref{U15}), (\ref{U16}) and (\ref{U12}), (\ref{U10})
adjusts the inequality (\ref{KK2}) as well as Theorem \ref{TheoremKMM}.

\begin{theorem}
\label{TheoremKUM} Let $\lambda^{\varepsilon}\in\left[  \pi^{2}-\delta,\pi
^{2}\right]  $, $\varepsilon\in\left(  0,\varepsilon_{0}\right]  $ and
$\gamma\in\left(  \beta_{1},\beta_{2}\right)  $. The solution (\ref{KK1}) of
the problem (\ref{K2}) with $\sigma=-\gamma$ and $F^{\varepsilon}\in
W_{-\gamma}^{1}\left(  \Pi_{+}^{\varepsilon}\right)  ^{\ast}$ admits the
estimate%
\begin{equation}
\left\Vert \widetilde{u}^{\varepsilon};W_{\gamma}^{1}\left(  \Pi
_{+}^{\varepsilon}\right)  \right\Vert +\left\vert a_{+}^{\varepsilon
}\right\vert +\left(  \pi^{2}-\lambda^{\varepsilon}\right)  ^{1/4}\left\vert
b_{+}^{\varepsilon}\right\vert \leq C\left\Vert F^{\varepsilon};W_{-\gamma
}^{1}\left(  \Pi_{+}^{\varepsilon}\right)  ^{\ast}\right\Vert \label{U18}%
\end{equation}
where $C$ does not depend on $\lambda^{\varepsilon}$, $\varepsilon$ and
$F^{\varepsilon}$.
\end{theorem}

\section{Justification of asymptotics\label{sect6}}

\subsection{The global asymptotic approximation\label{sect6.1}}

The reformulation (\ref{T2}) of the criterion (\ref{32}) implicates the
coefficients $S_{11}^{\varepsilon}$ and $S_{10}^{\varepsilon}$ in the
decomposition (\ref{31}) of the special solution $Z_{1}^{\varepsilon}$ of the
problem (\ref{1}), (\ref{2}) and this section is devoted to the justification
of the formal asymptotic expansions (\ref{44}). We emphasize that the similar
expansions (\ref{60}) of other entries in the augmented scattering matrix
$S^{\varepsilon}$ can be verified in the same way but actually we had used in
Section \ref{sect4.1} much simpler relation (\ref{T3}) only.

In Section \ref{sect3} we applied the method of matched asymptotic expansions
and our immediate objective becomes to compose a global approximation solution
from the inner and outer expansions (\ref{47}) and (\ref{Z1}). To this end, we
employ several smooth cut-off functions:%
\begin{align}
X_{\varepsilon}\left(  x\right)   &  =1\text{ for }x_{1}\leq l+1/\varepsilon
,\text{ \ \ }X_{\varepsilon}\left(  x\right)  =0\text{ for }x_{1}%
\geq2l+1/\varepsilon,\label{J1}\\
\chi_{\infty}\left(  x\right)   &  =1\text{ for }x_{1}\geq2l\text{
},\text{\ \ \ \ \ \ \ \ }\chi_{\infty}\left(  x\right)  =0\text{ for }%
x_{1}\leq3l/2,\nonumber\\
\chi_{\varepsilon}\left(  r\right)   &  =1\text{ for }r\leq2\varepsilon,\text{
\ \ \ \ \ \ \ \ \ }\chi_{\varepsilon}\left(  r\right)  =0\text{ for }%
r\geq3\varepsilon,\nonumber
\end{align}
where $r=\left(  \left\vert x_{1}-l\right\vert ^{2}+x_{2}^{2}\right)  ^{1/2}.$
We set%
\begin{align}
\mathfrak{Z}^{\varepsilon}  &  =\chi_{\infty}\mathfrak{Z}^{out}+X_{\varepsilon
}\mathfrak{Z}^{in}-\chi_{\infty}X_{\varepsilon}\mathfrak{Z}^{mat},\label{J2}\\
\mathfrak{Z}^{out}\left(  x\right)   &  =w_{1}^{\varepsilon-}\left(  x\right)
+S_{11}^{0}w_{1}^{\varepsilon+}\left(  x\right)  +\varepsilon^{1/2}S_{01}%
^{0}w_{0}^{\varepsilon+}\left(  x\right)  ,\label{J3}\\
\mathfrak{Z}^{in}\left(  x\right)   &  =\varepsilon^{-1/2}Z_{1}^{0}\left(
x\right)  +\varepsilon^{1/2}\left(  \left(  1-\chi_{\varepsilon}\left(
r\right)  \right)  \widehat{Z}_{1}^{\prime}\left(  x\right)  +\chi
_{\varepsilon}\left(  r\right)  Z_{1}^{\prime}\left(  l,0\right)  \right)
,\label{J4}\\
\mathfrak{Z}^{mat}\left(  x\right)   &  =\varepsilon^{-1/2}\left(
4\mu\right)  ^{-1/4}\cos\left(  \pi x_{2}\right)  \left(  1+i+S_{11}%
^{0}\left(  1-i\right)  +\varepsilon x_{1}\sqrt{\mu}\left(  1-i+S_{11}%
^{0}\left(  1+i\right)  \right)  \right)  +\label{J5}\\
&  +\varepsilon^{1/2}S_{01}^{0}\left(  2\pi\right)  ^{-1/2}e^{i\pi x_{1}%
}.\nonumber
\end{align}

This construction needs explanations. First, the expressions (\ref{J3}) and
(\ref{J4}) of the outer and inner types are multiplied with the cut-off
functions $\chi_{\infty}$ and $X_{\varepsilon}$ whose support overlap so that
the sum (\ref{J5}) of terms matched in Section \ref{sect3.2} attend the global
approximation (\ref{J2}) twice, i.e. in $\chi_{\infty}\mathfrak{Z}^{out}$ and
$X_{\varepsilon}\mathfrak{Z}^{in}$, but we compensate for this duplication by
subtracting $\chi_{\infty}X_{\varepsilon}\mathfrak{Z}^{mat}$. Moreover, the
formula for commutators $\left[  \Delta,\chi_{\infty}X_{\varepsilon}\right]
=\left[  \Delta,\chi_{\infty}\right]  +\left[  \Delta,X_{\varepsilon}\right]
$ demonstrates that
\begin{gather}
\left(  \Delta+\lambda^{\varepsilon}\right)  \mathfrak{Z}^{\varepsilon}%
=\chi_{\infty}\left(  \Delta+\lambda^{\varepsilon}\right)  \mathfrak{Z}%
^{out}+X_{\varepsilon}\left(  \Delta+\lambda^{\varepsilon}\right)
\mathfrak{Z}^{in}-\chi_{\infty}X_{\varepsilon}\left(  \Delta+\lambda
^{\varepsilon}\right)  \mathfrak{Z}^{mat}+\label{J6}\\
+\left[  \Delta,\chi_{\infty}\right]  \left(  \mathfrak{Z}^{out}%
-\mathfrak{Z}^{mat}\right)  +\left[  \Delta,X_{\varepsilon}\right]  \left(
\mathfrak{Z}^{in}-\mathfrak{Z}^{mat}\right)  :=\nonumber\\
:=\mathcal{F}^{\varepsilon}=\chi_{\infty}\mathcal{F}^{out}+X_{\varepsilon
}\mathcal{F}^{in}-\chi_{\infty}X_{\varepsilon}\mathcal{F}^{mat}+\mathcal{F}%
^{oma}+\mathcal{F}^{ima}.\nonumber
\end{gather}
Second, the function $Z_{1}^{0}$ is properly defined by the formula (\ref{50})
in the whole waveguide but $Z_{1}^{\prime}$ needs an extension from $\Pi
_{+}^{0}$ onto $\Pi_{+}^{\varepsilon}$ denoted by $\widehat{Z}_{1}^{\prime}$
in (\ref{J4}). Since the Neumann datum (\ref{52}) in the problem (\ref{51}),
$p=1$, has a jump at the point $(0,l)\in\partial\Pi_{+}^{0}$, the solution
$Z_{p}^{\prime}$ gets a singular behavior near this point. A simple
calculation based on the Kondratiev theory \cite{Ko} (see also \cite[Ch.
2]{NaPl}) demonstrates that
\begin{equation}
Z_{1}^{\prime}\left(  x\right)  =\pi^{-1}G_{1}^{\prime}r\left(  \ln
r\cos\varphi-\varphi\sin\varphi\right)  +\breve{Z}_{1}^{\prime}\left(
x\right)  \label{J7}%
\end{equation}
where $\left(  r,\varphi\right)  \in\mathbb{R}_{+}\times\left(  0,\pi\right)
$ are polar coordinates in fig. \ref{f4}, b and $\breve{Z}_{1}^{\prime}$ is a
smooth function in the closed rectangle $\overline{\Pi_{+}^{0}\left(
R\right)  }$ of any fixed length $R$. We emphasize that the solution
$Z_{1}^{\prime}$ has no singularities at the corner points $\left(
0,0\right)  $ and $\left(  0,1\right)  ,$ cf. \cite[Ch. 2]{NaPl}, but third
derivatives of $\breve{Z}_{1}^{\prime}$ are not bounded when $r\rightarrow+0$.
The extension $\widehat{Z}_{1}^{\prime}$ in (\ref{J7}) where $\breve{Z}%
_{1}^{\prime}$ is smoothly continued through the segment $\left\{  x:x_{1}%
\in\left[  0,l\right]  ,\text{\ }x_{2}=0\right\}  $.

Finally, we mention that the correction term $Z_{1}^{\prime}$ in (\ref{47})
was determined in Section \ref{sect3.2} up to the addendum $C_{1}^{0}%
\cos\left(  \pi x_{1}\right)  $ but putting $C_{1}^{0}=0$ in the expansion
(\ref{5222}) defines uniquely the function $Z_{1}^{\prime}$ as well as its
value $Z_{1}^{\prime}\left(  l,0\right)  $ according to (\ref{J7}). Notice
that we also must take $S_{11}^{\prime}=0$ by virtue of (\ref{56}). The
extension $\widehat{Z}_{1}^{\prime}$ of $Z_{1}^{\prime}$ is smooth everywhere
in a neighborhood of $\overline{\Pi_{+}^{0}}$ except at the point $\left(
l,0\right)  $ where it inherits a singularity from (\ref{J7}). Using the
partition of unity $\left\{  1-\chi_{\varepsilon},\chi_{\varepsilon}\right\}
$ makes the last term in (\ref{J4}) smooth in $\overline{\Pi_{+}^{\varepsilon
}}$ but produces additional discrepancy in the Helmholtz equation (\ref{1}).

\subsection{Estimating discrepancies\label{sect6.2}}

First of all, we observe that $\mathcal{F}^{out}=0$ in $\Pi_{+}^{0}$ according
to definition of waves in (\ref{13}) and (\ref{30}). In view of the factor
$\mathfrak{\chi}_{\infty}$ from (\ref{J1}) the first term on the right of
(\ref{J6}) vanishes. Moreover, the Taylor formulas (\ref{43}) and (\ref{46})
assure that%
\[
\left\vert \mathfrak{Z}^{out}\left(  x\right)  -\mathfrak{Z}^{mat}\left(
x\right)  \right\vert +\left\vert \nabla\mathfrak{Z}^{out}\left(  x\right)
-\nabla\mathfrak{Z}^{mat}\left(  x\right)  \right\vert \leq c\varepsilon^{3/2}%
\]
on the rectangle $\left[  3l/2,2l\right]  \times\left[  0,1\right]  $ where
supports of coefficients in the commutator $\left[  \Delta,\chi_{\infty
}\right]  $ are located in. Hence%
\begin{equation}
\left\vert \mathcal{F}^{oma}\left(  x\right)  \right\vert \leq c\varepsilon
^{3/2},\ \ \ \ \mathcal{F}^{oma}\left(  x\right)  =0\text{ \ \ for }x_{1}%
\geq2l. \label{J9}%
\end{equation}
Let us consider the sum
\begin{equation}
\mathcal{F}^{inm}\left(  x\right)  =X_{\varepsilon}\mathcal{F}^{in}%
-X_{\varepsilon}\mathfrak{\chi}_{\infty}\mathcal{F}^{mat}. \label{J10}%
\end{equation}
Outside the finite domain $\Pi_{+}^{\varepsilon}\left(  3l/2\right)  $ it is
equal to%
\begin{gather}
\varepsilon^{-1/2}X_{\varepsilon}\left(  x\right)  \left(  \left(  \Delta
+\pi^{2}\right)  Z_{1}^{0}\left(  x\right)  +\varepsilon\left(  \Delta+\pi
^{2}\right)  Z_{1}^{\prime}\left(  x\right)  -\mathfrak{\chi}_{\infty}\left(
\Delta+\pi^{2}\right)  \mathfrak{Z}^{mat}\left(  x\right)  \right)
+\label{J11}\\
+\varepsilon^{-1/2}X_{\varepsilon}\left(  x\right)  \left(  \lambda
^{\varepsilon}-\pi^{2}\right)  \left(  Z_{1}^{0}\left(  x\right)  +\varepsilon
Z_{1}^{\prime}\left(  x\right)  -\mathfrak{\chi}_{\infty}\left(  x\right)
\mathfrak{Z}^{mat}\left(  x\right)  \right)  =\nonumber\\
=0+\varepsilon^{3/2}\mu(Z_{1}^{0}\left(  x\right)  -\mathfrak{\chi}_{\infty
}\left(  x\right)  \left(  4\mu\right)  ^{-1/4}\cos\left(  \pi x_{2}\right)
\left(  1+i+S_{11}^{0}\left(  1-i\right)  \right)  -\nonumber\\
-\varepsilon(Z_{1}^{\prime}\left(  x\right)  -\mathfrak{\chi}_{\infty}\left(
x\right)  (\left(  4\mu\right)  ^{-1/4}\cos\left(  \pi x_{2}\right)
x_{1}\sqrt{\mu}\left(  1-i+S_{11}^{0}\left(  1+i\right)  \right)  +S_{01}%
^{0}\left(  2\pi\right)  ^{-1/2}e^{i\pi x_{1}}).\nonumber
\end{gather}
where formulas (\ref{42}) and (\ref{T1}) are taken into account. We now use
the representations (\ref{50}) and (\ref{5222}) to conclude that%
\begin{equation}
\left\vert \mathcal{F}^{inm}\left(  x\right)  \right\vert \leq c\varepsilon
^{3/2}e^{-x_{1}\sqrt{3}\pi}\text{ \ \ \ for }x_{1}\geq3l/2. \label{J12}%
\end{equation}
Inside $\Pi_{+}^{\varepsilon}\left(  3l/2\right)  $ we have%
\[
\mathcal{F}^{inm}=-\varepsilon^{2}\mu\mathfrak{Z}^{in}+\varepsilon
^{1/2}\left(  1-\chi_{\varepsilon}\right)  \left(  \Delta+\pi^{2}\right)
\widehat{Z}_{1}^{\prime}-\varepsilon^{1/2}\left[  \mathcal{\vartriangle}%
,\chi_{\varepsilon}\right]  (\widehat{Z}_{1}^{\prime}-Z_{1}^{\prime}\left(
l,0\right)  ).
\]
The inequality%
\[
\varepsilon^{2}\mu\left\vert \mathfrak{Z}^{in}\left(  x\right)  \right\vert
\leq c\varepsilon^{3/2}\text{ \ \ \ in }\overline{\Pi_{+}^{\varepsilon}\left(
3l/2\right)  }%
\]
is evident. Because of the singularity $O\left(  r\left\vert \ln r\right\vert
\right)  $ in (\ref{J7}), estimates of other two terms in (\ref{J11}) involve
the weight function%
\begin{equation}
\rho\left(  x\right)  =r+\left(  1+\left\vert \ln r\right\vert \right)
\label{J14}%
\end{equation}
cf. the Hardy inequality (\ref{hardy}). Since $\widehat{Z}_{1}^{\prime}%
=Z_{1}^{\prime}$ in $\Pi_{+}^{0}$ satisfies the Helmholtz equation from
(\ref{51}), we have%
\begin{align*}
\varepsilon^{1/2}\left(  1-\chi_{\varepsilon}\left(  r\right)  \right)
\left(  \mathcal{\Delta}+\pi^{2}\right)  \widehat{Z}_{1}^{\prime}\left(
x\right)   &  =0,\text{ \ \ }x\in\Pi_{+}^{0}\left(  3l/2\right) \\
\varepsilon^{1/2}\left\vert \left(  1-\chi_{\varepsilon}\left(  r\right)
\right)  \left(  \mathcal{\Delta}+\pi^{2}\right)  \widehat{Z}_{1}^{\prime
}\left(  x\right)  \right\vert  &  \leq c\varepsilon^{1/3}\left\vert
x_{1}\right\vert \left(  \varepsilon+r\right)  ^{-2}\rho\left(  x\right)
,\text{ \ \ }x\in\varpi_{+}^{\varepsilon}=\Pi_{+}^{\varepsilon}\setminus
\Pi_{+}^{0}%
\end{align*}
Observing that, according to the third line in (\ref{J1}), coefficients in the
commutator $\left[  \Delta,\chi_{\varepsilon}\right]  =2\nabla\chi
_{\varepsilon}\cdot\nabla+\Delta\chi_{\varepsilon}$ gets orders $\varepsilon
^{-1}$ and $\varepsilon^{-2},$ respectively but vanish outside the set
$\Cap^{\varepsilon}=\left\{  x\in\Pi_{+}^{\varepsilon}:2\varepsilon
<r<3\varepsilon\right\}  $, we conclude that
\begin{align}
\varepsilon^{1/2}\left[  \Delta,\chi_{\varepsilon}\right]  \left(  \widehat
{Z}_{1}^{\prime}\left(  x\right)  -Z_{1}^{\prime}\left(  l,0\right)  \right)
&  =0,\text{ \ \ }x\in\Pi_{+}^{\varepsilon}\left(  3l/2\right)  \setminus
\Cap^{\varepsilon},\label{J16}\\
\varepsilon^{1/2}\left\vert \left[  \Delta,\chi_{\varepsilon}\left(
\tau\right)  \right]  \left(  \widehat{Z}_{1}^{\prime}\left(  x\right)
-Z_{1}^{\prime}\left(  l,0\right)  \right)  \right\vert  &  \leq
c\varepsilon^{1/2}\left(  \varepsilon^{-1}\left\vert \ln r\right\vert
+\varepsilon^{-2}r\left\vert \ln r\right\vert \right)  \leq\nonumber\\
&  \leq c\varepsilon^{1/2}\left\vert x_{1}\right\vert \left(  \varepsilon
+r\right)  ^{-2}\rho\left(  x\right)  ,\text{ }x\in\Cap^{\varepsilon
}.\nonumber
\end{align}

Finally, we mention that the support of the term $\mathcal{F}^{ima}$ in
(\ref{J6}) belongs to the rectangle $\left[  l+1/\varepsilon,2l+1/\varepsilon
\right]  \times\left[  0,1\right]  ,$ see the first line of (\ref{J1}), where
the remainder $\widetilde{Z}_{1}^{\prime}\left(  x\right)  $ in (\ref{5222})
gets the exponential small order $O\left(  e^{-\sqrt{3}\pi/\varepsilon
}\right)  $, and hence%
\begin{equation}
\left\vert \left[  \Delta,X_{\varepsilon}\right]  \left(  \mathfrak{Z}%
^{in}\left(  x\right)  -\mathfrak{Z}^{mat}\left(  x\right)  \right)
\right\vert =\varepsilon^{1/2}\left\vert \left[  \Delta,X_{\varepsilon
}\right]  \widetilde{Z}_{1}^{\prime}\left(  x\right)  \right\vert \leq
c\varepsilon^{1/2}e^{-3\sqrt{\pi}/\varepsilon}. \label{J17}%
\end{equation}

It remains to consider discrepancies in the Neumann condition (\ref{2}). Since
the cut-off functions $X_{\varepsilon}$ and $\chi_{\infty}$ can be taken
dependent on the longitudinal coordinate $x_{1}$ only, the asymptotic solution
(\ref{J2}) satisfies the homogeneous Neumann condition everywhere on
$\partial\Pi_{+}^{0}$, except on the sides $\Upsilon^{\varepsilon}$ and
$\upsilon^{\varepsilon}$ of the rectangle (\ref{0}), cf. Sections
\ref{sect3.2} and \ref{sect3.3}. Furthermore, $Z_{1}^{0}$ does not depend on
$x_{1}$ and $Z_{1}^{\prime}$ is multiplied in (\ref{J4}) with the cut-off
function $\mathfrak{\chi}_{\varepsilon}$ in the radial variable $r.$ Thus,
$\partial_{1}\mathfrak{Z}^{\varepsilon}=0$ on the short side $\upsilon
^{\varepsilon}.$ Regarding the trace $\mathcal{G}^{\varepsilon}$ of
$\partial_{\nu}\mathfrak{Z}^{\varepsilon}=-\partial_{2}\mathfrak{Z}^{as}$ on
the long side $\Upsilon^{\varepsilon}$ we take the formulas (\ref{51}%
)-(\ref{52}) into account and, similarly to (\ref{J14}) and (\ref{J16}),
obtain%
\begin{equation}
\left\vert \mathcal{G}^{\varepsilon}\left(  x_{1},-\varepsilon\right)
\right\vert \leq c\varepsilon^{3/2}\left(  \varepsilon+r\right)  ^{-1}.
\label{J18}%
\end{equation}
We emphasize that differentiation in $x_{2}$ eliminates $\ln r$ in the first
term of (\ref{J7}).

\subsection{Comparing the approximate and true solutions.\label{sect6.3}}

First of all, we observe that $\mathfrak{Z}^{as}\left(  x\right)
=\mathfrak{Z}^{out}\left(  x\right)  $ as $x_{1}>2l$ by virtue of the
definition (\ref{J1}) of $X_{\varepsilon}$ and $\mathfrak{\chi}_{\infty}$.
Thus, in view of (\ref{31}) and (\ref{J2}), (\ref{J3}) the difference
$\mathcal{R}^{\varepsilon}=Z^{\varepsilon}-\mathfrak{Z}^{\varepsilon}$ loses
the incoming waves $w_{p}^{\varepsilon-}$ and falls into the space $W_{\gamma
}^{1}\left(  \Pi_{+}^{\varepsilon}\right)  _{out}$. Moreover, the
decomposition (\ref{KK1}) of $\mathcal{R}^{\varepsilon}\left(  x\right)  $
contains the coefficients $a_{+}^{\varepsilon}=\widehat{S}_{10}^{\varepsilon}$
and $b_{+}^{\varepsilon}=\widehat{S}_{11}^{\varepsilon}$ defined in
(\ref{S11}) and (\ref{S01}). The integral identity (\ref{K2}) with
$\sigma=-\gamma$ serving for $\mathcal{R}^{\varepsilon}$, involves the
functional%
\begin{equation}
F^{\varepsilon}\left(  v^{\varepsilon}\right)  =\left(  \mathcal{F}%
^{\varepsilon},v^{\varepsilon}\right)  _{\Pi_{+}^{\varepsilon}}-\left(
\mathcal{G}^{\varepsilon},v^{\varepsilon}\right)  _{\Upsilon_{\varepsilon}}
\label{J20}%
\end{equation}
where $\mathcal{F}^{\varepsilon}$ is given in (\ref{J6}) and $\mathcal{G}%
^{\varepsilon}=-\partial_{\nu}\mathfrak{Z}^{\varepsilon}$. If we prove the
inclusion $F^{\varepsilon}\in W_{-\gamma}^{1}\left(  \Pi_{+}^{\varepsilon
}\right)  ^{\ast}$, then the estimate (\ref{U18}) adjusted by the weighting
factor $\left(  \pi^{2}-\lambda^{\varepsilon}\right)  ^{1/4}=\varepsilon
^{1/2}\mu^{1/4}$ demonstrates that
\[
\left\vert \widehat{S}_{10}^{\varepsilon}\right\vert +\varepsilon
^{1/2}\left\vert \widehat{S}_{11}^{\varepsilon}\right\vert \leq c\left\Vert
F^{\varepsilon},W_{-\gamma}^{1}\left(  \Pi_{+}^{\varepsilon}\right)
\right\Vert .
\]
We fix some test function $v^{\varepsilon}\in W_{\gamma}^{1}\left(  \Pi
_{+}^{\varepsilon}\right)  $. The classical one-dimensional Hardy inequality%
\[
\int_{0}^{l}r^{-2}\left\vert \ln\frac{r}{l}\right\vert ^{2}\left\vert V\left(
r\right)  \right\vert ^{2}rdr\leq4\int_{0}^{l}\left\vert \frac{dV}{dr}\left(
r\right)  \right\vert ^{2}rdr,\text{ \ \ }V\in C_{0}^{\infty}\left[
0,l\right)  ,
\]
in a standard way, cf. \cite[Ch.1, \S 4]{MaNaPl}, leads to the relation%
\begin{equation}
\left\Vert \rho^{-1}v^{\varepsilon};L^{2}\left(  \Pi_{+}^{\varepsilon}\left(
2l\right)  \right)  \right\Vert ^{2}\leq c\left\Vert v^{\varepsilon}%
;H^{1}\left(  \Pi_{+}^{\varepsilon}\left(  2l\right)  \right)  \right\Vert
^{2}\leq c_{\gamma}\left\Vert v^{\varepsilon};W_{\gamma}^{1}\left(  \Pi
_{+}^{\varepsilon}\right)  \right\Vert ^{2} \label{hardy}%
\end{equation}
where $\rho$ is the weight factor (\ref{J14}). Moreover, introducing the new
weight factor $\rho_{1}\left(  x\right)  =r\left(  1+\left\vert \ln
r\right\vert \right)  ^{2}$, we derive the weighted trace inequality%
\begin{align}
\int_{\Upsilon^{\varepsilon}}\rho_{1}^{-1}\left\vert v^{\varepsilon
}\right\vert ^{2}dx_{1}  &  =\int_{\Pi_{+}^{\varepsilon}\left(  l\right)
}\dfrac{\partial}{\partial x_{2}}\left(  \mathfrak{\chi}_{0}\rho_{1}%
^{-1}\left\vert v^{\varepsilon}\right\vert ^{2}\right)  dx\leq\nonumber\\
&  \leq c\int_{\Pi_{+}^{\varepsilon}\left(  l\right)  }\left(  \left\vert
\dfrac{\partial v^{\varepsilon}}{\partial x_{2}}\right\vert \rho_{1}%
^{-1}\left\vert v^{\varepsilon}\right\vert +\left(  1+\dfrac{\partial
}{\partial x_{2}}\rho_{1}^{-1}\right)  \left\vert v^{\varepsilon}\right\vert
^{2}\right)  dx\leq\nonumber\\
&  \leq c\int_{\Pi_{+}^{\varepsilon}\left(  l\right)  }\left(  \left\vert
\nabla v^{\varepsilon}\right\vert ^{2}\rho^{-2}\left\vert v^{\varepsilon
}\right\vert ^{2}\right)  dx\leq c_{\gamma}\left\Vert v^{\varepsilon
};W_{\gamma}^{1}\left(  \Pi_{+}^{\varepsilon}\right)  \right\Vert ^{2}.
\label{J22}%
\end{align}
Here, we took into account that $\left\vert \nabla\rho_{1}\left(  x\right)
^{-1}\right\vert \leq c\rho\left(  x\right)  ^{-2}$ and used a cut-off
function $\chi_{0}\in C^{1}\left(  \mathbb{R}\right)  $, $\chi_{0}\left(
x_{2}\right)  =1$ for $x_{2}\leq1/3$ and $\chi_{0}\left(  x_{2}\right)  =1$
for $x_{2}\geq2/3$. The inclusion $F^{\varepsilon}\in W_{-\gamma}^{1}\left(
\Pi_{+}^{\varepsilon}\right)  ^{\ast}$ is obvious because $\mathcal{F}%
^{\varepsilon}$ has a compact support. To estimate the norm $\left\Vert
\mathcal{F}^{\varepsilon};W_{-\gamma}^{1}\left(  \Pi_{+}^{\varepsilon}\right)
^{\ast}\right\Vert $ we apply inequalities obtained in the previous section.
Since $\gamma\in(0,\sqrt{3}\pi)$, the estimates (\ref{J17}) gives%
\[
\left\vert \left(  \mathcal{F}^{ima},v^{\varepsilon}\right)  _{\Pi
_{+}^{\varepsilon}}\right\vert \leq c\varepsilon^{1/2}e^{\left(  \gamma
-\sqrt{3}\pi\right)  /\varepsilon}\int_{l+1/\varepsilon}^{2l+1/\varepsilon
}\int_{0}^{1}e^{-\gamma x_{1}}\left\vert v^{\varepsilon}\left(  x\right)
\right\vert dx\leq c\varepsilon^{3/2}\left\Vert v^{\varepsilon};W_{-\gamma
}^{1}\left(  \Pi_{+}^{\varepsilon}\right)  \right\Vert .
\]
By the formula (\ref{J9}), we have%
\[
\left\vert \left(  \mathcal{F}^{oma},v^{\varepsilon}\right)  _{\Pi
_{+}^{\varepsilon}}\right\vert \leq c\varepsilon^{3/2}\left\Vert
v^{\varepsilon};L^{1}\left(  \Pi_{+}^{\varepsilon}\left(  3l/2\right)
\right)  \right\Vert \leq c\varepsilon^{3/2}\left\Vert v^{\varepsilon
};W_{-\gamma}^{1}\left(  \Pi_{+}^{\varepsilon}\right)  \right\Vert .
\]
Recalling (\ref{J10}) and (\ref{J12}), (\ref{J14}), (\ref{J16}) yields%
\begin{gather*}
\left\vert \left(  \mathcal{F}^{inm},v^{\varepsilon}\right)  _{\Pi
_{+}^{\varepsilon}}\right\vert \leq\left(  c\varepsilon^{3/2}\int_{\Pi
_{\infty}\left(  3l/2\right)  }e^{-x_{1}\sqrt{3}\pi}\left\vert v^{\varepsilon
}\left(  x\right)  \right\vert dx+\varepsilon^{1/2}\int_{\varpi_{+}%
^{\varepsilon}}\dfrac{\left\vert x_{1}\right\vert \rho\left(  x\right)
}{\left(  \varepsilon+r\right)  ^{2}}\left\vert v^{\varepsilon}\left(
x\right)  \right\vert dx+\right. \\
+\left.  \varepsilon^{1/2}\int_{\Cap^{\varepsilon}}\dfrac{\rho\left(
x\right)  }{\left(  \varepsilon+r\right)  ^{2}}\left\vert v^{\varepsilon
}\left(  x\right)  \right\vert dx\right)  \leq\\
\leq c(\varepsilon^{3/2}\left(  \int_{3l/2}^{+\infty}e^{2\left(  \gamma
-\sqrt{3}\pi\right)  x_{1}}dx_{1}\right)  ^{1/2}\left(  \int_{\Pi
_{+}^{\varepsilon}\left(  3l/2\right)  }e^{-2\gamma x_{1}}\left\vert
v^{\varepsilon}\left(  x\right)  \right\vert ^{2}dx\right)  ^{1/2}+\\
+\varepsilon^{1/2}\left(  \int_{\varpi_{+}^{\varepsilon}}\left\vert
x_{1}\right\vert ^{2}\dfrac{\left(  \rho\left(  x\right)  \right)  ^{4}%
}{\left(  \varepsilon+r\right)  ^{4}}dx+\int_{\Cap^{\varepsilon}}%
\dfrac{\left(  \rho\left(  x\right)  \right)  ^{4}}{\left(  \varepsilon
+r\right)  ^{4}}dx\right)  ^{1/2}\left\Vert \rho^{-1}v^{\varepsilon}%
;L^{1}\left(  \Pi_{+}^{\varepsilon}\left(  l\right)  \right)  \right\Vert
\leq\\
\leq\varepsilon^{3/2}\left(  1+\left\vert \ln\varepsilon\right\vert \right)
^{2}\left\Vert v^{\varepsilon};W_{-\gamma}^{1}\left(  \Pi_{+}^{\varepsilon
}\right)  \right\Vert .
\end{gather*}

Finally, we derive from (\ref{J18}) and (\ref{J22}) the following estimate of
the last scalar product in (\ref{J20}):%
\[
\left\vert \left(  \mathcal{G},v^{\varepsilon}\right)  _{\Upsilon
_{\varepsilon}}\right\vert \leq c\varepsilon^{3/2}\left(  \int_{\Upsilon
_{\varepsilon}}\left(  \varepsilon+r\right)  ^{-2}\rho_{1}dx_{1}\right)
^{1/2}\left\Vert \rho_{1}^{-1/2}v^{\varepsilon};L^{2}\left(  \Upsilon
_{\varepsilon}\right)  \right\Vert \leq c\varepsilon^{3/2}\left(  1+\left\vert
\ln\varepsilon\right\vert \right)  ^{3/2}\left\Vert v^{\varepsilon}%
;W_{-\gamma}^{1}\left(  \Pi^{\varepsilon}\right)  \right\Vert .
\]
Collecting the obtained inequalities we conclude that the functional in
(\ref{J20}) meets the estimate%
\begin{equation}
\left\Vert F^{\varepsilon};W_{-\gamma}^{1}\left(  \Pi^{\varepsilon}\right)
^{\ast}\right\Vert \leq c\varepsilon^{3/2}\left(  1+\left\vert \ln
\varepsilon\right\vert \right)  ^{2}. \label{J23}%
\end{equation}

\subsection{Asymptotics of the augmented scattering matrix\label{sect6.4}}

We are in position to formulate the main technical result in the paper. Since
the norm $||\mathcal{R}^{\varepsilon};W_{\gamma}^{1}\left(  \Pi_{+}%
^{\varepsilon}\right)  _{out}||$ in the space with detached asymptotics
contains the coefficients $\widehat{S}_{10}^{\varepsilon}$ and $\widehat
{S}_{11}^{\varepsilon}$ in the representation (\ref{KK1}) of $\mathcal{R}%
^{\varepsilon}$ estimates of the asymptotic remainders in (\ref{S01}) and
(\ref{S11}) follow directly from (\ref{J23}) and (\ref{U18}), (\ref{42}). A
similar outcome for $S_{00}^{\varepsilon}$ can be obtained by repeating word
by word calculations in the previous section based on the formal asymptotic
(\ref{60}), (\ref{64}). We only mention that the discrepancy (\ref{Z02}) on
the small side $\upsilon^{\varepsilon}$ of the box $\varpi_{+}^{\varepsilon}$
can be considered as follows:%
\begin{align*}
\left(  2\pi\right)  ^{1/2}\left\vert \sin\left(  \pi l\right)  \int
_{-\varepsilon}^{0}v^{\varepsilon}\left(  l,x_{2}\right)  dx_{2}\right\vert
&  \leq c\varepsilon^{1/2}\left(  1+\left\vert \ln\varepsilon\right\vert
\right)  \left\Vert r^{-1/2}\left(  1+\left\vert \ln r\right\vert \right)
^{-1}v^{\varepsilon};L^{2}\left(  \upsilon^{\varepsilon}\right)  \right\Vert
\leq\\
&  \leq c\varepsilon\left(  1+\left\vert \ln\varepsilon\right\vert \right)
\left\Vert v^{\varepsilon};W_{-\gamma}^{1}\left(  \Pi_{+}^{\varepsilon
}\right)  \right\Vert
\end{align*}
where a weighted trace inequality of type (\ref{J22}) is applied.

\begin{theorem}
\label{TheoremASY}Remainders in the asymptotic forms (\ref{S11})-(\ref{S01})
enjoy the estimate
\begin{equation}
\left\vert \widehat{S}_{11}^{\varepsilon}\right\vert +\varepsilon
^{-1/2}\left\vert \widehat{S}_{01}^{\varepsilon}\right\vert +\left\vert
S_{00}^{\varepsilon}-1\right\vert \leq c\varepsilon\left(  1-\left\vert
\ln\varepsilon\right\vert \right)  ^{2}. \label{est}%
\end{equation}

\end{theorem}

We finally note that formulas (\ref{T3}), (\ref{T999}) is the only issue in
the paper which deals with $S_{00}^{\varepsilon}$ but requires much less
accurate information.

\subsection{Dependence on $\bigtriangleup l$ and $\bigtriangleup\mu
$\label{sect6.5}}

Let $\varepsilon$ be fixed small and positive. We take $l=\pi k+\bigtriangleup
l$ and make the change of coordinates%
\begin{equation}
x\rightarrow\mathbf{x=}\left(  \mathbf{x}_{1},\mathbf{x}_{2}\right)  =\left(
x_{1},x_{2}\right)  ,\ \ \ \ \ \ \mathbf{x}_{1}=\left(  1-\chi_{k}\left(
x_{1}\right)  \right)  x_{1}+\chi_{k}\left(  x_{1}\right)  \left(
x_{1}-\bigtriangleup l\right)  \label{J51}%
\end{equation}
where $\chi_{k}$ is a smooth cut-off function, $\chi_{k}\left(  x_{1}\right)
=1$ for $\left\vert x_{1}-\pi k\right\vert <\pi/3$ and $\chi_{k}\left(
x_{1}\right)  =0$ for $\left\vert x_{1}-\pi k\right\vert >2\pi/3$. If
$\bigtriangleup l$ is small, this change is nonsingular. Moreover, it
transforms $\Pi_{l}^{\varepsilon}$ into $\Pi_{\pi k}^{\varepsilon}$ and turns
the Helmholtz operator $\bigtriangleup+\pi^{2}-\varepsilon^{2}\left(
\mu+\bigtriangleup\mu\right)  $ into the second-order differential operator
$L^{\varepsilon}\left(  \bigtriangleup l,\bigtriangleup\mu;\mathbf{x}%
,\nabla_{\mathbf{x}}\right)  $ whose coefficients depend smoothly on
$\bigtriangleup\mu$ and $\bigtriangleup l.$ Clearly, $L^{\varepsilon}\left(
0,0;\mathbf{x},\nabla_{\mathbf{x}}\right)  =\bigtriangleup_{\mathbf{x}}%
+\pi^{2}-\varepsilon^{2}\mu$. Owing the Fourier method, we can rewrite the
element $S_{11}^{\varepsilon}=S_{11}^{\varepsilon}\left(  \bigtriangleup
\mu,\bigtriangleup l\right)  $ of the augmented scattering matrix as the
integral%
\begin{equation}
S_{11}^{\varepsilon}\left(  \bigtriangleup\mu,\bigtriangleup l\right)
=\alpha_{11}^{\varepsilon}\left(  \bigtriangleup\mu\right)  \int_{Q_{k}}%
Z_{11}^{\varepsilon}\left(  \bigtriangleup\mu,\bigtriangleup l;x\right)  dx
\label{J53}%
\end{equation}
over the rectangle $Q_{k}=\left(  \pi\left(  k+1\right)  ,\pi\left(
k+2\right)  \right)  \times\left(  0,1\right)  $ where $\mathbf{x}=x$
according to (\ref{J51}). Due to the general result in the perturbation theory
of linear operators, see, e.g., \cite{HilPhi, Kato}, the special solution
$Z_{11}^{\varepsilon}\left(  x\right)  =$ $Z_{11}^{\varepsilon}\left(
\bigtriangleup\mu,\bigtriangleup l;x\right)  $ rewritten in the coordinates
$\mathbf{x}$ depend smoothly on $\left(  \bigtriangleup\mu,\bigtriangleup
l\right)  \in\overline{\mathbb{B}_{\rho}}.$The coefficient $\alpha
_{11}^{\varepsilon}\left(  \bigtriangleup\mu\right)  $ in (\ref{J53}) is also
a smooth function whose exact form is of no need. Thus, the element
(\ref{J53}) inherits this smooth dependence while the remainder $\widetilde
{S}_{11}^{\varepsilon}\left(  \bigtriangleup\mu,\bigtriangleup l\right)  $ in
the representation (\ref{S11}) gets the same property according to the formula
(\ref{58}) for $S_{11}^{0}\left(  \bigtriangleup\mu,\bigtriangleup l\right)  $
written in the variables (\ref{T1}).

Similar operations apply to $S_{01}^{\varepsilon}\left(  \bigtriangleup
\mu,\bigtriangleup l\right)  $ and $\widetilde{S}_{01}^{\varepsilon}\left(
\bigtriangleup\mu,\bigtriangleup l\right)  .$

Finally, recalling our examination in Section \ref{sect5.5} and Theorem
\ref{TheoremASY}, we formulate the result.

\begin{proposition}
\label{PropositionANAL} The remainders in the asymptotic formulas (\ref{S11}),
(\ref{58}) and (\ref{S01}), (\ref{anti}) satisfy the inequality%
\[
\left\vert \nabla_{\left(  \bigtriangleup\mu,\bigtriangleup l\right)
}\widetilde{S}_{11}^{\varepsilon}\left(  \bigtriangleup\mu,\bigtriangleup
l\right)  \right\vert +\varepsilon^{-1/2}\left\vert \nabla_{\left(
\bigtriangleup\mu,\bigtriangleup l\right)  }\widetilde{S}_{01}^{\varepsilon
}\left(  \bigtriangleup\mu,\bigtriangleup l\right)  \right\vert \leq
c\varepsilon\left(  1+\left\vert \ln\varepsilon\right\vert ^{2}\right)
,\text{ \ \ }\left(  \bigtriangleup\mu,\bigtriangleup l\right)  \in
\overline{\mathbb{B}_{\rho}}.
\]

\end{proposition}

\section{The uniqueness assertions\label{sect7}}

\subsection{Eigenvalues in the vicinity of the threshold $\pi^{2}%
.$\label{sect7.1}}

Let us adapt a trick from \cite[\S 7]{na489} for the box-shaped perturbation
(\ref{0}) and conclude with the uniqueness mentioned in Theorem
\ref{TheoremEX}.

Assume that there exists an infinitesimal sequence $\left\{  \varepsilon
_{k}\right\}  _{k\in\mathbb{N}}$ such that the problem (\ref{9}), (\ref{10})
in the semi-infinite waveguide $\Pi_{l_{k}}^{\varepsilon}$ has two eigenvalues
$\lambda_{1}^{\varepsilon_{k}}$ and $\lambda_{2}^{\varepsilon_{k}}$ while
\begin{equation}
\varepsilon_{k}\rightarrow+0,\text{ \ \ }l_{k}\rightarrow l_{0}%
>0,\ \ \ \ \ \ \ \ \lambda_{j}^{\varepsilon_{k}}=\pi^{2}+\widehat{\lambda}%
_{j}^{\varepsilon_{k}}\in(0,\pi^{2}],\text{ \ \ }\widehat{\lambda}%
_{j}^{\varepsilon_{k}}\rightarrow0,\text{ \ \ }j=1,2. \label{U0}%
\end{equation}
In what follows we write $\varepsilon$ instead of $\varepsilon_{k}.$ The
corresponding eigenfunctions $u_{1}^{\varepsilon}$ and $u_{2}^{\varepsilon}$
are subject to the normalization and orthogonality conditions%
\begin{equation}
\left\Vert u_{j}^{\varepsilon};L^{1}\left(  \Pi_{+}^{\varepsilon}\left(
2l\right)  \right)  \right\Vert =1,\text{ \ \ }\left(  u_{1}^{\varepsilon
},u_{2}^{\varepsilon}\right)  _{\Pi_{+}^{\varepsilon}}=0, \label{U00}%
\end{equation}
cf. (\ref{M0}). Repeating with evident changes our arguments in Section
\ref{sect5.3} we observe that the restrictions $u_{j0}^{\varepsilon}$ of
$u_{j}^{\varepsilon}$ onto $\Pi_{+}^{0}$ converge to $u_{j0}^{0}$ weakly in
$W_{-\gamma}^{1}\left(  \Pi_{+}^{0}\right)  $ and strongly in $L^{2}\left(
\Pi_{+}^{0}\left(  2l\right)  \right)  $. Furthermore, the limits satisfy the
formula (\ref{M4}) and the following integral identity, see (\ref{M5}):%
\begin{equation}
\left(  \nabla u_{j0}^{0},\nabla v\right)  _{\Pi_{+}^{0}}=\pi^{2}\left(
u_{j0}^{0},v\right)  _{\Pi_{+}^{0}},\text{ \ \ }v\in C_{c}^{\infty}%
(\overline{\Pi_{+}^{0}}). \label{U1}%
\end{equation}
Any solution in $W_{-\gamma}^{1}\left(  \Pi_{+}^{0}\right)  $ with $\gamma
\in\left(  \beta_{1},\beta_{2}\right)  $ of the homogeneous Neumann problem
(\ref{U1}) in the semi-strip $\Pi_{+}^{0}=\left(  0,+\infty\right)
\times\left(  0,1\right)  $ is a linear combination of two bounded solutions
(\ref{491}) and (\ref{492}):%
\begin{equation}
u_{j}^{0}\left(  x\right)  =c_{j1}\cos\left(  \pi x_{1}\right)  +c_{j2}%
\cos\left(  \pi x_{2}\right)  . \label{U2}%
\end{equation}
Let us prove that $c_{11}=c_{21}=0$ in (\ref{U2}). Since the trapped mode
$u_{j}^{\varepsilon}$ has an exponential decay at infinity, the Green formula
in $\Pi_{\infty}\left(  3l/2\right)  $ with it and the bounded function
$e^{\pm ix_{1}\sqrt{\lambda_{j}^{\varepsilon}}}$ assures that
\begin{equation}
\int_{0}^{1}e^{\pm ix_{1}\sqrt{\lambda_{j}^{\varepsilon}}}\left.  \left(
\partial_{1}u_{j}^{\varepsilon}\left(  x\right)  \mp i\sqrt{\lambda
_{j}^{\varepsilon}}u_{j}^{\varepsilon}\left(  x\right)  \right)  \right\vert
_{x_{1}=3l/2}dx_{2}=0. \label{U3}%
\end{equation}
The local estimate (\ref{M5}) in $\varpi^{\prime}\ni\left(  3l/2,x_{2}\right)
,$ $x_{2}\in\left(  0,1\right)  ,$ and formulas in (\ref{U0}), (\ref{U00})
allow us to compute the limit of the left-hand side of (\ref{U3}) and obtain
that%
\begin{equation}
e^{\pm i3l\pi/2}\int_{0}^{1}\left(  \dfrac{\partial u_{j0}^{0}}{\partial
x_{1}}\left(  \frac{3}{2}l,x_{2}\right)  \pm i\pi u_{j0}^{0}\left(  \frac
{3}{2}l,x_{2}\right)  \right)  dx_{2}=0. \label{U4}%
\end{equation}
Inserting (\ref{U2}) into (\ref{U4}), we see that $c_{j1}=0,$ indeed.

\begin{remark}
\label{RemarkUU1}. If $\widetilde{\lambda}_{j}^{\varepsilon}>0$ and
$\lambda_{j}^{\varepsilon}>\pi^{2}$ in (\ref{U0}), one may use the Green
formula in $\Pi_{\infty}\left(  3l/2\right)  $ with four bounded functions
$e^{\pm ix_{1}\sqrt{\lambda_{j}^{\varepsilon}}}$ and $e^{\pm ix_{1}%
\sqrt{\lambda_{j}^{\varepsilon}-\pi^{2}}}\cos\left(  \pi x_{2}\right)  $. In
this way one derives the equalities $c_{j1}=c_{j2}=0$ (see \cite[\S 7]{na489}
for details) and concludes that a small neighborhood of the threshold $\pi
^{2}$ can contain only eigenvalues indicated in (\ref{U0}).$\ $The same
reasoning show that the problem (\ref{9}), (\ref{10}) cannot get an eigenvalue
$\lambda^{\varepsilon}\rightarrow+0$ as $\varepsilon\rightarrow+0$%
.\ \ \ $\boxtimes$
\end{remark}

Since $c_{j1}=0,$ the limit normalization (\ref{M4}) shows that $u_{j}%
^{0}\left(  x\right)  =l^{-1/2}\cos\left(  \pi x_{2}\right)  ,$\ $j=1,2.$
Moreover, Theorem \ref{TheoremKA} (2) applied to the trapped mode
$u_{j}^{\varepsilon}\in H^{1}\left(  \Pi_{+}^{\varepsilon}\right)  \subset
W_{-\gamma}^{1}\left(  \Pi_{+}^{\varepsilon}\right)  $ gives the formula%
\begin{equation}
u_{j}^{\varepsilon}\left(  x\right)  =B_{j}^{\varepsilon}e^{-x_{1}\sqrt
{\pi^{2}-\lambda_{j}^{\varepsilon}}}\cos\left(  \pi x_{2}\right)
+\widetilde{u}_{j}^{\varepsilon}\left(  x\right)  ,\ \ \ \ \ \ \left\vert
b_{j}^{\varepsilon}\right\vert +\left\Vert \widetilde{u}_{j}^{\varepsilon
};W_{\gamma}^{1}\left(  \Pi_{+}^{\varepsilon}\right)  \right\Vert \leq
c\left\Vert u_{j}^{\varepsilon};W_{-\gamma}^{1}\left(  \Pi_{+}^{\varepsilon
}\right)  \right\Vert , \label{U6bis}%
\end{equation}
where $\gamma\in\left(  \beta_{1},\beta_{2}\right)  ,$ $c$ is independent of
$\varepsilon$ according to the content of Section \ref{sect5.5} and the waves
$w_{0}^{\varepsilon\pm}$ in (\ref{13}) and $v_{1}^{\varepsilon+}$ in
(\ref{27}) do not appear in the expansion of $u_{j}^{\varepsilon}$ due to the
absense of decay at infinity. Since the right-hand side of (\ref{U6bis}) is
uniformly bounded in $\varepsilon=\varepsilon_{k}$ , $k\in\mathbb{N}$ (see
Section \ref{sect5.5} again), we have%
\[
bB_{j}^{\varepsilon}\rightarrow l^{-1/2},\text{ \ \ \ \ \ }\widetilde{u}%
_{j0}^{\varepsilon}\rightarrow0\text{\ weakly in }W_{\gamma}^{1}\left(
\Pi_{+}^{\varepsilon}\right)  \text{ }%
\]
along a subsequence of $\left\{  \varepsilon_{k}\right\}  _{k\in\mathbb{N}}.$
Moreover, the last equality in (\ref{U00}) turns into%
\[
0=\left(  u_{1}^{\varepsilon},u_{2}^{\varepsilon}\right)  _{\Pi_{+}%
^{\varepsilon}}=\int_{\Pi_{+}^{\varepsilon}}\cos^{2}\left(  \pi x_{2}\right)
e^{-x_{1}\Lambda_{\varepsilon}}dx+\left(  u_{1}^{\varepsilon}-\widetilde
{u}_{1}^{\varepsilon},\widetilde{u}_{2}^{\varepsilon}\right)  _{\Pi
_{+}^{\varepsilon}}+\left(  \widetilde{u}_{1}^{\varepsilon},u_{2}%
^{\varepsilon}\right)  _{\Pi_{+}^{\varepsilon}}=\frac{1}{2}\Lambda
_{\varepsilon}^{-1}b_{1}^{\varepsilon}\overline{b_{2}^{\varepsilon}}+O\left(
1\right)  .
\]
We multiply this relation with $\Lambda_{\varepsilon}=\sqrt{\pi^{2}%
-\lambda_{1}^{\varepsilon}}+\sqrt{\pi^{2}-\lambda_{2}^{\varepsilon}%
}\rightarrow0,$ as $\varepsilon\rightarrow0$ and derive from (\ref{U7}) the
absurd formula $o(1)=b_{1}^{\varepsilon}\overline{b_{2}^{\varepsilon}%
}\rightarrow l^{-2}.$ Thus, there can exist at most one eigenvalue indicated
in (\ref{U0}).

\subsection{Absence of eigenfunctions which are odd in $x_{1}.$\label{sect7.2}%
}

In Section \ref{sect1.3} we have changed the original problem (\ref{1}),
(\ref{2}) in $\Pi^{\varepsilon}$ for the Neumann problem (\ref{9}), (\ref{10})
in the half $\Pi_{+}^{\varepsilon}$ of the waveguide while assuming that an
eigenfunction is even in $x_{1}$. Replacing (\ref{10}) by the mixed boundary
conditions%
\begin{equation}
\partial_{\nu}u^{\varepsilon}\left(  x\right)  =0,\text{ \ }x\in\partial
\Pi_{+}^{\varepsilon},\text{ \ \ }x_{1}>0,\text{ \ \ }u_{+}^{\varepsilon
}\left(  x\right)  =0,\text{ }x\in\partial\Pi_{+}^{\varepsilon},\text{
\ \ }x_{1}=0, \label{U8bis}%
\end{equation}
we deal with the alternative, namely an eigenfunction is odd in $x_{1}$ and
therefore vanishes at $\Gamma^{\varepsilon}=\left\{  x:x_{1}=0,\ x_{2}%
\in\left(  -\varepsilon,1\right)  \right\}  .$ The variational formulation of
the problem (\ref{9}), (\ref{U8bis}),%
\[
\left(  \nabla u_{+}^{\varepsilon},\nabla v^{\varepsilon}\right)  _{\Pi
_{+}^{\varepsilon}}=\lambda_{+}^{\varepsilon}\left(  u_{+}^{\varepsilon
},v^{\varepsilon}\right)  _{\Pi_{+}^{\varepsilon}},\text{ \ }v^{\varepsilon
}\in H_{0}^{1}\left(  \Pi_{+}^{\varepsilon};\Gamma^{\varepsilon}\right)
\text{,}%
\]
involves a subspace of functions in $H_{0}^{1}\left(  \lambda\right)  $ which
are null on $\Gamma^{\varepsilon}.$ Evident modifications of considerations in
Sections \ref{sect5} and \ref{sect6} adapt all our results to the mixed
boundary value problem (\ref{9}), (\ref{U8bis}). The only but important
difference is that the solutions (\ref{491}), (\ref{492}) of the limit Neumann
problem in the half-strip $\Pi_{+}^{0}=\mathbb{R}_{+}\mathbb{\times}\left(
0,1\right)  $ now turn into the following ones:%
\[
u_{0}^{0}\left(  x\right)  =i\sin\left(  \pi x_{1}\right)  ,\text{ \ \ }%
u_{1}^{0}\left(  x\right)  =x_{1}\cos\left(  \pi x_{2}\right)
\]
Thus, supposing that, for an infinitesimal sequence $\left\{  \varepsilon
_{k}\right\}  _{k\in\mathbb{N}},$ the problem (\ref{9}), (\ref{U8bis}) in
$\Pi_{l_{k}+}^{\varepsilon}$ has an eigenvalue $\lambda_{1}^{\varepsilon_{k}}$
with the properties (\ref{U0}) at $j=1,$ we obtain that a non-trivial limit
$u_{01}^{0}$, cf. (\ref{U2}), of the corresponding eigenfunction
$u_{1}^{\varepsilon_{k}}$ becomes
\[
u_{10}^{0}\left(  x\right)  =c_{11}\sin\left(  \pi x_{1}\right)  +c_{12}%
x_{1}\sin\left(  \pi x_{2}\right)  .
\]
Now, in contrast to Section \ref{sect7.1}, we may insert $u_{1}^{\varepsilon
_{k}}$ into the Green formula in $\Pi_{\infty}\left(  3l/2\right)  $ together
with one of three bounded functions $e^{\pm i\sqrt{\lambda_{1}^{\varepsilon
_{k}}}x_{1}}$ and $e^{-\sqrt{\pi^{2}-\lambda_{1}^{\varepsilon_{k}}}}%
\cos\left(  \pi x_{2}\right)  $. Similarly to (\ref{U3}) and (\ref{U4}), these
possibilities allow us to conclude that $c_{11}=0$, $c_{12}=0$ and, hence,
$u_{10}^{0}=0$. The observed contradiction and Remark \ref{RemarkUU1} which
remains true for the problem (\ref{9}), (\ref{U8bis}), confirm the absence of
eigenvalues in a small neighborhood of the threshold $\lambda^{\varepsilon
}=\pi^{2}.$

\subsection{Absence of eigenvalues at a distance from the
threshold.\label{sect7.3}}

For any $\lambda\in\left(  0,\pi^{2}\right)  $, the limit Neumann problem in
$\Pi_{+}^{0}=\mathbb{R}_{+}\mathbb{\times}\left(  0,1\right)  $ has the
solutions%
\begin{align}
Z_{0}^{00}\left(  \lambda,x\right)   &  =\left(  2k\left(  \lambda\right)
\right)  ^{-1/2}\left(  e^{-ik\left(  \lambda\right)  x_{1}}+e^{ik\left(
\lambda\right)  x_{1}}\right)  ,\label{U11bis2}\\
Z_{1}^{00}\left(  \lambda,x\right)   &  =\left(  2k_{1}\left(  \lambda\right)
\right)  ^{-1/2}\left(  (e^{k_{1}\left(  \lambda\right)  x_{1}}+ie^{-k_{1}%
\left(  \lambda\right)  x_{1}}\right)  \cos\left(  \pi x_{2}\right)  +i\left(
e^{k_{1}\left(  \lambda\right)  x_{1}}-ie^{-k_{1}\left(  \lambda\right)
x_{1}}\right)  \cos\left(  \pi x_{2}\right)  =\nonumber\\
&  =\left(  2k_{1}\left(  \lambda\right)  \right)  ^{-1/2}\left(  1+i\right)
(e^{k_{1}\left(  \lambda\right)  x_{1}}+e^{-k_{1}\left(  \lambda\right)
x_{1}}\cos\left(  \pi x_{2}\right) \nonumber
\end{align}
where $k\left(  \lambda\right)  =\sqrt{\lambda}$, $k_{1}\left(  \lambda
\right)  =\sqrt{\pi^{2}-\lambda}$ and, thus, the augmented scattering matrix
takes the form%
\begin{equation}
S^{00}=\left(
\begin{array}
[c]{cc}%
1 & 0\\
0 & i
\end{array}
\right)  . \label{U12bis}%
\end{equation}
To fulfill the criterion (\ref{33}), the perturbation $\varpi_{+}%
^{\varepsilon}$ in the waveguide $\Pi_{+}^{\varepsilon}$ has to turn the
right-hand bottom element of the matrix (\ref{U12bis}) into $-1$ that cannot
be made, e.g., for any $\lambda\in\left[  0,\pi^{2}-c\sqrt{\varepsilon
}\right]  ,$ $c>0$ and a small $\varepsilon$. The latter fact may be easily
verified by either constructing asymptotics in Sections \ref{sect3},
\ref{sect6}, or applying a perturbation argument as in Section \ref{sect5.5}.
Notice that we have succeeded in Sections \ref{sect3.2}, \ref{sect4.2} to
construct asymptotics (\ref{S11}) with the main term $S_{11}^{0}=-1$ because
the spectral parameter (\ref{42}) stays too close to the threshold and the
augmented scattering matrix in $\Pi_{+}^{0}$ is not continuous at $\lambda
=\pi^{2}$, cf. Section \ref{sect5.5}. In the case of the mixed boundary value
problem (\ref{9}), (\ref{U8bis}) evident changes in solutions (\ref{U11bis2})
give the matrix $S^{00}=diag\left\{  -1,-i\right\}  $ instead of
(\ref{U12bis}) but our final conclusion remains the same.

\subsection{Interferences on the uniqueness\label{sect7.4}}

The material of the previous three sections proves the last assertion in
Theorem \ref{TheoremEX}. The interval $\left(  0,\pi^{2}\right)  $ where the
eigenvalue (\ref{T10}) is unique in the waveguide $\Pi_{l}^{\varepsilon}$ with
$l=l_{k}\left(  \varepsilon\right)  $ and a fixed $\varepsilon\in\left(
0,\varepsilon_{k}\right)  $ can be enlarged up to $\left(  0,\pi^{2}%
+c\sqrt{\varepsilon}\right)  ,$ $c>0,$ due to Remark \ref{RemarkUU1}.
Moreover, enhancing our consideration in Section \ref{sect7.3} by dealing with
the exponential waves $e^{\pm x_{1}\sqrt{4\pi^{2}-\lambda}}\cos\left(  2\pi
x_{2}\right)  $ and the augmented scattering matrix of size $3\times3$, cf.
\cite{na489}, confirms that any $\lambda\in\left[  \pi^{2}+c\sqrt{\varepsilon
},4\pi^{2}-c\sqrt{\varepsilon}\right)  $ cannot be an eigenvalue as well. We
will discuss the higher thresholds $\pi^{2}k^{2}$ with $k=2,3,...$ in Section
\ref{sect8.2}.

To confirm Theorem \ref{TheoremUN}, we use a similar reasoning. In this
way,,we recall the asymptotic formulas (\ref{S01}), (\ref{anti}) and
(\ref{est}) and observe that $S_{01}^{\varepsilon}$ cannot vanish for a small
$\varepsilon$ when the length parameter (\ref{T11}) stays outside the segment%
\begin{equation}
\left[  \pi k-c\varepsilon\left(  1+\left\vert \ln\varepsilon\right\vert
\right)  ^{2},\pi k+c\varepsilon\left(  1+\left\vert \ln\varepsilon\right\vert
\right)  ^{2}\right]  \label{U13bis}%
\end{equation}
If $\lambda^{\varepsilon}$ belongs to (\ref{U13bis}), the uniqueness of the
solution $\left(  \bigtriangleup\mu,\bigtriangleup l\right)  $ of the abstract
equation (\ref{T7}) which is equivalent to the criterion in Theorem
\ref{TheoremA} follows from the contraction principle.

\section{Available generalizations\label{sect8}}

\subsection{Eigenvalues in the discrete spectrum\label{sect8.1}}

Let us consider the mixed boundary value problem (\ref{1}), (\ref{7}). As in
Section \ref{sect1.3} we reduce it to the half (\ref{11}) of the perturbed
waveguide $\Pi^{\varepsilon}=\Pi\cup\varpi^{\varepsilon}$, cf. (\ref{9}),
(\ref{10}):%
\begin{align}
-\Delta u_{+}^{\varepsilon}\left(  x\right)   &  =\lambda_{+}^{\varepsilon
}u_{+}^{\varepsilon}\left(  x\right)  ,\text{ \ \ }x\in\Pi_{+}^{\varepsilon
},\text{ }u_{+}^{\varepsilon}\left(  x_{1},1\right)  =0,\text{\ \ }%
x_{1}>0,\label{X1}\\
\partial_{\nu}u_{+}^{\varepsilon}\left(  x\right)   &  =0,\text{ \ }%
x\in\partial\Pi_{+}^{\varepsilon},\text{ \ \ }x_{2}<1.\nonumber
\end{align}
If $\lambda^{\varepsilon}\in\left(  0,\pi^{2}/4\right)  $ stays below the
continuous spectrum $\ \wp_{co}^{M}=\left[  \pi^{2}/4,+\infty\right)  $ of the
problem (\ref{X1}), there is no oscillating wave but deal with the exponential
waves%
\[
v_{1/2}^{\varepsilon\pm}\left(  x\right)  =\left(  k_{1/2}^{\varepsilon
}\right)  ^{-1/2}e^{\pm k_{1/2}^{\varepsilon}x_{1}}\cos\left(  \frac{\pi}%
{2}x_{2}\right)  ,\text{ \ \ }k_{1/2}^{\varepsilon}=\sqrt{\frac{\pi^{2}}%
{4}-\lambda^{\varepsilon}}%
\]
and, similarly to (\ref{27}), (\ref{30}), compose the linear combinations%
\[
w_{1/2}^{\varepsilon\pm}\left(  x\right)  =2^{-1/2}\left(  v_{1/2}%
^{\varepsilon+}\left(  x\right)  \mp v_{1/2}^{\varepsilon-}\left(  x\right)
\right)  .
\]
The conditions (\ref{28}), (\ref{29}) and (\ref{19}) with $p,q=1/2$ are
satisfied and we may determine the augmented scattering matrix $S^{\varepsilon
}$ which now is a scalar. Theorem \ref{TheoremA} remains valid and, therefore,
the equality%
\begin{equation}
S^{\varepsilon}=-1 \label{X2}%
\end{equation}
states a criterion for the existence of a trapped mode. Constructing
asymptotics of $S^{\varepsilon}$ and solving the equation (\ref{X2}) yield the
relation (\ref{8}) for an eigenvalue in the discrete spectrum of the problems
(\ref{X1}) and (\ref{1}), (\ref{7}). Repeating arguments from Sections
\ref{sect5} - \ref{sect7} proves estimates of the asymptotic remainders as
well as the uniqueness of the eigenvalue $\lambda_{+}^{\varepsilon}\in
\wp_{d_{i}}^{M},$ however, for any $l>0.$ The latter conclusion requires to
explain a distinction between analysis of isolated and embedded eigenvalues.

The main difference is caused by the application of the criterion (\ref{33})
which in the case of the scalar $S^{\varepsilon}$ changes just into one
equation%
\begin{equation}
\operatorname{Re}S^{\varepsilon}=-1 \label{X3}%
\end{equation}
which is equivalent to (\ref{X2}) because $\left\vert S^{\varepsilon
}\right\vert =1.$ As a result we may satisfy (\ref{X3}) by choosing
$\bigtriangleup\mu$ and do not need the additional parameter $\bigtriangleup
l$ in (\ref{T1}) which was used in Section \ref{sect4} to solve the system
(\ref{T2}) of two transcendental equations. In other words, the absence of
oscillating waves crucially restricts a possible position of $S_{11}%
^{\varepsilon}=S^{\varepsilon}$ to the unit circle on the complex plane while
the entry $S_{11}^{\varepsilon}$ in the previous unitary matrix
$S^{\varepsilon}$ of size $2\times2,$ see Section \ref{sect2.2}, can step
aside from $\mathbb{S}$ and a fine tuning by means of $\bigtriangleup l$ is
necessary to assure the equality (\ref{33}).

\subsection{Higher thresholds\label{sect8.2}}

\begin{figure}
[ptb]
\begin{center}
\ifcase\msipdfoutput
\includegraphics[
height=1.0179in,
width=4.5057in
]%
{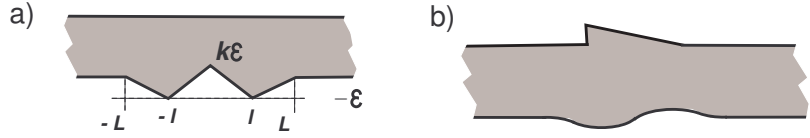}%
\else
\includegraphics[
height=1.0179in,
width=4.5057in
]%
{C:/Users/Desktop/Documents/Google Drive/Articoli/Articoli-in-Preparazione/CardoneDuranteNazarov1/graphics/CaDuNa6__6.pdf}%
\fi
\caption{Other type of piecewise smooth perturbations}%
\label{f6}%
\end{center}
\end{figure}

A straightforward modification of our approach may be used for an attempt to
construct embedded eigenvalues near the thresholds $\pi^{2}k^{2},$ $k=2,3,...$
of the continuous spectrum$\wp_{co}$ of the problem (\ref{9}), (\ref{10}) in
$\Pi_{+}^{\varepsilon}.$ At the same time, the number of oscillating outgoing
waves at the threshold $\pi^{2}k^{2}$ equals $k$ and, therefore, size of the
augmented scattering matrix becomes $\left(  k+1\right)  \times\left(
k+1\right)  .$ In this case the fine tuning needs at least $k$ free
parameters, cf. \cite{na489, na546}, instead of only one $\bigtriangleup l$ in
Section \ref{sect4}. Additional parameters can be easily introduced when the
perturbed wall is a broken line like in fig. \ref{f6}, a, with $l,$ $L$ and
$k$. The amplification of the augmented scattering matrix does not affect the
criterion (\ref{33}) in Theorem \ref{TheoremA}, in Sections \ref{sect3} and
\ref{sect4}.%

In the mirror symmetry with respect to the line $\left\{  x:x_{1}=0\right\}  $
is denied, see fig. \ref{f6}, b, then we have to analyze the problem
(\ref{1}), (\ref{2}) in the intact waveguide $\Pi^{\varepsilon}$ where the
augmented scattering matrix gets rise of size even in the case $\lambda
^{\varepsilon}\leq\pi^{2}.$ In this sense the box-shaped perturbation is
optimal because it demonstrates all technicalities but reduces the
computational details to the necessary minimum. A preliminary assessment
predicts that embedded eigenvalues of the problem (\ref{9}), (\ref{10}) in
$\Pi_{+}^{\varepsilon}$ do not appear near any threshold $\pi^{2}k^{2}$ with
$k>1$ but we are not able to verify this fact rigorously.

\subsection{The Dirichlet boundary condition.\label{sect8.3}}

All procedures described above can be applied to detect eigenvalues of the
Helmholtz equation (\ref{1}) in the quantum waveguide $\Pi^{\varepsilon}$, cf.
\cite{Sim}, with the Dirichlet condition (\ref{5}). However, the asymptotic
structures must be modified a bit due to the following observation. The
correction term $Z^{\prime}$ in the inner asymptotic expansion%
\[
Z^{\varepsilon}\left(  x\right)  =\sin\left(  \pi x_{2}\right)  +\varepsilon
Z^{\prime}\left(  x\right)  +....,
\]
cf. (\ref{62}), must be found out from the mixed boundary value problem in the
semi-strip%
\begin{align*}
-\Delta Z^{\prime}\left(  x\right)   &  =\pi^{2}Z^{\prime}\left(  x\right)
,\text{\ }x\in\Pi_{+}^{0},\text{ \ \ \ }-\partial_{1}Z^{\prime}\left(
0,x_{2}\right)  =0,\text{\ }x_{2}\in\left(  0,1\right)  ,\\
Z^{\prime}\left(  x_{1},1\right)   &  =0,\text{\ }x_{1}>0,\text{
\ \ \ }Z^{\prime}\left(  x_{1},1\right)  =\pi^{2},\text{\ }x_{1}\in\left(
0,l\right)  ,\text{ \ \ \ \ }Z^{\prime}\left(  x_{1},0\right)  =0,\text{\ }%
x_{1}>l.
\end{align*}
According to the Kondratiev theory \cite{Ko}, see also \cite[Ch. 2]{NaPl}, a
solution of this problem admits the representation%
\begin{equation}
Z^{\prime}\left(  x\right)  =\left(  C^{0}+x_{1}C^{1}\right)  \sin\left(  \pi
x_{2}\right)  +\widetilde{Z}^{\prime}\left(  x\right)  \label{X10}%
\end{equation}
where $\widetilde{Z}^{\prime}\left(  x\right)  $ has the decay $O\left(
e^{-\sqrt{3}\pi x_{1}}\right)  ,$ is smooth everywhere in $\overline{\Pi
_{+}^{0}}$ except at the point $P=\left(  l,0\right)  $ and behaves as%
\begin{equation}
Z^{\prime}\left(  x\right)  =\pi\varphi+O\left(  r\right)  ,\text{
\ \ }r\rightarrow0, \label{X11}%
\end{equation}
while $\left(  r,\varphi\right)  \in\mathbb{R}_{+}\times\left(  0,\pi\right)
$ is the polar coordinates system centered at $P$. The singularity in
(\ref{X11}) leads the function $\widetilde{Z}^{\prime}$ out from the Sobolev
space $H^{1}\left(  \Pi_{+}^{0}\right)  $. Nevertheless, the solution
$Z^{\prime}$ still lives in appropriate Kondratiev space with a weighted norm
so that the coefficient $C^{1}$ in (\ref{X10}) can be computed by inserting
$Z^{\prime}\left(  x\right)  $ and $\sin\left(  \pi x_{2}\right)  $ into the
Green formula in $\Pi_{+}^{0}\left(  R\right)  $. To compensate for the
singularity, one may construct a boundary layer as a solution of the Dirichlet
problem in the unbounded domain (\ref{39}) in fig. \ref{f4}, a.

The above commentary exhibits all changes in the asymptotic analysis in
Section \ref{sect3.2}. As for the justification scheme in Section \ref{sect6},
it should be noted that, due to the Dirichlet condition (\ref{5}), the
inequality (\ref{hardy}) of Hardy's type takes the form%
\[
\left\Vert r^{-1}v^{\varepsilon};L^{2}\left(  \Pi_{+}^{\varepsilon}\left(
2l\right)  \right)  \right\Vert ^{2}\leq c\left\Vert v^{\varepsilon}%
;H^{1}\left(  \Pi_{+}^{\varepsilon}\left(  2l\right)  \right)  \right\Vert
^{2}%
\]
and sheds the factor $1+\left\vert \ln r\right\vert $ in the weight function
(\ref{J14}). As a result, the factor $\left(  1+\left\vert \ln\varepsilon
\right\vert \right)  ^{2}$ occurring in (\ref{8}) and (\ref{est}) for the
Neumann case, disappears from the asymptotic remainder in (\ref{6}) for the
Dirichlet condition.

\end{document}